\documentclass[a4paper, 10pt]{article}

\usepackage[utf8]{inputenc}
\usepackage{amsmath,amsfonts,amsthm,amssymb,mathrsfs,enumerate,enumitem,graphicx,float,hyperref,verbatim,comment}
\hypersetup{colorlinks=true,linkcolor=blue}   

\usepackage[letterpaper,top=2.5cm, bottom=2.5cm, left=2.5cm, right=2.5cm]{geometry}

\usepackage{pgf,tikz}
\usetikzlibrary{shapes,positioning,intersections,quotes,patterns}
\usetikzlibrary[topaths]
\usepackage{circuitikz}
\usepackage{color}
\usepackage{pgfplots}
\pgfplotsset{compat=1.7}
\usepackage{caption}
\usepackage{subcaption}

\def\R{{\mathbb R}} 
\def\N{{\mathbb N}}

\def\C{{\mathbb C}}

\def\S{{\mathbb{S}}}

\def\dist{{\operatorname{dist}}}  

\def\ds{\displaystyle}
\newcommand{\norm}[1]{\left\Vert#1\right\Vert}
\renewcommand{\leq}{\leqslant}                 
\renewcommand{\geq}{\geqslant}
\renewcommand{\tilde}{\widetilde}
\renewcommand{\hat}{\widehat}
\renewcommand{\bar}{\overline}
\renewcommand{\div}{\hbox{div\,}}

\renewcommand{\i}{\mathrm{\bf i}}                 
\renewcommand{\O}{\mathcal{O}}

\newtheorem{theorem}{Theorem}[section]
\newtheorem{remark}[theorem]{Remark}
\newtheorem{lemma}[theorem]{Lemma}

\newtheorem{corollary}[theorem]{Corollary}
\newtheorem{proposition}[theorem]{Proposition}

\numberwithin{equation}{section}

\def\keywords{
	\vspace{1ex}
	\noindent
	\if@twocolumn
	\small{\bf  Keywords}\/---$\!$    \else
	\begin{center}\small\ {\bf Keywords}\end{center}\quotation\small
	\fi}
\def\endkeywords{\vspace{0.6em}\par\if@twocolumn\else\endquotation\fi
	\normalsize\rm}

\title{Elliptic $L^p$ Carleman estimates near the boundary, quantitative unique continuation, and applications to controllability and stabilization\footnote{L. T. has received financial support from the Analysis Group of the Institut de Mathématiques de Bordeaux and the NumOpTes project (ANR-22-CE46-000). }}

\author{Kévin Le Balc'h\footnote{Sorbonne Université, CNRS, Université Paris Cité, Inria, Laboratoire Jacques-Louis Lions (LJLL), F-75005 Paris, France 
\texttt{kevin.le-balc-h@inria.fr}}
\and
Lotfi Thabouti\footnote{Université de Bordeaux, CNRS, Bordeaux INP, F-33400 Talence, France. \texttt{lotfi.thabouti@math.cnrs.fr}}
}  
\begin{document}
	
	\maketitle
	
	\begin{abstract}
    This article deals with quantitative unique continuation results up to and from the boundary for second-order elliptic equations with variable coefficients via $L^p$ Carleman estimates. The coefficients in the leading part are assumed to be Lipschitz continuous, and the potentials in the lower order term are allowed to belong to almost scale-invariant $L^q$ spaces. As applications, following the method of Lebeau and Robbiano, we obtain a spectral estimate which yields the interior null controllability of the associated parabolic equation, with an explicit dependence of the control cost on the potential. We also prove the logarithmic decay of the energy for damped wave equations with unbounded dampings and potentials. 
	\end{abstract}


	\setcounter{tocdepth}{1}
	\tableofcontents
	%
	%
	\section{Introduction}
   Let $d \geq 3$, $\Omega \subset \R^d$ be a bounded smooth domain. To a second order elliptic operator $-\div(A(x) \nabla \cdot)$, a potential $V$ and a source term $f$, we associate the elliptic equation
    \begin{equation}\label{Elliptic-equation}
        - \div (A(x) \nabla u) + V u = f, \qquad \text{in } \Omega.
    \end{equation}
    We say that problem \eqref{Elliptic-equation} with $f=0$ satisfies the unique continuation property in $\Omega$ if any solution $ u \in H^1(\Omega)$ that vanishes on a nonempty open subset of $\Omega$ vanishes identically in $\Omega$. For the operators considered here, this property is now well understood and we briefly review the existing results below. Its quantitative counterpart, which is the object of this paper, asserts that a solution of \eqref{Elliptic-equation} with small source term $f$, which is small on an open subset $ \omega \subset  \Omega$, or on a part $\Gamma \subset \partial \Omega$ of the boundary, remains small on a larger subset $ \O \subset \Omega$, provided it satisfies an a priori bound. This is expressed via an estimate involving a modulus of continuity $m$. From an interior observation on $\omega \subset \Omega$, it takes the form 
    \begin{equation}\tag{QUC$_{\omega}$} \label{QUC-omega}
        \|u\|_{\O} \leq m (\|u\|_{\omega}, \|f\|_{\Omega}, \|u\|_{\Omega}), 
    \end{equation} 
    while from a boundary observation on $\Gamma \subset \partial \Omega$, that is, from the Cauchy data of $u$ on $\Gamma$, it takes the form 
    \begin{equation}\tag{QUC$_{\Gamma}$}\label{QUC-Gamma}
        \|u\|_{\O} \leq m (\|u\|_{\Gamma }, \| \partial_{\nu} u\|_{\Gamma } , \|f\|_{\Omega}, \|u\|_{\Omega}), 
    \end{equation} 
    Here, the norms are taken in suitable function spaces, and $m$ tends to $0$ when the observation and source terms tend to $0$. We refer to \eqref{QUC-omega} and \eqref{QUC-Gamma} as quantitative unique continuation estimates from the interior and from the boundary, respectively. Such estimates have numerous applications, for instance to spectral estimates for sums of eigenfunctions, null-controllability of parabolic equations, stabilization of damped waves, and stability estimates for inverse problems.  

    The first unique continuation results for elliptic operator with non-analytic coefficients go back to Carleman \cite{carleman1939probleme} in dimension two, and were extended to higher dimensions, notably by Aronszajn and his collaborators \cite{aronszajn1962unique}. The key ingredient of their approach is a weighted $L^2$ estimate, now known as Carleman estimates. The simplest such estimate takes the form
    \begin{equation*}
    \tau^{\frac{3}{2}} \|e^{\tau \varphi} u \|_{L^2(\Omega)} + \tau ^{\frac{1}{2}} \|e^{\tau \varphi } \nabla u \|_{L^2(\Omega)} \leq  C \|e^{\tau \varphi} \div(A\nabla u) \|_{L^2(\Omega)}, \qquad \tau\gg 1, 
    \end{equation*}
   for smooth function $u$ compactly supported in $\Omega$, and for an appropriately chosen weight function $\varphi$. Such local $L^2$ Carleman estimates yield estimates of the form of \eqref{QUC-omega} with $\O \Subset \Omega$, but not up to the boundary. In the 1990s, Lebeau and Robbiano \cite{lebeau1995controle}, and Fursikov and Imanuvilov \cite{fursikov1996controllability}, developed $L^2$ Carleman estimates which remain valid up to the boundary of the domain. These allowed them to quantify the unique continuation properties up to the boundary, that is, \eqref{QUC-omega} with $\O$ reaching $ \partial \Omega$, and from the boundary, that is, \eqref{QUC-Gamma}, to derive controllability results for the heat equations. These techniques have since been extended to a wide range of parabolic and hyperbolic equations, and we refer the reader to the recent books \cite{leellipticV1, leellipticV2} for a more comprehensive presentation of $L^2$ Carleman estimates and their applications. Such estimates have been obtained for Lipschitz continuous leading coefficients $A$ and, more recently, for coefficients with jumps at interfaces \cite{LeRousseau2010carleman, LeRousseau2013interface, Vessella2017Jumps}, but they allow only mildly unbounded potentials $V$.
   %

   Establishing \eqref{QUC-Gamma} and \eqref{QUC-omega} becomes more difficult if one seeks to work with unbounded potentials $V$ in or near the scale invariant $L^{\frac{d}{2}}$ space. Indeed, $L^2$ Carleman estimates are insufficient and do not allow one to reach this critical integrability class, see for instance \cite{Saut1982-sur-l'unitcite}, where unique continuation is obtained for $V\in L^q$ with $q > \frac{2d}{3}$. For the unique continuation property, the key breakthrough is due to Jerison and Kenig \cite{jerison1985unique}, who replaced $L^2$ Carleman estimates by replaced by $L^p$ estimates of the form
   \begin{equation*}
       \|e^{\tau \varphi } u \|_{L^{\frac{2d}{d-2}}(\R^d)} \leq  C \|e^{\tau \varphi} \Delta u \|_{L^{\frac{2d}{d+2}}(\R^d)},  
   \end{equation*}
   for all $u \in C_c^\infty(\R^d\setminus \{0\})$, where $\varphi(x)=-\ln (|x|)$ and $C$ does not depend $\tau$, provided $\tau$ stays away from a discrete set. This result was extended to operators with smooth variable leading coefficients by Sogge \cite{sogge1989oscillatory,sogge1990strong}, while first-order terms in the scale-invariant $L^d$ space were handling by Wolff \cite{wolff1992property}, who introduced for this purpose an osculation argument. Of particular relevance to the present work, Koch and Tataru \cite{koch2001carleman} extended those $L^p$ Carleman estimates to Lipschitz continuous leading coefficients, with both zero and first-order terms, and derived strong unique continuation results. They also constructed counterexamples \cite{Koch-Tataru-02} showing that, on the scale of Lebesgue spaces, $L^{\frac{d}{2}}$  and $L^d$ are the sharp spaces  for the zero- and first-order terms, respectively; see also the work of Dos Santos Ferreira \cite{ferreira2003sharp} for $L^p$ Carleman estimates for second-order operators of real principal type. All of these results are based on local $L^p$ Carleman estimates, that is, estimates for compactly supported functions, which cannot be used near the boundary. In particular, they can be applied to derive estimates of the form of \eqref{QUC-omega} with $\O \Subset \Omega$ with $L^\frac{d}{2}$ potential, but not up to the boundary. 

   Global $L^p$ Carleman estimates valid up to the boundary were obtained by Dehman, Ervedoza, and the second author \cite{Dehman-Ervedoza-Thabouti-2023}, for the Laplace operator in a smooth bounded domain of $\R^d$, $d\geq 3$, and for weight functions satisfying Hörmander's sub-ellipticity condition and having $\partial \Omega$ as a level set. As an application, they obtained estimates of the form \eqref{QUC-omega} up to the boundary for solutions of  \eqref{Elliptic-equation} with $A=\operatorname{Id}$, $V\in L^{q_0}(\Omega)$ with $q_0 > \frac{d}{2} $, and $f=W_1 \nabla u + \div(W_2 u)$ with $W_1 \in L^{q_1}(\Omega; \C^d)$ and $W_2 \in L^{q_2}(\Omega; \C^d)$, $q_1, q_2 > \frac{3d-2}{2}$ (with an additional restriction when both $W_1$ and $W_2$ are nonzero), with an explicit dependence of the modulus of continuity $m$ on the norms of $V$, $W_1$ and $W_2$. Since these estimates involve an interior observation and require $\partial\Omega$ to be a level set of the weight, therefore, they cannot be applied to derive \eqref{QUC-Gamma}. Caro, Ervedoza and the second author \cite{caro2024quantitative} then combined $L^p$ Carleman estimates with Wolff's argument to treat first order terms $W_1 \in L^{q_1}(\Omega; \C^d)$  and $W_2 \in L^{q_2}(\Omega; \C^d)$, $q_1, q_2 >d$, that is, arbitrarily near to the scale-invariant space $L^d$, but only for local estimates of the form \eqref{QUC-omega}, with $\O \Subset \Omega$.  
   \medskip

   Our aim is to develop $L^p$ Carleman estimates near the boundary and to derive from them quantitative unique continuation estimates of the form \eqref{QUC-Gamma}, from the boundary, and \eqref{QUC-omega}, from the interior up to the boundary, for solutions of \eqref{Elliptic-equation} with $A$ Lipschitz continuous and $V \in L^q$, $q>\frac{d}{2}$. As applications, we obtain a spectral inequality and the null-controllability of the associated parabolic equation for real valued potentials $V \in L^{q}(\Omega)$ with $q > \frac{d+1}{2}$, with an explicit dependence of the control cost on $\|V\|_{L^q(\Omega)}$, as well as the logarithmic decay of the energy of damped wave equations with unbounded dampings in $L^p(\Omega)$, $p>d$, and unbounded potentials in $L^q(\Omega),$ $q>\frac{d+1}{2}$. 

   \section{Main results and discussion}
   Let $A(x) = (a^{ij}(x))_{1 \leq i,j \leq d}$ be a symmetric matrix-valued function satisfying, for given constants $\Lambda_0 \in (0,1)$ and $\Lambda _1 >0$,  the unifom ellipticity condition 
\begin{equation}
\label{Ellipticity-A}
\Lambda_0 |\xi|^2 \leq a^{ij}(x) \xi_i \xi_j \leq \Lambda_0^{-1} |\xi|^2, \qquad \forall x \in \overline{\Omega}, \forall \xi \in \R^d, 
 \end{equation}
 and the Lipschitz condition 
 \begin{equation}\label{Lipschitz-condition-A}
 |A(x)- A(y)| \leq \Lambda_1 |x-y|, \qquad \forall x,y \in \Omega.
 \end{equation}
\subsection{Carleman estimates}
Our first main result concerns local $L^p$ Carleman estimates near the boundary for variable coefficients.
\begin{theorem} \label{Thm-Carleman-General}
Let $d\geq 3$, $r_0 \in (0, \tfrac{1}{2}]$, $\Omega=B^+_{2r_0}:= B(0,2r_0) \cap \{x_1>0\}$ and $\Gamma=B_0(2r_0) \cap \{x_1>0\}$. Let $A$ satisfy \eqref{Ellipticity-A}--\eqref{Lipschitz-condition-A}, let $\varphi \in C^{2}([0,2r_0])$ satisfy 
\begin{equation}
    \label{Weight-function-varphi}
    \varphi(0) = 0, \qquad 0 < \alpha_1 \leq - \varphi'(x_1) \leq \alpha_2, \qquad \varphi''(x_1) \geq \beta >0 \qquad \text{on } [0,2r_0], 
\end{equation}
and set 
\begin{equation} \label{Weight-function-Phi}
    \Phi(x) = \varphi(x_1) - \tfrac{\varepsilon}{2} \, |x'|^2. 
\end{equation}

There exist $C>0$, $\varepsilon_2 \in (0, 1]$ and $\tau_0 \geq 1$, depending only on $d, \Lambda_0, \Lambda_1, \alpha_1, \alpha_2, \beta$ and $r_0$, such that for all $0 < \varepsilon \leq \varepsilon_2$, all $\tau \geq  \tau_0$, and all $u \in \overline{C_c^{\infty}}^{\, \Omega}(B_0(r_0))$ satisfying 
\begin{equation}\label{Elliptic-Eq-general}
-\div(A (\varepsilon x ) \nabla u ) = f_2 + f_{2_*'}+ \div(F) \qquad \text{in } \Omega, 
\end{equation}
with $f_2 \in L^2(\Omega)$, $f_{2_*'} \in L^{\frac{2d}{d+2}}(\Omega)$ and $F\in L^2(\Omega; \C^d))$, the following estimates hold
\begin{multline} \label{Carleman-General-1}
    \tau^{\frac{3}{2}} \|e^{\tau \Phi} u\|_{L^2(\Omega)}+ \tau^{\frac{1}{2}} \|e^{\tau \Phi} \nabla u\|_{L^2(\Omega)} \leq C\left(   \|e^{\tau \Phi} f_2\|_{L^2(\Omega)}+ \tau^{\frac{3}{4}-\frac{1}{2d}} \|e^{\tau \Phi} f_{2_*'}\|_{L^\frac{2d}{d+2}(\Omega)}+ \tau \|e^{\tau \Phi} F\|_{L^2(\Omega)}
    \right. \\ \left. 
    + \tau^{\frac{5}{4}}\|e^{\tau \Phi} u\|_{L^2(\Gamma)} +\tau^{\frac{3}{4}} \|e^{\tau \Phi} u\|_{H^{\frac{1}{2}}(\Gamma)}+\tau^{\frac{1}{4}- \theta} \|e^{\tau \Phi} (A(\varepsilon\cdot) \nabla u-F)\cdot \nu \|_{H^\theta(\Gamma)},
    \right) 
\end{multline}
and 
\begin{multline} \label{Carleman-General-2}
    \tau^{\frac{3}{4}+\frac{1}{2d}} \|e^{\tau \Phi} u\|_{L^\frac{2d}{d-2}(\Omega)} \leq C\left(   \|e^{\tau \Phi} f_2\|_{L^2(\Omega)}+ \tau^{\frac{3}{4}+\frac{1}{2d}} \|e^{\tau \Phi} f_{2_*'}\|_{L^\frac{2d}{d+2}(\Omega)}+ \tau \|e^{\tau \Phi} F\|_{L^2(\Omega)} 
    \right. \\ \left. 
   + \tau^{\frac{5}{4}}\|e^{\tau \Phi} u\|_{L^2(\Gamma)} +\tau^{\frac{3}{4}+\frac{1}{2d}} \|e^{\tau \Phi} u\|_{H^{\frac{1}{2}}(\Gamma)}+\tau^{\frac{1}{4}+\frac{1}{2d}- \theta} \|e^{\tau \Phi} (A(\varepsilon\cdot) \nabla u-F)\cdot \nu \|_{H^\theta(\Gamma)},
    \right) 
\end{multline}
where $\nu = e_1$ and $\theta \in \{0,-1/2\}$. 
\end{theorem}
The $L^p$ Carleman estimates \eqref{Carleman-General-1}--\eqref{Carleman-General-2} are boundary counterpart of the estimates of \cite[Theorem 1.1]{Dehman-Ervedoza-Thabouti-2023} for the Laplace operator: the interior observation terms are replaced by boundary observation terms on $\Gamma$, involving the Dirichlet trace terms in $L^2(\Gamma)$ and $H^{\frac{1}{2}}(\Gamma)$ and the Neumann trace in $H^{\theta}(\Gamma)$. In contrast to \cite[Theorem 1.1]{Dehman-Ervedoza-Thabouti-2023}, which holds with an arbitrary smooth bounded domain of $\R^d$, $d\geq3$, and general weights satisfying H\"ormander's sub-ellipticity conditions and having $\partial \Omega$ as a level set, we work here in a half-ball with the explicit weight \eqref{Weight-function-Phi}. This local setting is sufficient to prove the quantitative unique continuation from and up to the boundaries, which cannot be achieved with weights having $\partial \Omega$ a level set. The level set condition was imposed in \cite{Dehman-Ervedoza-Thabouti-2023} in order to perform a change of variables reducing the general case to a strip. Here, on the contrary, the weight on $\{x_1=0\}$ is $\Phi(0,x')=-\frac{\varepsilon}{2}|x'|^2$, which provides a convexification with respect to $\{x_1=0\}$.

A further contribution, in this $L^p$ setting with non-vanishing traces on the boundary, is the extension to variable coefficient operators. Our Carleman estimates are stated for the rescaled operator $\operatorname{div}\big(A(\varepsilon\cdot)\nabla \cdot\big)$ with $A$ Lipschitz continuous and $\varepsilon$ small enough, which is sufficient to prove quantitative unique continuation results for operators $-\div(A \nabla \cdot)$ with Lipschitz continuous $A$. Indeed, once $\varepsilon$ is fixed sufficiently small, the change of variable $u_\varepsilon(x) = u(\varepsilon x)$ gives the corresponding estimates for $\div(A\nabla \cdot)$, in a half-ball $B^+_{2\varepsilon r_0}$, for functions supported in $B_0(\varepsilon x)$, with the rescaled weight $\Phi(\cdot/\varepsilon)$. These estimates have the same powers of $\tau$, and constants depending on $\varepsilon$;  see the proof of Lemma \eqref{lemma-local-at-the-boundary}.

When $f_{2_*'}=0$ in \eqref{Carleman-General-1}, estimates of the same type, with source terms $f_2 + \div(F)$ of weaker regularity and on a smooth compact Riemannian manifold with boundary, are proved by patching microlocal estimates in the book of Le Rousseau, Lebeau, and Robbiano \cite[Theorem 14.35]{leellipticV2}, where the optimality of the factor $\tau^{3/4}$ in front of the boundary terms is also established \cite[Section~14.5]{leellipticV2}. Related global estimates with a $C^2$ coefficients, an interior observation and a nonhomogeneous Dirichlet term, were obtained earlier by Imanuvilov and Puel in \cite{imanuvilov2003global}. Note that, even in this Hilbertian case, our estimates allow Lipschitz continuous coefficients, whereas \cite{leellipticV2} and \cite{imanuvilov2003global} assume smooth and $C^2$ coefficients, respectively. The powers of $\tau$ in front of the boundary terms in \eqref{Carleman-General-1} are the same as  in \cite[Theorem 14.35]{leellipticV2}; they are therefore optimal in the Hilbertian framework, and we expect them to be optimal in the $L^p$ framework as well. 

The term $\tau^{\frac{5}{4}}\|e^{\tau\Phi}u\|_{L^2(\Gamma)}$ in \eqref{Carleman-General-2} is not produced by the $L^{\frac{2d}{d-2}}$ estimates. The corresponding estimate of Proposition \ref{Prop-Carleman-Strip-2} contains no $L^2$ trace term. It comes from the $L^2$ estimates of Proposition \ref{Prop-Carleman-Strip-2}, which are essential to handle perturbations. 
\medskip

The proof of \eqref{Carleman-General-1}--\eqref{Carleman-General-2}
also yields another two Carleman estimates, stated in the following theorem, which corresponds to the reverse sign condition $\varphi'\geq\alpha_1>0$, that is, to $\partial_\nu\Phi<0$ on $\Gamma$. These estimates allows us to propagate smallness from the interior up to the boundary, see Theorem \ref{Est-local-up-to-boundary}. The proof of this Theorem is given in Section \ref{Sec-proof-Thm2}.  

\begin{theorem} \label{Thm2-Carleman-General}
Let $d\geq 3$, $r_0 \in (0, \tfrac{1}{2}]$, $\Omega=B^+_{2r_0}:= B(0,2r_0) \cap \{x_1>0\}$ and $\Gamma=B_0(2r_0) \cap \{x_1>0\}$. Let $A$ satisfy \eqref{Ellipticity-A}--\eqref{Lipschitz-condition-A}, let $\varphi \in C^{2}([0,2r_0])$ satisfy 
\begin{equation}
    \label{Weight-function-varphi-Thm2}
    \varphi(0) = 0, \qquad 0 < \alpha_1 \leq \varphi'(x_1) \leq \alpha_2, \qquad \varphi''(x_1) \geq \beta >0 \qquad \text{on } [0,2r_0], 
\end{equation}
and set 
\begin{equation} \label{Weight-function-Phi-Thm2}
    \Phi(x) = \varphi(x_1) - \tfrac{\varepsilon}{2} \, |x'|^2. 
\end{equation}

There exist $C>0$, $\varepsilon_2 \in (0, 1]$ and $\tau_0 \geq 1$, depending only on $d, \Lambda_0, \Lambda_1, \alpha_1, \alpha_2, \beta$ and $r_0$, such that for all $0 < \varepsilon \leq \varepsilon_2$, all $\tau \geq  \tau_0$, and all $u \in \overline{C_c^{\infty}}^{\, \Omega}(B_0(r_0))$ satisfying \eqref{Elliptic-Eq-general} with $f_2 \in L^2(\Omega)$, $f_{2_*'} \in L^{\frac{2d}{d+2}}(\Omega)$ and $F\in L^2(\Omega; \C^d))$, the following estimates hold
\begin{multline} \label{Carleman-General-1-Thm2}
    \tau^{\frac{3}{2}} \|e^{\tau \Phi} u\|_{L^2(\Omega)}+ \tau^{\frac{1}{2}} \|e^{\tau \Phi} \nabla u\|_{L^2(\Omega)} \leq C\left(   \|e^{\tau \Phi} f_2\|_{L^2(\Omega)}+ \tau^{\frac{3}{4}-\frac{1}{2d}} \|e^{\tau \Phi} f_{2_*'}\|_{L^\frac{2d}{d+2}(\Omega)} 
    \right. \\ \left. 
   + \tau \|e^{\tau \Phi} F\|_{L^2(\Omega)} +\tau^{\frac{3}{4}} \|e^{\tau \Phi} u\|_{H^{\frac{1}{2}}(\Gamma)}
    \right), 
\end{multline}
and 
\begin{multline} \label{Carleman-General-2-Thm2}
    \tau^{\frac{3}{4}+\frac{1}{2d}} \|e^{\tau \Phi} u\|_{L^\frac{2d}{d-2}(\Omega)} \leq C\left(   \|e^{\tau \Phi} f_2\|_{L^2(\Omega)}+ \tau^{\frac{3}{4}+\frac{1}{2d}} \|e^{\tau \Phi} f_{2_*'}\|_{L^\frac{2d}{d+2}(\Omega)} + \tau \|e^{\tau \Phi} F\|_{L^2(\Omega)}
    \right. \\ \left. 
    +\tau^{\frac54}\|e^{\tau\Phi}u\|_{L^2(\Gamma)} +\tau^{\frac{3}{4}+\frac{1}{2d}} \|e^{\tau \Phi} u\|_{H^{\frac{1}{2}}(\Gamma)}
    \right). 
\end{multline}
\end{theorem}
\subsection{Interpolation inequalities} 
A consequence of the Carleman estimates of Theorem \eqref{Thm-Carleman-General} is the following interpolation inequality from the boundary, which is an estimate of the form \eqref{QUC-Gamma} with a Hölder-type modulus of continuity.  
\begin{theorem}[Interpolation inequality from the boundary]\label{Thm-interpolation-inequality}
 Let $d\geq 3$, $\Omega \subset \R^d$ be a bounded domain of class $C^{1,1}$, $V \in L^q(\Omega)$ with $q>\frac{d}{2}$, and $A$ satisfying \eqref{Ellipticity-A}--\eqref{Lipschitz-condition-A}. Let $\mathcal{U}$ be an open subset of $\Omega$ and let $\mathcal U _1$ be an open neighborhood of $\overline{\mathcal U}$ in $\R^d$. Let $y \in \partial \Omega$ and let $\mathcal{V}$ be an open neighborhood of $y$ in $\R^d$. 

Then there exist $C>0$ and $\alpha \in (0,1)$, depending only on $d, \Lambda_0, \Lambda_1, q, \Omega, \mathcal{U}, \mathcal{U}_1, $ and $\mathcal{V}$, such that
\begin{equation}\label{interpolation-inequality}
    \|u\|_{H^1(\mathcal{U})} \leq  C e^{
	  C \norm{V}_{L^{q}(\Omega)}^{\gamma_d(q)}} \|u\|_{H^1(\Omega)}^{1-\alpha} \left(
\|f\|_{L^2(\Omega)}
+ \|u|_{\partial\Omega}\|_{H^\frac{1}{2}(\mathcal{V} \cap \partial\Omega)}
+ \| (A\nabla u )\cdot \nu\|_{H^{-\frac{1}{2}}(\mathcal{V} \cap \partial\Omega)}
\right)^{\alpha},
\end{equation}
for all $u \in H^1(\Omega)$ such that $u=0$ on $\partial \Omega \cap {\mathcal{U}}_1 $ and solving \eqref{Elliptic-equation} with $f\in L^2(\Omega)$, where 
\begin{equation} \label{Conditions-gamma(q)}
	\gamma_d(q)= \left\{
		\begin{array}{lll}
			\ds \frac{1}{\ds \frac{3}{2} \left(1- \frac{d}{2q}\right) + \frac{1}{2q} } & \text{ if } q \geq d,  
			 \\
			\ds \frac{1}{ \ds \left(\frac{3}{4} + \frac{1}{2d}\right)  \left(2- \frac{d}{q}\right) } & \text{ if } \frac{d}{2} < q \leq d.
		\end{array}
		\right.
\end{equation}
\end{theorem}

When $V=0$, the boundary interpolation inequality \eqref{interpolation-inequality} was first proved by Lebeau and Robbiano in \cite{lebeau1995controle}, Le Rousseau and Robbiano \cite{LeRousseau2010carleman} for an isotropic diffusion coefficient with a jump at an interface, and Le Rousseau and Lerner \cite{LeRousseau2013interface} for a jump in the whole diffusion matrix.  Here, for Lipschitz continuous $A$, a potential $V\in L^q(\Omega)$ with $q>\frac{d}{2}$ is allowed, with an explicit dependence of the constant on $\|V\|_{L^q(\Omega)}$, as in \cite{Dehman-Ervedoza-Thabouti-2023, caro2024quantitative}. The proof follows the propagation scheme of \cite[Theorem 5.13]{leellipticV1}, the new ingredient here is the $L^p$ Carleman estimate near the boundary of Theorem \ref{Thm-Carleman-General}.

Moreover, the proof of Theorem~\ref{Thm-interpolation-inequality} can be adapted to the case where first-order terms and a divergence source term are added to \eqref{Elliptic-equation}, that is, to solutions of
\begin{equation*}
-\div(A\nabla u)+Vu=f+W_1\cdot\nabla u+\div(W_2u)+\div F\qquad\text{in }\Omega,
\end{equation*}
with $F\in L^2(\Omega;\C^d)$, $W_1\in L^{q_1}(\Omega;\C^d)$ and $W_2\in L^{q_2}(\Omega;\C^d)$, where $q_1,q_2>\frac{3d-2}2$ and, if $W_1$ and $W_2$ are both nonzero, $\frac1{q_1}+\frac1{q_2}<\frac{4(d-1)}{d(3d-2)}$, as in~\cite[Theorem~1.2]{Dehman-Ervedoza-Thabouti-2023}. For the sake of simplicity, we do not consider these terms here. 

\medskip 

The Carleman estimates Theorem \ref{Thm2-Carleman-General} yields the following interpolation initiated from an interior observation, which is an estimate of the form \eqref{QUC-omega} up to the boundary. 

\begin{theorem}[Interpolation inequality from the interior, up to the boundary]\label{Thm-interpolation-from-interior}
Let $d\geq 3$, $\Omega \subset \R^d$ be a bounded domain of class $C^{1,1}$, $V \in L^q(\Omega)$ with $q>\frac{d}{2}$, and $A$ satisfying \eqref{Ellipticity-A}--\eqref{Lipschitz-condition-A}. Let $\omega$ and $\mathcal{U}$ be nonempty open subsets of $\Omega$ and let $\mathcal U _1$ be an open neighborhood of $\overline{\mathcal U}$ in $\R^d$. 

Then there exist $C>0$ and $\alpha \in (0,1)$, depending only on $d, \Lambda_0, \Lambda_1, q, \Omega, \mathcal{U}, \mathcal{U}_1, $ and $\omega$, such that
\begin{equation}\label{Est-local-up-to-boundary}
    \|u\|_{H^1(\mathcal{U})} \leq  C e^{
	  C \norm{V}_{L^{q}(\Omega)}^{\gamma_d(q)}} \|u\|_{H^1(\Omega)}^{1-\alpha} \left(
\|f\|_{L^2(\Omega)}
+ \|u \|_{L^2(\omega)}\right)^{\alpha},
\end{equation}
for all $u \in H^1(\Omega)$ such that $u=0$ on $\partial \Omega \cap {\mathcal{U}}_1 $ and solving \eqref{Elliptic-equation} with $f\in L^2(\Omega)$, where $\gamma_d$ is given by \eqref{Conditions-gamma(q)}. 
\end{theorem}

Theorem \ref{Thm-interpolation-from-interior} is an analogue of \cite[Theorem 5.11]{leellipticV1} for Lipschitz continuous coefficients and potentials $V \in L^q$, $q >d/2$. For $A=\operatorname{Id}$, estimates of the form \eqref{QUC-omega} up to the boundary, with an explicit dependence on $\|V\|_{L^q(\Omega)}$, were obtained in \cite[Theorem 1.2]{Dehman-Ervedoza-Thabouti-2023}.

Unlike Theorem \ref{Thm-interpolation-inequality}, the observation in Theorem \ref{Thm-interpolation-from-interior} is purely interior and involves only the $L^2$-norm of $u$ on $\omega$, while the set $\mathcal{U}$ may reach the boundary. Its proof combines the Carleman estimate of Theorem \ref{Thm2-Carleman-General} near the boundary with a propagation of smallness argument from $\omega$ to $\mathcal{U}$. Theorem \ref{Thm-interpolation-from-interior} is also an ingredient of the proof of Theorem \ref{Thm-interpolation-inequality}: the Carleman estimate of Theorem \ref{Thm-Carleman-General} gives the smallness of $u$ near the observed part $\mathcal V\cap\partial\Omega$ of the boundary, and Theorem \ref{Thm-interpolation-from-interior} propagates it to $\mathcal U$. 
\bigskip

The method of Lebeau and Robbiano \cite{lebeau1995controle}, which we use below, requires an interpolation inequality for the elliptic operator 
\begin{equation*}
    \operatorname{Q}: = - \partial_s ^2 - \div_x(A(x) \nabla_x \cdot) + V(x) \qquad \text{in } Z := (0,S_0) \times \Omega, 
\end{equation*}
in the subdomain $Y:= (\beta, S_0-\beta) $, where $S_0 >0$ and $\beta \in (0,S_0/2)$, with an observation $\{0\}\times \omega\subset\partial Z$. This inequality follows from Thereom \ref{Thm-interpolation-inequality} in dimension $d+1$, with $\Omega$ replaced by $Z$, after rounding its edges away from $\overline Y$ and $\{0\}\times\overline{\omega_0}$, $\mathcal U$ by $Y$, $\mathcal U_1$ by $(\beta/2,S_0-\beta/2)\times\R^d$, and $\mathcal V$ by $(-\beta,\beta)\times\omega_0$, where $\omega_0\Subset\omega$. Note that $\mathcal V\cap\partial Z=\{0\}\times\omega_0$, and that the normal derivative on this set is $-\partial_s$.  This gives the following result.
\begin{corollary}\label{corollary-interpolations-Q}
Let $d\geq 3$, $\Omega$ and $A$ be as in Theorem \ref{Thm-interpolation-inequality}, let $V \in L^q(\Omega)$ with $q>\frac{d+1}{2}$, and let $\omega$ and $\mathcal{U}$ be a nonempty open subset of $\Omega$. 

Then there exist $C>0$ and $\alpha \in (0,1)$, depending only on $d, \Lambda_0, \Lambda_1, q, \Omega, S_0, \beta $ and $\omega$, such that
\begin{equation}\label{interpolations-Q}
    \|u\|_{H^1(Y)} \leq  C e^{
	  C \norm{V}_{L^{q}(\Omega)}^{\gamma_{d+1}(q)}} \|u\|_{H^1(Z)}^{1-\alpha} \left(
\|\operatorname{Q}u\|_{L^2(\Omega)} +\|u|_{s=0} \|_{H^1(\omega)}
+ \| \partial_s u|_{s=0} \|_{L^2(\omega)}\right)^{\alpha},
\end{equation}
for all $u \in H^2(Z)$ such that $u=0$ on $(0,S_0) \times \partial \Omega$, where $\gamma_{d+1}$ is defined by \eqref{Conditions-gamma(q)} with $d$ replaced by $d+1$. 
\end{corollary}
\bigskip

Several possible extensions are in order.  
\begin{itemize}
    \item \textit{Dependence on the potential.} Our interpolation inequalities hold for $V\in L^{q}(\Omega)$, $q>d/2$, with constants depending explicitly on the potential through $\exp \left(C \|V\|_{L^q(\Omega)}^{\gamma_d(q)}\right)$. This dependence is closely related to quantitative strong unique continuation for Schröding operators. For solutions of $\Delta u + V u=0$ with $V\in L^q$, Davey \cite{davey2020quantitative} prove that the maximal order of vanishing is bounded by $C(1+ \|V\|_L^q)^\mu$, with $\mu = \gamma_d(q)$ for $q >d$ and $\mu = \gamma_d(q) + \varepsilon$ for $\frac{d}{2} < q \leq d$. For $q=\infty$, Meshkov's construction \cite{Mes91} shows that the exponent $\gamma_d(\infty) =\frac{2}{3}$ is sharp for complex-valued potentials; see also the optimality results of Duyckaerts, Zhang and Zuazua \cite{duyckaerts2008optimality} for observability inequalities. To our knowledge, the optimality of $\gamma_d(q)$ for finite $q$ is open.  
\item \textit{Coefficients with jumps.} We expect our results to extend to anisotropic diffusion matrices that are Lipschitz continuous on both sides of an interface and jump across it. In the $L^2$ framework, Carleman estimates for such operators were proved by Le Rousseau and Lerner \cite{LeRousseau2013interface}, by a microlocal factorization of the conjugated operator, and by Di Cristo, Francini, Lin, Vessella and Wang \cite{Vessella2017Jumps}, who replaced microlocal analysis by a Fourier transform in the tangential variables and an explicit factorization of the conjugated operator, which are also the main tools of our approach. Combining these techniques with ours should give $L^p$ Carleman estimates for this class of operators, hence quantitative unique continuation for the potentials $V \in L^q$, $q> \frac{d}{2}$, and the corresponding three-region inequalities, which are the key tool for estimating the size of an inclusion from a single boundary measurement \cite{Vessella2016Three-region}. 

\item \textit{First-order terms.} First-order terms $W_1\cdot\nabla y+\div(W_2y)$, with $W_j\in L^{q_j}$ and $q_j\in(d,2d]$, $j \in \{1,2\}$ could probably be handled as well, by combining a refined version of Theorem \ref{Thm-Carleman-General} with the argument of Wolff \cite{wolff1992property}. This would, however, require adapting that argument near the boundary, which is precisely the difficulty that led Caro, Ervedoza and the second author, in \cite{caro2024quantitative}, to restrict themselves to interior observation of the form \eqref{QUC-omega} with $\O \Subset \Omega$.  
\end{itemize}

\subsection{Spectral inequality and null-controllability} 
Let $V \in L^q(\Omega; \R)$, $q> \frac{d+1}{2}$, and let $V_-=\max\{-V, 0\}$ denote its negative part. Since $q > \frac{d}{2}$, Hölder's and Sobolev's inequalities give, for all $u \in H^1_0(\Omega)$,
\begin{equation}\label{quantitative-lower-bound-P}
    \int_\Omega V_- |u|^2 d x \leq \frac{1}{2} \int_\Omega A \nabla u \cdot \nabla \bar u d x + R_V \|u\|^2_{L^2(\Omega)}, \quad R_V:=C(1+ \|V_-\|_{L^q}^{\frac{2q}{2q-d}}),
\end{equation}
where $C$ depends only on $d$, $q$, $\Omega$ and $\Lambda_0$, and, in the same way, $\int_\Omega |V| |u|^2 d x< \infty $. Hence the quadratic form $u \mapsto \int_\Omega A \nabla u \cdot \nabla \bar u d x  + \int_\Omega |V| |u|^2 d x $ is closed on $H^1_0(\Omega)$ and bounded from below by - $R_V\|u\|^2_{L^2(\Omega)}$. The associated operator $\operatorname{P}:=-\div(A\nabla \cdot) +V$ with domain
\begin{equation*}
D(\operatorname{P)}= \{ u \in H^1_0(\Omega); \quad -\div(A\nabla u) +Vu \in L^2(\Omega) \},
\end{equation*}
 is self-adjoint on $L^2$ and satisfies $\operatorname{P}\geq - R_V $ in the sense of quadratic forms. Since its form domain $H^1_0(\Omega)$ is compactly embedded in $L^2(\Omega)$, $\operatorname{P}$ has a compact resolvent: there exist an orthonormal basis $(\varphi_j)_{j\in \N}$ of $L^2(\Omega)$ of real valued eigenfunctions and a non-decreasing sequence of eigenvalues $(\mu_j)_{j \in \N},$ with $\mu_j \geq -R_V$ and  $\mu_j \to  +\infty$, such that 
\begin{equation}
    \label{eigenfunction}
    \operatorname{P} \varphi_j= \mu_j \varphi_j. 
\end{equation}
\begin{theorem}[Spectral inequality] \label{thm-spectral-inquality}
    Let $d\geq 3$, $\Omega$ and $A$ be as in Theorem \ref{Thm-interpolation-inequality}, let $V \in L^q(\Omega; \R)$ with $q>\frac{d+1}{2}$, and let $\omega$ and $\mathcal{U}$ be a nonempty open subset of $\Omega$. There exists $C>0$, independent of $V$, such that, for all $\mu \geq 0$, 
\begin{equation}\label{Spectral-inequality}
    \|w\|_{L^2(\Omega)} \leq C e^{C \left(\|V\|_{L^q(\Omega)}^{\gamma_{d+1}(q)} + \sqrt{\mu}\right) } \|w\|_{L^2(\omega)}, \quad \text{for all } w = \sum_{ \mu_j \leq \mu } w_j \varphi_j. 
\end{equation}
\end{theorem}
By the method of Lebeau and Robbiano \cite{lebeau1995controle}, Theorem \ref{thm-spectral-inquality} yields the following null-controllability result, in any time $T>0$ : 

\begin{theorem}[Null-controllability]\label{thm-null-controllability} Under the assumptions of Theorem \ref{thm-spectral-inquality}, there exists $C>0$, depending only $\omega, \Omega, \Lambda_0, \Lambda_1,$ and $ q $ such that for all $T>0$ and $y_0\in L^2(\Omega)$, there exists a control $v\in L^2((0,T)\times \omega)$ satisfying 
\begin{equation}\label{control-cost}
    \|v\|_{L^2(0,T;L^2(\omega))}
\leq Ce^{ C\left(1+\frac1T
+\|V\|_{L^q(\Omega)}^{\gamma_{d+1}(q)}+TR_V\right) 
}\|y_0\|_{L^2(\Omega)}
\end{equation}
for which the solution $y\in C([0,T]; L^2(\Omega))$ of 
\begin{equation}\label{Heat-equation}
\begin{cases}
\partial_t y - \div (A \nabla  y) + V y = v \mathbf{1}_{\omega}  & \text{in } (0,T)\times\Omega, \\[0.3em]
y = 0 & \text{on } (0,T)\times\partial\Omega, \\[0.3em]
y(0,\cdot) = y_0(x) & \text{in } \Omega,
\end{cases}
\end{equation}
satisfies $y(T,\cdot) = 0.$
\end{theorem}

The null-controllability of the heat equation with unbounded, possibly time-dependent, potentials has mainly been obtained through global $L^2$ Carleman estimates.  Imanuvilov and Yamamoto \cite{ImanuvilovYamamoto2003} proved it for $V\in L^\infty(0,T; L^q(\Omega))$ with $q>\frac{2d}{3}$, without explicit dependence on $V$, while Duyckaerts, Zhang and Zuazua \cite[Theorem 2.1]{duyckaerts2008optimality}, obtained for $q\geq d$, an observation constant of the order of $\exp(C[1+\frac{1}{T}+ T \|V\|_{L^q(\Omega)}+ \|V\|_{L^q(\Omega)}^{\frac{1}{3/2- d/q}} ])$; see also \cite{PhungWang2013} for observations on measurable sets in time, with $q>d$. For time-independent potentials, Theorem \ref{thm-null-controllability} relaxes the condition $q>\frac{2d}{3}$ to $q>\frac{d+1}{2}$, with an explicit dependence on $\|V\|_{L^q(\Omega)}$. Moreover, for $d\leq q<\infty$, the exponent $\gamma_{d+1}(q)$ in the $T$-independent factor of \eqref{control-cost} is smaller than $\frac{1}{3/2-d/q}$, both exponents being equal to $\frac{2}{3}$ when $q=\infty$. We do not claim the optimality of $\gamma_{d+1}(q)$: for $q<\infty$, this question is open, already for the exponent of \cite{duyckaerts2008optimality}; see \cite[Open problems: 8.4]{duyckaerts2008optimality}. In terms of local singularities, our results allow $V(x)=c|x-x_0|^{-a}$ with $c\in\R$ and $a<2d/(d+1)$, instead of $a< 3/2$ in \cite{ImanuvilovYamamoto2003}.

\paragraph{Divergence-form reduction and observation from sets of positive measure.} Theorem \ref{thm-spectral-inquality} requires the observation set $\omega$ to be open. Without analyticity assumptions, spectral inequalities with observation on sets of positive measure were established by Burq and Moyano \cite{burq2022propagation} for the divergence-form operators with Lipschitz coefficients and $V=0$, using the propagation-of-smallness results of Logunov and Malinnikova \cite{LogunovMalinnikova2018}; see also \cite{apraiz2014observability} and the references therein for the analytic case. For Schr\"odinger operators with real-valued potentials, such inequalities were then established by Le Balc'h and Martin \cite{LeBalcMartin2025} and by Zhu \cite{ZhuVariousPotentials2026}, with applications to the null-controllability of parabolic equations from sets of positive measure. Both works lift the spectral projector to an elliptic equation in one additional variable, as in the proof of Theorem \ref{thm-spectral-inquality}, and then remove the potential through a positive multiplier. More precisely, if $\phi>0$ solves $-\div(A\nabla \phi) + V \phi=0$, then $u=\phi v$ satisfies $-\div(A\nabla u)+Vu=0$ if and only if
\begin{equation*}
-\div(\phi^2A\nabla v)=0.
\end{equation*}
When $V\in L^\infty$, one may assume $V\geq0$ by shifting the spectrum. The multiplier $\phi$ can then be chosen to be Lipschitz continuous and bounded from below, and the results of \cite{LogunovMalinnikova2018} apply to $v$. 

This argument seems to extend to real-valued potentials $V\in L^q$ with $q>d$, by constructing positive multipliers $\phi= \phi(x)$ on sufficiently small balls. Indeed, elliptic regularity yields $\phi\in W^{2,q}_{\mathrm{loc}}$, and hence $W^{1,\infty}_{\mathrm{loc}}$. For $q \leq d$, however, this regularity does not in general guarantee that $\phi$ is Lipschitz continuous. Consequently, the transformed coefficient matrix $\phi^2 A$ need not be Lipschitz continuous,and new ideas appear to be required to treat this range of integrability.    
\bigskip

We conclude this subsection with some possible extensions of Theorems \ref{thm-spectral-inquality} and \ref{thm-null-controllability}. 
\begin{itemize}
    \item \emph{Measurable sets in time.} Can the control in Theorem \ref{thm-null-controllability} be supported in $\omega\times E$, where $E \subset (0,T)$ is a measurable set of positive measure? For potentials in $L^\infty(0,T; L^r(\Omega))$, $r>d$, Phung and Wang \cite{PhungWang2013} proved such a result by deducing the observability inequality from an interpolation inequality at one time through a telescoping series argument. In \cite{apraiz2014observability,burq2022propagation}, this interpolation inequality is derived from a spectral inequality, which suggests that Theorem \ref{thm-spectral-inquality} could be used in the same way, with an explicit dependence on $\|V\|_{L^q(\Omega)} $.
    \item Measurable sets in space. Can the open set in Theorem \ref{thm-spectral-inquality} be replaced by a set of positive measure, in particular in the propagation of smallness from sets of positive measure for solutions of $\Delta u + V u=0$ with $V \in L^q$, $q> \frac{d}{2}$, recently obtained by  Malinnikova and Zhu \cite{MalinnikovaZhu2026}. 
    \item \emph{More general lower-order terms.} For complex-valued potentials, or for first-order terms $W_1\cdot\nabla y+\div(W_2y)$ with $W_j\in L^{q_j}$, the operator $P$ is no longer self-adjoint, and one could rely on the non-self-adjoint version of the method of Lebeau and Robbiano developed by L\'eautaud \cite{leautaud2010spectral}. Both versions of this method, however, rely on the spectral analysis of the time-independent operator $\operatorname{P}$, and do not apply to time-dependent potentials $V(t,x)$, for instance in the mixed classes $L^p_tL^r_x$. Such potentials would require global parabolic $L^p$ Carleman estimates, extending the local ones of Koch and Tataru \cite{koch2009carleman}.
\end{itemize}

\subsection{Stabilization of wave equations with unbounded dampings and potentials}\label{Subsec-stabilization}
The interpolation inequality \eqref{Est-local-up-to-boundary} also yields the stabilization of the damped wave (or Klein--Gordon) equation
\begin{equation}\label{damped-wave-equation}
\left\{
\begin{aligned}
&\partial_t^2u-\div(A\nabla u)+a(x)\,\partial_tu+m(x)\,u=0 && \text{in }(0,+\infty)\times\Omega,\\
&u=0 && \text{on }(0,+\infty)\times\partial\Omega,\\
&u|_{t=0}=u_0,\qquad\partial_tu|_{t=0}=u_1 && \text{in }\Omega,
\end{aligned}
\right.
\end{equation}
where the damping $a$ and the potential $m$ are nonnegative functions such that
\begin{equation}\label{Assumptions-a-m}
a\in L^p(\Omega),\ p>d,\qquad m\in L^q(\Omega),\ q>\tfrac{d+1}2,\qquad a\geq a_0\ \text{ a.e. in }\omega,
\end{equation}
for some $a_0>0$ and some nonempty open set $\omega\subset\Omega$. Let $P=-\div(A\nabla\cdot)+m$, with Dirichlet boundary conditions, which is self-adjoint and positive on $L^2(\Omega)$. Since $a\in L^d(\Omega)$, $au_1\in L^2(\Omega)$ for $u_1\in H^1_0(\Omega)$, and \eqref{damped-wave-equation} is generated by
\begin{equation*}
\mathcal A=\begin{pmatrix}0&I\\-P&-a\end{pmatrix},\qquad
D(\mathcal A)=D(P)\times H^1_0(\Omega),
\end{equation*}
which is maximal dissipative on $\mathcal H=H^1_0(\Omega)\times L^2(\Omega)$.

Following~\cite{Lebeau-1996} (see also~\cite{leellipticV1}), we can obtain a resolvent estimate $\|(\mathcal A-i\lambda)^{-1}\|_{\mathcal L(\mathcal H)}\leq Ce^{C|\lambda|}$, $\lambda\in\R$, which in turn yields the energy decay estimate below, see~\cite[Theorem 3]{Burq-1998}. If $u$ solves $(P-\lambda^2+i\lambda a)u=g$, then $v(s,x)=e^{|\lambda|s}u(x)$ solves the elliptic equation $(-\partial_s^2+P)v=e^{|\lambda|s}(g-i\lambda au)$ in $(0,S_0)\times\Omega\subset\R^{d+1}$, with $v=0$ on $(0,S_0)\times\partial\Omega$, to which we apply \eqref{Est-local-up-to-boundary} in dimension $d+1$: this is why $m$ is taken in $L^q(\Omega)$ with $q>\frac{d+1}2$. The damping term $i\lambda au$ is not treated as a potential, which would make the constants depend on $\lambda$, but as a source term controlled by the dissipation: the identity $\lambda\|\sqrt a\,u\|^2_{L^2(\Omega)}=\operatorname{Im}(g,u)_{L^2(\Omega)}$ controls the observation $\|u\|_{L^2(\omega)}$, and, since $p>d$, $\|au\|_{L^2(\Omega)}\leq C\|\sqrt a\,u\|^\theta_{L^2(\Omega)}\|u\|^{1-\theta}_{H^1(\Omega)}$ with $\theta=\frac{2(p-d)}{2p-d}>0$, which is enough to close the estimate.

\begin{theorem}\label{Thm-stabilization}
For all $k\in\N\setminus\{0\}$, there exists $C>0$ such that
\begin{equation}\label{Log-decay}
\|\partial_tu(t)\|_{L^2(\Omega)}+\|u(t)\|_{H^1(\Omega)}\leq\frac{C}{[\log(2+t)]^k}\,\big\|\mathcal A^k(u_0,u_1)\big\|_{\mathcal H},\qquad t>0,
\end{equation}
for $u$ a solution to \eqref{damped-wave-equation} with $(u_0,u_1)\in D(\mathcal A^k)$. In particular, for $k=1$, the right-hand side of \eqref{Log-decay} is bounded by $\frac C{\log(2+t)}\big(\|u_1\|_{H^1(\Omega)}+\|u_0\|_{D(P)}\big)$.
\end{theorem}

\begin{remark}\label{Rem-stabilization}
When $A$, $a$ and $m$ are smooth, Theorem \ref{Thm-stabilization} is due to Lebeau~\cite{Lebeau-1996}, and the rate $1/\log(2+t)$ is optimal for some geometries. The assumption $a\geq a_0$ on $\omega$ cannot be dropped, since the energy is conserved if $a=0$. The novelty here is that $a$ and $m$ may be unbounded. With source terms in $H^{-1}$ in \eqref{Est-local-up-to-boundary}, which the Carleman estimates of Theorem \ref{Thm2-Carleman-General} allow, the assumption on $a$ could be weakened to $a\in L^{\frac{d}{2}}(\Omega)$, the damping term being then controlled in $H^{-1}(\Omega)$ by the dissipation.
\end{remark}

\subsection{Proof strategies}

\noindent{\bf Strategy for Theorem \ref{Thm-Carleman-General}.} As in \cite{Dehman-Ervedoza-Thabouti-2023}, the proof starts from a model geometric setting, namely a vertical strip $\mathcal S=(0,L)\times\R^{d-1}$, the Laplace operator, and a weight
function which is linear in the tangential variables,
$x\mapsto\varphi(x_1)+\ell\cdot x'$, with $|\ell|$ small (see Section \ref{Sec-Carleman-Strip} for the statements). For such a weight, the conjugated operator \eqref{Conjugated-Operator-strip} has coefficients independent of the tangential variables. A partial Fourier transform in $x'$ therefore reduces it to a family of ordinary differential equations, which can be factorized and solved explicitly.

The main difference with \cite{Dehman-Ervedoza-Thabouti-2023} concerns the reduction to the strip. There, a change of variables built from the weight itself, of the form $y_1=\varphi-\varphi(x_0)$, maps the boundary, which is a level set of the weight, onto a face of the strip.
Here, the weight \eqref{Weight-function-Phi} has a non-constant trace on $\Gamma$, which is used to propagate smallness from the boundary. Such a change of variables would therefore no longer preserve the flat boundary. We instead keep the geometry fixed and linearize the weight:
localizing in a ball of radius $(\varepsilon\tau)^{-1/2}$ centred at $x_0$, we write
\begin{equation*}
\begin{gathered}
\Phi(x)
=\big(\varphi(x_1)+\ell\cdot(x'-x_0')\big)
-\frac{\varepsilon}{2}|x_0'|^2+q(x'),
\\
\ell=-\varepsilon x_0',
\qquad
q(x')=-\frac{\varepsilon}{2}|x'-x_0'|^2.
\end{gathered}
\end{equation*}
The additive constant does not affect the conjugated operator, while the quadratic remainder satisfies $\tau\|q\|_{L^\infty}\leq C$ on the localization ball and is removed by the bounded conjugation $w\mapsto e^{-\tau q}w$.
This requirement determines the scale $(\varepsilon\tau)^{-1/2}$, which also allows us to freeze the rescaled Lipschitz coefficients. The local estimates so obtained are finally glued by a
partition of unity at the same scale, the commutator terms being absorbed for $\varepsilon$ small. This
two-step scheme, Carleman estimate for a tangentially linear weight and then localization and
linearization,  originates in the Hilbertian framework, see \cite{Vessella2017Jumps}; its implementation in the $L^p$ framework is, to our knowledge, new.

Accordingly, most of the article focuses on the proof of a Carleman estimate in the strip for the Laplace
operator with the weight $x\mapsto\varphi(x_1)+\ell\cdot x'$, uniformly with respect to $|\ell|\leq\ell_*$.
We do that in several steps. First, we construct explicitly a parametrix, see Section
\ref{Sec-Construction-parametrix}; the approach is inspired by
\cite{sogge1989oscillatory,koch2001carleman}, by works on Carleman estimates for operators with
discontinuous conductivities, for instance \cite{LeRousseau2013interface,Vessella2017Jumps}, and more
precisely by the recent work \cite{Dehman-Ervedoza-Thabouti-2023}. 

The next step is to estimate the operators appearing in the parametrix. As in \cite{Dehman-Ervedoza-Thabouti-2023}, we use the Stein--Tomas Fourier restriction theorem. Here, however, the restriction estimates must be applied to the positive level sets of $\operatorname{Re}\psi$, where
$$
\psi(\xi')=\big(|\xi'|^2-\tau^2|\ell|^2+2i\tau\,\ell\cdot\xi'\big)^{1/2}
$$
is the square root of the tangential symbol of the conjugated operator. These level sets are ellipsoids whose shape depends on $s=\tau|\ell|/\lambda$. We therefore need a restriction theorem which is uniform on this family when $s$ stays bounded (see Theorem \ref{Thm-Fourier-Rest-Sigma-s}). As $s\to\infty$, the curvature of the normalized ellipsoids degenerates and the restriction constants cannot remain uniform. For $|\ell|$ sufficiently small, the low-frequency region $\{\operatorname{Re}\psi<2\tau|\ell|\}$ lies in a region where the conjugated operator is microlocally elliptic. The corresponding operators are therefore estimated directly by Schur's test and
Bernstein's inequality. The remaining high-frequency region is treated as in \cite[Section 6]{Dehman-Ervedoza-Thabouti-2023}, using pointwise bounds on the kernels and $L^p$--$L^q$ estimates for Fourier multipliers.

Finally, in the last section, we pass to the variable-coefficient operator through a perturbation argument at the same scale $\rho=(\varepsilon\tau)^{-1/2}$, in the spirit of
\cite[Section 5]{koch2001carleman}. On a ball of radius $\rho$, the rescaled Lipschitz matrix $A(\varepsilon\cdot)$ differs from its frozen value by $O(\varepsilon\rho)=O(\varepsilon^{1/2}\tau^{-1/2})$.
The resulting perturbation terms are absorbed for $\varepsilon$ small, completing the local step before gluing.
\medskip

\noindent{\bf Strategy for Theorem \ref{Thm2-Carleman-General}.}
The proof follows that of Theorem \ref{Thm-Carleman-General}, with a weight whose normal derivative has the opposite sign, that is, $\varphi'\geq\alpha_1>0$ in the strip model. The reflection $x_1\mapsto L-x_1$ exchanges the two faces of the strip and reverses the sign of $\varphi'$. The reflected function has zero Cauchy data on $\{x_1=0\}$, since the localized function and the source terms vanish near the interior face, and its Dirichlet trace on $\Gamma$ is now prescribed on $\{x_1=L\}$. Hence, in the parametrix of Theorem \ref{Thm-Carleman-General}, the terms involving the boundary terms on $\{x_1=0\}$ vanish, while the backward integration from $\{x_1=L\}$, which used $w=0$ there, now involves the Dirichlet trace and produces one new operator $G^\sharp_\tau$, see \eqref{Parametrix-reverse}--\eqref{Def-G-sharp}. Equivalently, in the original variables, the conjugated operator is factorized in the reverse order, see Remark \ref{Remark-parametrix-for the second Carleman theorem}. In particular, no Neumann trace appears.  The operators acting on the source terms are those of the proof of Theorem \ref{Thm-Carleman-General}, composed with the reflection, and are estimated in the same way. The operator $G^\sharp_\tau$ is estimated from $H^{\frac{1}{2}}(\R^{d-1})$ by the same arguments as the operator $G_\tau$ in the proof of Theorem \ref{Thm-Carleman-General}; see Section \ref{Sec-proof-Thm2}. 
\medskip

\noindent{\bf Strategy for Theorem \ref{Thm-interpolation-inequality}.} First, Theorem \ref{Thm-Carleman-General}, applied near $y$ after flattening the boundary, rescaling and cutting off, gives Lemma \ref{lemma-local-at-the-boundary}, which bounds $u$ in $\mathcal{W}^* \cap \Omega$, with $\mathcal{W}^*=B_y(\varrho)$, in terms of $\|f\|_{L^2(\Omega)}$ and of the Cauchy data of $u$ on $\mathcal{V} \cap \Omega$. Second, we apply Theorem \ref{Thm-interpolation-from-interior} with $\omega = \mathcal{W }^* \cap \Omega$. Combining the two the two estimates gives Theorem \ref{Thm-interpolation-inequality}, with an interpolation exponent equal to the product of the two exponents. 
\medskip

\noindent{\bf Strategy for Theorem \ref{Thm-interpolation-from-interior}.}
We follow the scheme of \cite[Theorem 5.11]{leellipticV1}. From Theorem \ref{Thm2-Carleman-General}, we derive two local interpolation estimates without boundary term, in which $u$ is observed in $L^2$ on a set at positive distance from $\partial\Omega$: an interior one, where the cut-off function vanishes near the artificial flat face, and a boundary one near each point of $\partial\Omega\cap\mathcal U_1$, where the only boundary term in Theorem \ref{Thm2-Carleman-General}, which involves the Dirichlet trace, vanishes since $u=0$ on $\partial\Omega\cap\mathcal U_1$. Iterating the interior estimate gives a three-ball inequality and, since $\Omega$ is connected, chains of balls propagate the smallness of $u$ from a ball contained in $\omega$ to the set $\Omega_\varsigma$ of points at distance greater than $\varsigma$ from $\partial\Omega$, for every $\varsigma>0$. The boundary estimate then crosses the remaining layer in one step, and, by compactness, finitely many such steps, together with some $\Omega_\varsigma$, cover $\mathcal U$. Since the number of steps is bounded independently of $u$, $f$ and $V$, the exponents multiply to some $\alpha\in(0,1)$, and the factors $e^{C\|V\|_{L^q(\Omega)}^{\gamma_d(q)}}$, raised to powers at most one, are absorbed by changing $C$.
\medskip

\noindent{\bf Strategy for Theorem \ref{thm-spectral-inquality}.}
For a finite spectral sum $w$, we introduce an elliptic lifting $u(s,x)$ whose Cauchy data at $s=0$ are $u=0$ and $\partial_su=w$, and which solves the homogeneous augmented equation.  Corollary \ref{corollary-interpolations-Q} controls $u$ in the interior of the cylinder from $w|_\omega$.  Explicit upper and lower bounds on the scalar lifting functions, including for the finitely many negative eigenvalues, then convert this interpolation estimate
into the spectral inequality.
\medskip

\noindent{\bf Strategy for Theorem \ref{thm-null-controllability}.}
After shifting $\operatorname{P} $ to the nonnegative self-adjoint operator $\widetilde {\operatorname{P}} =\operatorname{P}+\operatorname{R}_V$, Theorem \ref{thm-spectral-inquality} gives, by the abstract Lebeau--Robbiano method~\cite{EgidiNakicEtAl2020}, an observability estimate for $(e^{-t\widetilde P})_{t\geq0}$, hence, by duality, the null-controllability of $\partial_t+\widetilde{\operatorname{P}}$. The change of unknown $y(t)=e^{R_Vt}z(t)$ then returns to $\operatorname{P}$, with the factor $e^{T\operatorname{R}_V}$ in the control cost.
\medskip

\noindent{\bf Strategy for Theorem \ref{Thm-stabilization}.}
The proof follows the classical scheme of Lebeau~\cite{Lebeau-1996}, described before Theorem \ref{Thm-stabilization} in Subsection \ref{Subsec-stabilization}, and is not detailed here.
\medskip

\noindent{\bf Outline.} Section \ref{Sec-Proof-interpolation} is devoted to the proof of the interpolation inequalities, and Section \ref{Sec-Proof-spectral} to the spectral inequality and the null-controllability result. The Carleman estimates are proved in Section \ref{Sec-Carleman-Laplacian} in the model case of the Laplacian, and in Section \ref{Sec-Carleman-variable} for variable coefficients. The last section contains the auxiliary Sobolev multiplier estimate used in the gluing argument.
\medskip

\noindent{\bf Notations.}  Let us finally introduce some of the notation that we will use throughout the article: 
	\begin{itemize}[topsep=0pt,itemsep=0pt,parsep=0pt]
		\item For every $x=(x_1,..,x_d)\in \R^{d}$, we set $x=(x_1,x')$, where $x'=(x_2,..,x_d) \in \R^{d-1}$.  
		\item For $x \in \R^d$ and $r > 0$, $B_{x}(r)$ denotes the ball centered at $x \in \R^d$ and of radius $r > 0$. 
		\item  The symbols $\nabla$ and $\Delta$ denote the gradient and the Laplacian with respect to $x=(x_1,..,x_d)$, whereas $\nabla'=(\partial_{2},..,\partial_{{d}})$ and $\Delta'= \sum_{j=2}^{d} \partial_{j}^2$ denote the tangential gradient and Laplacian with respect to $x'$. 
        \item On $\partial \Omega$, we write $\partial_{\nu_A}u= (A\nabla u) \cdot\nu$ for the conormal derivative, where $\nu$ is the the outward unit normal. 
		\item  The Fourier transform is always taken to be the Fourier transform with respect to $x'=(x_2,..,x_d)$; its dual variable is denoted by $\xi'=(\xi_2,..,\xi_d)$. For a function $f$ defined on $\R^d$ (or a vertical strip $(X_0, X_1) \times \R^{d-1}$) such that $f(x_1,\cdot ) \in \mathscr{S}(\R^{d-1})$, $\hat f(x_1, \cdot)$ denotes the partial Fourier transform with respect to $x'$, namely  
		\begin{equation}
		\label{Partial-Fourier-Transform}
		\hat f(x_1,\xi') = \frac{1}{(2\pi)^{\frac{d-1}{2}}} \int_{\R^{d-1}} e^{-\i  x' \cdot \xi'} f(x_1,x') \, dx', \qquad \qquad\xi' \in \R^{d-1}.
		\end{equation}
		\item For a measurable subset $E$ of $\R^d$, we denote by $|E|$ its Lebesgue measure. 
        \item For $\Omega$ and $\mathcal{V}$ open subsets of $\R^d$, we set
        \begin{equation*}
            \overline{C_c^\infty}^{\,\Omega}( \mathcal{V})
          :=\left\{u=\tilde u|_{\Omega\cap \mathcal{V}}\,;\,
\tilde u \in C_c^\infty(\R^d),\ \operatorname{supp}(\tilde u)\subset \mathcal{V} \right\}.
        \end{equation*}
\end{itemize}
\section{Proof of the interpolation inequalities} \label{Sec-Proof-interpolation}
The proof of Theorem \ref{Thm-interpolation-inequality} follows the same strategy of \cite[Theorem 5.13]{leellipticV1} and combine tow propagation steps. First, Theorem \ref{Thm-Carleman-General}, whose weight increases towards the observation part $\Gamma$ of the boundary, yields Lemma \ref{lemma-local-at-the-boundary} below, an interpolation estimate in a neighborhood of a point of $\Gamma$ from the Cauchy data on $\Gamma$. Second, Theorem \ref{Thm2-Carleman-General}, in which the normal derivative of the weight has the opposite sign, yields Theorem \ref{Thm-interpolation-from-interior}. Accordingly, we first prove Lemma \ref{Thm-interpolation-inequality} and Theorem \ref{Thm-interpolation-from-interior}, and then combine them to prove Theorem \ref{Thm-interpolation-inequality}. 
\medskip

In this part, we assume that $d\geq3$, $\Omega\subset\R^d$ is a bounded domain of class $C^{1,1}$, and $V\in L^q(\Omega)$ with $q>\frac{d}{2}$.
\medskip 

We begin with the following result, which is inspired by \cite[Lemma 5.12]{leellipticV1}. 
\begin{lemma}[Local quantitative unique continuation at the boundary] \label{lemma-local-at-the-boundary} Let $A$ satisfy \eqref{Ellipticity-A}--\eqref{Lipschitz-condition-A}. Let $y \in \partial \Omega$ and let $\mathcal{V}$ be an open neighbourhood of $y$ in $\R^d$. There exist $\varrho>0$, $C>0$ and $\alpha \in (0,1) $ such that, writing $\mathcal{W}^*= B_y(\varrho)$, 
\begin{equation}\label{Est-local-at-the-boundary}
       \|u\|_{H^1(\mathcal{W}^* \cap \Omega)} \leq  C e^{
	  C \norm{V}_{L^{q}(\Omega)}^{\gamma_d(q)}} \|u\|_{H^1(\Omega)}^{1-\alpha} \left(
\|f\|_{L^2(\Omega)}
+ \|u|_{\partial\Omega}\|_{H^\frac{1}{2}(\mathcal{V} \cap \partial\Omega)}
+ \| \partial_{\nu_A} u\|_{H^{-\frac{1}{2}}(\mathcal{V} \cap \partial\Omega)}
\right)^{\alpha},
\end{equation}
 for all $u\in H^1(\Omega)$ satisfying \eqref{Elliptic-equation} with $f \in L^2(\Omega)$, where $\gamma_d$ is given by \eqref{Conditions-gamma(q)}.   
\end{lemma}
The proof of this lemma is inspired by \cite{Hormander1983,Vessella2016Three-region}. 
\begin{proof} For $V \in L^q(\Omega)$, we consider decompositions of the form 
\begin{equation*}
    V = V_{\frac{d}{2}} +V_d + V_\infty, \qquad V_{\frac{d}{2}} \in L^{\frac{d}{2}}(\Omega), \quad V_d \in L^d(\Omega), \quad V_\infty \in L^\infty(\Omega), 
\end{equation*}
which will be specified in the quantification step. 
\medskip

\noindent\textit{Step 1: flattening the boundary, weight and cut-off function.} Since $\partial \Omega$ is of class $C^{1,1}$, after a translation and a rotation, we may assume that $y=0$ and that there exist a neighbourhood $\mathcal{U}$ of $0$ and $g\in C^{1,1}(\R^{d-1}) $ with $g(0)=0$ such that $\mathcal{U} \cap \Omega= \mathcal{U} \cap \{z_1 > g(z')\}$. The map $\Psi (z):=(z_1-g(z'),z')$ is a $C^{1,1}$ diffeomorphism of $\R^d$ with $\det D \Psi=1$, which maps $\mathcal{U} \cap \Omega$ onto $\Psi (\mathcal{U}) \cap \{x_1 >0\}$ and $\mathcal{U} \cap \partial \Omega$ onto $\Psi (\mathcal{U}) \cap \{x_1 = 0\}$. The function $u\circ \Psi^{-1}$ satisfies an equation of the form \eqref{Elliptic-equation} in $\Psi(\mathcal{U}) \cap \{x_1 >0\}$.  Its matrix is symmetric, uniformly elliptic and Lipschitz continuous, with constants depending only $\Lambda_0, \Lambda_1$ and $\|g\|_{C^{1,1}}$, and its potential and source terms have the same $L^p$ norms on corresponding sets. Moreover, on the corresponding subsets of the boundary, the norms $H^{\frac{1}{2}}$ and $H^{-\frac{1}{2}}$ for the Dirichlet traces and the conormal derivative are, respectively, comparable. We keep the notations $u, A, V, f$ and $\partial_{\nu_A}$ in the new coordinates. 

For $r>0$, set $B_r^+:=B_0(r)\cap\{x_1>0\}$ and $\Gamma_r:=B_0(r)\cap\{x_1=0\}$. Fix
\begin{equation}\label{r0-choice}
r_0\in\Big(0,\frac14\Big]\qquad\text{such that}\qquad\overline{B_0(2r_0)}\subset\Psi(\mathcal U),
\end{equation}
so that the flattened domain coincides with $B_{2r_0}^+$ in $B_0(2r_0)$, and let $r_{\mathcal V}>0$ be such that $B_0(r_{\mathcal V})\subset\Psi(\mathcal U\cap\mathcal V)$. Let $\varepsilon_2\in(0,1]$ be given by Theorem \ref{Thm-Carleman-General}, and fix
\begin{equation*}
\varepsilon:=\min\Big\{\varepsilon_2,\ \frac{r_{\mathcal V}}{r_0}\Big\}\in(0,1].
\end{equation*}
We consider the weight
\begin{equation*}
\Phi(x)=\varphi(x_1)-\frac\varepsilon2|x'|^2,\qquad\varphi(x_1)=-x_1+\frac{x_1^2}2,
\end{equation*}
 Set $\rho:=\frac{\varepsilon r_0^2}8$, let $\chi_1\in C^\infty(\R;[0,1])$ satisfy
\begin{equation}\label{Def-weight-1}
\chi_1(t)=\begin{cases}0,&t\leq-3\rho,\\ 1,&t\geq-2\rho,\end{cases}
\end{equation}
and define $\chi(x):=\chi_1(\Phi(x))$. Let $K:=\{0\leq x_1\leq4\rho,\ |x'|\leq\frac{\sqrt3}2r_0\}$. We claim that
\begin{equation}\label{support-chi}
\operatorname{supp}\chi\cap\overline{B_{2r_0}^+}\subset K\Subset B_0(r_0).
\end{equation}
Indeed, let $x\in\overline{B_{2r_0}^+}$ with $\Phi(x)\geq-3\rho$. Since $\varphi\leq0$ on $[0,2r_0]$, we have $\frac\varepsilon2|x'|^2\leq3\rho$, that is, $|x'|\leq\frac{\sqrt3}2r_0$. Moreover, $\varphi(x_1)\leq-\frac{3}{4}x_1$ for $0\leq x_1\leq\frac{1}{2}$, and $\varphi(x_1)\geq\Phi(x)\geq-3\rho$, so that $x_1\leq4\rho$. Finally, $\nabla\chi=\chi_1'(\Phi)\nabla\Phi$, so that
\begin{equation}\label{support-grad-chi}
\operatorname{supp}\nabla\chi\cap\overline{B_{2r_0}^+}\subset\mathcal A_1:=\{-3\rho\leq\Phi\leq-2\rho\}.
\end{equation}
\medskip

\noindent\textit{Step 2: scaled equation.} Set $u_\varepsilon(x)= u(\varepsilon x ), \,  A_\varepsilon(x)= A(\varepsilon x )$ and $v_\varepsilon= \chi  u_\varepsilon $ on $B_{2r_0}^+$. A direct computation gives 
\begin{equation}\label{equation-vepsilon}
-\div(A_\varepsilon\nabla v_\varepsilon)=\varepsilon^2\chi f(\varepsilon\cdot)-\varepsilon^2\chi V(\varepsilon\cdot)u_\varepsilon-2A_\varepsilon\nabla\chi\cdot\nabla u_\varepsilon-u_\varepsilon\div(A_\varepsilon\nabla\chi)\qquad\text{in }B_{2r_0}^+.
\end{equation}
Moreover, with $\nu=-e_1$ and $\Gamma:=\Gamma_{2r_0}$, the Cauchy data of $v_\varepsilon$ on $\Gamma$ are
\begin{equation}\label{Cauchy-data-v}
v_{\varepsilon|x_1=0}=\chi\,u_{\varepsilon|x_1=0},\qquad\partial_{\nu_{A_\varepsilon}}v_\varepsilon=\chi\,\partial_{\nu_{A_\varepsilon}}u_\varepsilon+(A_\varepsilon\nabla\chi\cdot\nu)\,u_{\varepsilon|x_1=0}. 
\end{equation}
By \eqref{support-chi}, both are supported in $K \cap \{x_1 = 0\} \Subset \Gamma _{r_0}$. 
\medskip

\noindent \textit{Step 3: applying the Carleman estimates.} We first justify the application of Theorem \ref{Thm-Carleman-General} to $v_\varepsilon$, which only belongs to $H^{1}(B^+_{2r_0}$. Extending $v_\varepsilon$ to a function of $H^1(\R^d)$, multiplying it by a function of $C_c^\infty(B_0(r_0))$ equal to $1$ in a neighbourhood of $K$, and mollifying, we obtain restrictions $W_n$ to $B_{2r_0}^+$ of functions of $C_c^{\infty}(B_0(r_0))$ such that $W_n \to v_\varepsilon$ in $H^{1}(B_{2r_0}^+)$. Let the right hand side of \eqref{equation-vepsilon} be written as $\tilde f_{2,\varepsilon} + \tilde f_{2_*',\varepsilon},$ with $\tilde f_{2,\varepsilon} \in L^2(B_{2r_0}^+)$ and $\tilde f_{2_*',\varepsilon}\in L^{\frac{2d}{d+2}}(B_{2r_0}^+)$. The two decompositions that we use are given below. Then
\begin{equation*}
-\div(A_\varepsilon\nabla W_n)=\tilde f_{2,\varepsilon}+\tilde f_{2_*',\varepsilon}+\div F_n,\qquad F_n:=-A_\varepsilon\nabla(W_n-v_\varepsilon)\xrightarrow[n\to\infty]{}0\ \text{ in }L^2(B_{2r_0}^+),
\end{equation*}
and $A_\varepsilon\nabla W_n+F_n=A_\varepsilon\nabla v_\varepsilon$, so that the boundary flux $(A_\varepsilon\nabla W_n+F_n)\cdot\nu$ equals $\partial_{\nu_{A_\varepsilon}}v_\varepsilon$ for every $n$. We apply Theorem \ref{Thm-Carleman-General} to $W_n$ and let $n\to\infty$. Since $W_n \to v_\varepsilon$ in $H^1(B^+_{2r_0})$, hence in $L^\frac{2d}{d-2}(B_{2r_0}^+)$, ${W}_{n|x_1=0} \to v_\varepsilon$, we obtain \eqref{Carleman-General-1} and \eqref{Carleman-General-2} for $v_\varepsilon$, with $  F =0$.  In both estimates, the boundary terms are
\begin{equation*}
\mathcal T(\tau):=\Big( \tau^{\frac{1}{2}}\|(e^{\tau\Phi}v_\varepsilon)_{|x_1=0}\|_{L^2(\Gamma)} +\|(e^{\tau\Phi}v_\varepsilon)_{|x_1=0}\|_{H^{\frac{1}{2}}(\Gamma)}+\|e^{\tau\Phi}\partial_{\nu_{A_\varepsilon}}v_\varepsilon\|_{H^{-\frac{1}{2}}(\Gamma)}\Big),
\end{equation*}
and we also set
\begin{equation*}
\mathcal R(\tau):=\varepsilon^2\|e^{\tau\Phi}\chi f(\varepsilon\cdot)\|_{L^2(B_{2r_0}^+)}+\big\|e^{\tau\Phi}\big(2A_\varepsilon\nabla\chi\cdot\nabla u_\varepsilon+u_\varepsilon\div(A_\varepsilon\nabla\chi)\big)\big\|_{L^2(B_{2r_0}^+)}.
\end{equation*}

For \eqref{Carleman-General-1}, we choose
\begin{align*}
&\tilde f_{2,\varepsilon}=\varepsilon^2\chi\big(f(\varepsilon\cdot)-V_\infty(\varepsilon\cdot)u_\varepsilon\big)-2A_\varepsilon\nabla\chi\cdot\nabla u_\varepsilon-u_\varepsilon\div(A_\varepsilon\nabla\chi),\\
&\tilde f_{2_*',\varepsilon}=-\varepsilon^2\chi\big(V_{\frac{d}{2}}(\varepsilon\cdot)+V_d(\varepsilon\cdot)\big)u_\varepsilon.
\end{align*}
By H\"older's inequality, since $\chi u_\varepsilon=v_\varepsilon$,
\begin{equation*}
\|e^{\tau\Phi}\tilde f_{2_*',\varepsilon}\|_{L^{\frac{2d}{d+2}}(B_{2r_0}^+)}\leq\varepsilon^2\|V_{\frac{d}{2}}(\varepsilon\cdot)\|_{L^{\frac{d}{2}}(B_{2r_0}^+)}\|e^{\tau\Phi}v_\varepsilon\|_{L^{\frac{2d}{d-2}}(B_{2r_0}^+)}
+\varepsilon^2\|V_d(\varepsilon\cdot)\|_{L^d(B_{2r_0}^+)}\|e^{\tau\Phi}v_\varepsilon\|_{L^2(B_{2r_0}^+)},
\end{equation*}
and we obtain
\begin{multline}\label{Carleman-applied-1}
\tau^{\frac{3}{2}}\|e^{\tau\Phi}v_\varepsilon\|_{L^2(B_{2r_0}^+)}+\tau^{\frac{1}{2}}\|e^{\tau\Phi}\nabla v_\varepsilon\|_{L^2(B_{2r_0}^+)}\leq C\Big(\mathcal R(\tau)+\tau^{\frac{3}{4}} \mathcal T(\tau)+\varepsilon^2\|V_\infty(\varepsilon\cdot)\|_{L^\infty(B_{2r_0}^+)}\|e^{\tau\Phi}v_\varepsilon\|_{L^2(B_{2r_0}^+)}\\
+\tau^{\frac{3}{4}-\frac1{2d}}\varepsilon^2\|V_d(\varepsilon\cdot)\|_{L^d(B_{2r_0}^+)}\|e^{\tau\Phi}v_\varepsilon\|_{L^2(B_{2r_0}^+)}+\tau^{\frac{3}{4}-\frac1{2d}}\varepsilon^2\|V_{\frac{d}{2}}(\varepsilon\cdot)\|_{L^{\frac{d}{2}}(B_{2r_0}^+)}\|e^{\tau\Phi}v_\varepsilon\|_{L^{\frac{2d}{d-2}}(B_{2r_0}^+)}\Big).
\end{multline}
For \eqref{Carleman-General-2}, we put the $V_d$ term in $\tilde f_{2,\varepsilon}$, since $\varepsilon^2\|e^{\tau\Phi}\chi V_d(\varepsilon\cdot)u_\varepsilon\|_{L^2(B_{2r_0}^+)}\leq\varepsilon^2\|V_d(\varepsilon\cdot)\|_{L^d(B_{2r_0}^+)}\|e^{\tau\Phi}v_\varepsilon\|_{L^{\frac{2d}{d-2}}(B_{2r_0}^+)}$:
\begin{align*}
\tilde f_{2,\varepsilon}&=\varepsilon^2\chi\big(f(\varepsilon\cdot)-(V_d(\varepsilon\cdot)+V_\infty(\varepsilon\cdot))u_\varepsilon\big)-2A_\varepsilon\nabla\chi\cdot\nabla u_\varepsilon-u_\varepsilon\div(A_\varepsilon\nabla\chi),\\
\tilde f_{2_*',\varepsilon}&=-\varepsilon^2\chi V_{\frac{d}{2}}(\varepsilon\cdot)u_\varepsilon,
\end{align*}
and we obtain
\begin{multline}\label{Carleman-applied-2}
\tau^{\frac{3}{4}+\frac1{2d}}\|e^{\tau\Phi}v_\varepsilon\|_{L^{\frac{2d}{d-2}}(B_{2r_0}^+)}\leq C\Big(\mathcal R(\tau)+\tau^{\frac{3}{4}+\frac1{2d}}\mathcal T(\tau)+\varepsilon^2\|V_\infty(\varepsilon\cdot)\|_{L^\infty(B_{2r_0}^+)}\|e^{\tau\Phi}v_\varepsilon\|_{L^2(B_{2r_0}^+)}\\
+\varepsilon^2\|V_d(\varepsilon\cdot)\|_{L^d(B_{2r_0}^+)}\|e^{\tau\Phi}v_\varepsilon\|_{L^{\frac{2d}{d-2}}(B_{2r_0}^+)}+\tau^{\frac{3}{4}+\frac1{2d}}\varepsilon^2\|V_{\frac{d}{2}}(\varepsilon\cdot)\|_{L^{\frac{d}{2}}(B_{2r_0}^+)}\|e^{\tau\Phi}v_\varepsilon\|_{L^{\frac{2d}{d-2}}(B_{2r_0}^+)}\Big).
\end{multline}
We combine \eqref{Carleman-applied-1} and \eqref{Carleman-applied-2}, $C$ denoting the larger of the two constants, and assume that
\begin{equation}\label{Constraints-1}
4C\varepsilon^2\|V_\infty(\varepsilon\cdot)\|_{L^\infty(B_{2r_0}^+)}\leq\tau^{\frac{3}{2}},\qquad4C\varepsilon^2\|V_d(\varepsilon\cdot)\|_{L^d(B_{2r_0}^+)}\leq\tau^{\frac{3}{4}+\frac1{2d}},\qquad4C\varepsilon^2\|V_{\frac{d}{2}}(\varepsilon\cdot)\|_{L^{\frac{d}{2}}(B_{2r_0}^+)}\leq1.
\end{equation}
Absorbing them in the left-hand side, we obtain, for all $\tau\geq\tau_0$ satisfying \eqref{Constraints-1},
\begin{equation}\label{Carleman-combined-fixed}
\tau^{\frac{3}{2}}\|e^{\tau\Phi}v_\varepsilon\|_{L^2(B_{2r_0}^+)}+\tau^{\frac{1}{2}}\|e^{\tau\Phi}\nabla v_\varepsilon\|_{L^2(B_{2r_0}^+)}+\tau^{\frac{3}{4}+\frac1{2d}}\|e^{\tau\Phi}v_\varepsilon\|_{L^{\frac{2d}{d-2}}(B_{2r_0}^+)}
\leq C\big(\mathcal R(\tau)+\tau^{\frac{3}{4}+\frac{1}{2d}}\mathcal T(\tau)\big).
\end{equation}
\smallskip

\begin{figure}[ht]\centering
\begin{tikzpicture}[xscale=1.1,yscale=0.8,
  lab/.style={font=\scriptsize,inner sep=1pt}]
\fill[green!10] plot[domain=-2:2,samples=80] (\x,{2-0.5*\x*\x}) -- cycle;
\fill[blue!22] plot[domain=-2.4495:2.4495,samples=120] (\x,{3-0.5*\x*\x})
  -- plot[domain=2:-2,samples=120] (\x,{2-0.5*\x*\x}) -- cycle;
\fill[green!45] (-1,0) rectangle (1,0.5);
\draw[orange!80!black,thick,dashed] plot[domain=-2.4495:2.4495,samples=120] (\x,{3-0.5*\x*\x});
\draw[blue,thick] plot[domain=-2:2,samples=120] (\x,{2-0.5*\x*\x});
\draw[teal,thick] plot[domain=-1.4142:1.4142,samples=120] (\x,{1-0.5*\x*\x});
\draw[densely dashed] (-2.4495,0) rectangle (2.4495,4);
\draw[very thick] (-3.4,0) -- (3.4,0) node[right]{\small$x'$};
\draw[orange!90!black,line width=2.5pt] (-2.8284,0) -- (2.8284,0);
\draw[->] (-3.4,0) -- (-3.4,0.9) node[left]{\small$x_1$};
\fill (0,0) circle (1.3pt) node[below]{\small$0$};
\node[below,text=orange!80!black] at (2.1,-0.05) {\small$\Gamma_{r_0}$};
\node[lab] at (0,0.25) {$\mathcal W$};
\node[lab] at (1.75,0.97) {$\mathcal A_0$};
\node[below left,lab] at (2.4495,4) {$K$};
\node[above,lab,text=orange!80!black] at (0,3) {$\Phi=-3\rho$};
\node[below,lab,text=blue] at (0,2) {$\Phi=C_1$};
\node[below,lab,text=teal] at (0,1) {$\Phi=C_2$};
\end{tikzpicture}
\caption{Geometry of the proof of Lemma \ref{lemma-local-at-the-boundary}, in the variables of Step~2. The cut-off function $\chi=\chi_1(\Phi)$ is supported in $\{\Phi\geq-3\rho\}\subset K$, equals $1$ in the light green region $\{\Phi\geq C_1\}$, and its gradient is supported in $\mathcal A_0$ (blue), where $\Phi\leq C_1$. The Cauchy data of $v_\varepsilon$ are supported in $K\cap\{x_1=0\}\subset\Gamma_{r_0}$ (orange), where $\Phi\leq0$. On $\mathcal W$ (green), $\Phi\geq C_2$. In the original variables, $\Gamma_{r_0}$ is mapped into $\mathcal V\cap\partial\Omega$, and the image of $\mathcal W$ contains $B_y(\varrho)\cap\Omega$.}
\label{Geometry-local-at-the-boundary}
\end{figure}
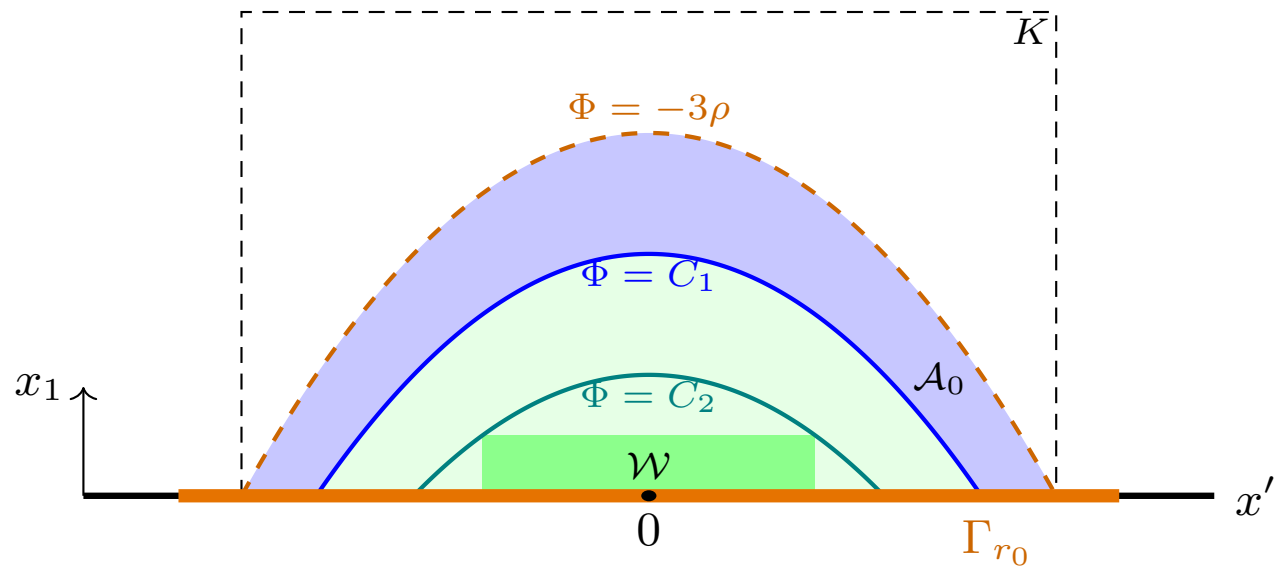

\noindent\textit{Step 4: geometric bounds for the weight.} We now bound each term in $\mathcal{R}(\tau)$ and $\mathcal{T}(\tau)$. 
\smallskip

 \textit{Cut-off region.} By \eqref{support-grad-chi}, the second term of $\mathcal R(\tau)$ is supported in $\mathcal A_0:=\mathcal A_1\cap B_{2r_0}^+$, on which $\Phi\leq-2\rho=:C_1$. Since $\chi$ is fixed and $A_\varepsilon$ is bounded and Lipschitz continuous, $\|\nabla\chi\|_{L^\infty(B_{2r_0}^+)}+\|\div(A_\varepsilon\nabla\chi)\|_{L^\infty(B_{2r_0}^+)}\leq C$, hence
\begin{equation}\label{cutoff-gradient}
\big\|e^{\tau\Phi}\big(2A_\varepsilon\nabla\chi\cdot\nabla u_\varepsilon+u_\varepsilon\div(A_\varepsilon\nabla\chi)\big)\big\|_{L^2(B_{2r_0}^+)}\leq Ce^{C_1\tau}\|u_\varepsilon\|_{H^1(B_{2r_0}^+)}.
\end{equation}

\textit{Interior region.} Set $\mathcal W:=\{0<x_1<\frac\rho2,\ |x'|<\sqrt{\rho/\varepsilon}\}$, where $\sqrt{\rho/\varepsilon}=\frac{r_0}{2\sqrt2}$, so that $\mathcal W\subset B_{2r_0}^+$. On $\mathcal W$, $\varphi(x_1)\geq-x_1\geq-\frac\rho2$ and $\frac\varepsilon2|x'|^2\leq\frac\rho2$, hence $\Phi\geq-\rho=:C_2$. In particular, $\chi=1$ and $v_\varepsilon=u_\varepsilon$ on $\mathcal W$. Therefore, for $\tau\geq1$,
\begin{equation}\label{interior-fixed}
e^{C_2\tau}\|u_\varepsilon\|_{H^1(\mathcal W)}\leq\tau^{\frac{3}{2}}\|e^{\tau\Phi}v_\varepsilon\|_{L^2(B_{2r_0}^+)}+\tau^{\frac{1}{2}}\|e^{\tau\Phi}\nabla v_\varepsilon\|_{L^2(B_{2r_0}^+)}.
\end{equation}

\textit{Source term.} Since $\Phi\leq0$ on $B_{2r_0}^+$, we have 
\begin{equation}\label{source-fixed}
\varepsilon^2\|e^{\tau\Phi}\chi f(\varepsilon\cdot)\|_{L^2(B_{2r_0}^+)}\leq\varepsilon^2\|f(\varepsilon\cdot)\|_{L^2(B_{2r_0}^+)}.
\end{equation}

\textit{Boundary terms.}  By \eqref{Cauchy-data-v}, on $\Gamma$,
\begin{equation*}
e^{\tau\Phi}v_\varepsilon=m_\tau u_\varepsilon,\qquad e^{\tau\Phi}\partial_{\nu_{A_\varepsilon}}v_\varepsilon=m_\tau\partial_{\nu_{A_\varepsilon}}u_\varepsilon+\tilde m_\tau u_\varepsilon,\qquad m_\tau:=e^{\tau\Phi}\chi,\quad\tilde m_\tau:=e^{\tau\Phi}A_\varepsilon\nabla\chi\cdot\nu,
\end{equation*}
where $m_\tau$ and $\tilde m_\tau$ are supported in $K\cap\{x_1=0\}\Subset\Gamma_{r_0}$, and satisfy $\|m_\tau\|_{W^{1,\infty}(\Gamma)}+\|\tilde m_\tau\|_{W^{1,\infty}(\Gamma)}\leq C\tau$ for $\tau\geq1$, since $\Phi\leq0$ and $|\nabla\Phi|\leq C$ on $\Gamma$. We then obtain
\begin{equation*}
\mathcal T(\tau)\leq C\tau \Big(\|u_{\varepsilon|x_1=0}\|_{H^{\frac{1}{2}}(\Gamma_{r_0})}+\|\partial_{\nu_{A_\varepsilon}}u_\varepsilon\|_{H^{-\frac{1}{2}}(\Gamma_{r_0})}\Big).
\end{equation*}
We set $C_3:=\rho$. Since $\tau^{\frac{7}{4}+\frac{1}{2d}}\leq Ce^{C_3\tau}$ for $\tau\geq1$,
\begin{equation}\label{boundary-estimates}
\tau^{\frac{3}{4}+\frac{1}{2d}}\mathcal T(\tau)\leq Ce^{C_3\tau}\Big(\|u_{\varepsilon|x_1=0}\|_{H^{\frac{1}{2}}(\Gamma_{r_0})}+\|\partial_{\nu_{A_\varepsilon}}u_\varepsilon\|_{H^{-\frac{1}{2}}(\Gamma_{r_0})}\Big).
\end{equation}
Note that $C_1=-2\rho<C_2=-\rho<C_3=\rho$.
\medskip

\noindent\textit{Step 5: combined estimate and rescaling.}  Inserting \eqref{cutoff-gradient}--\eqref{boundary-estimates} into \eqref{Carleman-combined-fixed}, and using \eqref{interior-fixed}, we obtain, for all $\tau\geq\tau_0$ satisfying \eqref{Constraints-1},
\begin{multline}\label{Est-before-optimizing-fixed}
e^{C_2\tau}\|u_\varepsilon\|_{H^1(\mathcal W)}\leq Ce^{C_1\tau}\|u_\varepsilon\|_{H^1(B_{2r_0}^+)}\\
+Ce^{C_3\tau}\Big(\varepsilon^2\|f(\varepsilon\cdot)\|_{L^2(B_{2r_0}^+)}+\|u_{\varepsilon|x_1=0}\|_{H^{\frac{1}{2}}(\Gamma_{r_0})}+\|\partial_{\nu_{A_\varepsilon}}u_\varepsilon\|_{H^{-\frac{1}{2}}(\Gamma_{r_0})}\Big).
\end{multline}
We now return to the variable $z=\varepsilon x$. We have $u_{\varepsilon|x_1=0}=u_{|x_1=0}(\varepsilon\cdot)$ and $\partial_{\nu_{A_\varepsilon}}u_\varepsilon=\varepsilon(\partial_{\nu_A}u)(\varepsilon\cdot)$ and recall that $\varepsilon>0$ is fixed in terms only of the geometry and of the neighbourhood $\mathcal V$. Since $\varepsilon\Gamma_{r_0}\subset\Gamma_{r_{\mathcal V}}$, we obtain, after returning to the original coordinates,
\begin{align*}
&\|u_\varepsilon\|_{H^1(B_{2r_0}^+)}\leq C\|u\|_{H^1(\Omega)},\quad && \varepsilon^2\|f(\varepsilon\cdot)\|_{L^2(B_{2r_0}^+)}\leq C\|f\|_{L^2(\Omega)},\\
&\|u_{\varepsilon|x_1=0}\|_{H^{\frac{1}{2}}(\Gamma_{r_0})}\leq C\|u_{|\partial\Omega}\|_{H^{\frac{1}{2}}(\mathcal V\cap\partial\Omega)},\quad&&\|\partial_{\nu_{A_\varepsilon}}u_\varepsilon\|_{H^{-\frac{1}{2}}(\Gamma_{r_0})}\leq C\|\partial_{\nu_A}u\|_{H^{-\frac{1}{2}}(\mathcal V\cap\partial\Omega)}.
\end{align*}
Moreover, $\varepsilon\mathcal W\supset B_{r_1}^+$ with $r_1:=\varepsilon\min\{\frac\rho2,\frac{r_0}{2\sqrt2}\}$. Since $\Psi$ is bi-Lipschitz, there exists $\varrho>0$ such that $B_y(\varrho)\cap\Omega\subset\Psi^{-1}(\varepsilon\mathcal W)$, and then $\|u\|_{H^1(B_y(\varrho)\cap\Omega)}\leq C\|u_\varepsilon\|_{H^1(\mathcal W)}$. Setting $\mathcal W^*:=B_y(\varrho)$, we obtain
\begin{equation}\label{Est-before-optimizing}
e^{C_2\tau}\|u\|_{H^1(\mathcal W^*\cap\Omega)}\leq Ce^{C_1\tau}\|u\|_{H^1(\Omega)}+Ce^{C_3\tau}\mathcal B(u),
\end{equation}
where
\begin{equation*}
\mathcal B(u):=\|f\|_{L^2(\Omega)}+\|u_{|\partial\Omega}\|_{H^{\frac{1}{2}}(\mathcal V\cap\partial\Omega)}+\|\partial_{\nu_A}u\|_{H^{-\frac{1}{2}}(\mathcal V\cap\partial\Omega)}.
\end{equation*}
Dividing \eqref{Est-before-optimizing} by $e^{C_2\tau}$ gives
\begin{equation}\label{Est-divided}
\|u\|_{H^1(\mathcal W^*\cap\Omega)}\leq Ce^{-(C_2-C_1)\tau}\|u\|_{H^1(\Omega)}+Ce^{(C_3-C_2)\tau}\mathcal B(u).
\end{equation}

\medskip
\noindent\textit{Step 6: quantification with respect to the potential.}
For $p\in\{\frac{d}{2},d,\infty\}$, since $\varepsilon\leq1$ and $\Psi$ preserves the Lebesgue measure,
\begin{equation*}
\varepsilon^2\|V_p(\varepsilon\cdot)\|_{L^p(B_{2r_0}^+)}=\varepsilon^{2-\frac{d}{p}}\|V_p\|_{L^p(B_{2\varepsilon r_0}^+)}\leq\|V_p\|_{L^p(\Omega)}.
\end{equation*}
Hence \eqref{Constraints-1} holds as soon as
\begin{equation}\label{Constraints-for-V}
4C\|V_\infty\|_{L^\infty(\Omega)}\leq\tau^{\frac{3}{2}},\qquad4C\|V_d\|_{L^d(\Omega)}\leq\tau^{\frac{3}{4}+\frac1{2d}},\qquad4C\|V_{\frac{d}{2}}\|_{L^{\frac{d}{2}}(\Omega)}\leq1.
\end{equation}
For each $\tau\geq\tau_0$, we choose the decomposition of $V$ as follows; we may assume that $V\neq0$.
\smallskip

\textit{Case $d\leq q\leq\infty$.} Set $V_{\frac{d}{2}}=0$, $V_\infty=V\mathbf 1_{\{|V|\leq\lambda\}}$ and $V_d=V\mathbf 1_{\{|V|>\lambda\}}$, with $\lambda:=\|V\|_{L^q(\Omega)}\tau^{\frac{d}{q}\left(\frac{3}{4}-\frac1{2d}\right)}$. Then $\|V_\infty\|_{L^\infty(\Omega)}\leq\lambda$. If $q=\infty$, then $\lambda=\|V\|_{L^\infty(\Omega)}$ and $V_d=0$. If $q<\infty$, since $q\geq d$,
\begin{equation*}
\|V_d\|_{L^d(\Omega)}\leq\lambda^{1-\frac qd}\|V\|_{L^q(\Omega)}^{\frac qd}=\|V\|_{L^q(\Omega)}\tau^{\left(\frac{3}{4}-\frac1{2d}\right)\left(\frac{d}{q}-1\right)}.
\end{equation*}
In both cases \eqref{Constraints-for-V} holds as soon as
\begin{equation*}
4C\|V\|_{L^q(\Omega)}\leq\tau^{\frac{3}{2}\left(1-\frac d{2q}\right)+\frac1{2q}}.
\end{equation*}

\smallskip
\textit{Case $\frac{d}{2}<q\leq d$.} Set $V_\infty=0$, $V_{\frac{d}{2}}=V\mathbf 1_{\{|V|>\lambda\}}$ and $V_d=V\mathbf 1_{\{|V|\leq\lambda\}}$, with $\lambda:=\|V\|_{L^q(\Omega)}\tau^{\frac{d}{q}\left(\frac{3}{4}+\frac1{2d}\right)}$. Then
\begin{align*}
\|V_{\frac{d}{2}}\|_{L^{\frac{d}{2}}(\Omega)}&\leq\lambda^{1-\frac{2q}d}\|V\|_{L^q(\Omega)}^{\frac{2q}d}=\|V\|_{L^q(\Omega)}\tau^{-\left(\frac{3}{4}+\frac1{2d}\right)\left(2-\frac{d}{q}\right)},\\
\|V_d\|_{L^d(\Omega)}&\leq\lambda^{1-\frac qd}\|V\|_{L^q(\Omega)}^{\frac qd}=\|V\|_{L^q(\Omega)}\tau^{\left(\frac{3}{4}+\frac1{2d}\right)\left(\frac{d}{q}-1\right)},
\end{align*}
so that \eqref{Constraints-for-V} holds as soon as $4C\|V\|_{L^q(\Omega)}\leq\tau^{\left(\frac{3}{4}+\frac1{2d}\right)\left(2-\frac{d}{q}\right)}$.

\smallskip
In all cases, with $\gamma_d(q)$ given by \eqref{Conditions-gamma(q)}, that is,
\begin{equation*}
\gamma_d(q)=\frac1{\left(\frac{3}{4}+\frac1{2d}\right)\left(2-\frac{d}{q}\right)}\quad\text{if }\frac{d}{2}<q\leq d,\qquad\gamma_d(q)=\frac1{\frac{3}{2}\left(1-\frac d{2q}\right)+\frac1{2q}}\quad\text{if }d\leq q\leq\infty,
\end{equation*}
the conditions \eqref{Constraints-1} hold for all
\begin{equation}\label{tau-V}
\tau\geq\tau(V):=\max\Big\{\tau_0,\big(4C\|V\|_{L^q(\Omega)}\big)^{\gamma_d(q)}\Big\}\leq C\big(1+\|V\|_{L^q(\Omega)}^{\gamma_d(q)}\big).
\end{equation}

Finally, using \cite[Lemma 5.4]{leellipticV1} and \eqref{Est-divided}, there exist $C>0$ and $\alpha \in (0,1)$ such that 
\begin{equation*}
\|u\|_{H^1(\mathcal W^*\cap\Omega)}
\leq
C
e^{C\tau(V)}
\|u\|_{H^1(\Omega)}^{1-\alpha}
\left(
\|f\|_{L^2(\Omega)}
+
\|u|_{\partial\Omega}\|_{H^{\frac{1}{2}}(\mathcal V\cap\partial\Omega)}
+
\|\partial_{\nu_A}u|_{\partial\Omega}
\|_{H^{-\frac{1}{2}}(\mathcal V\cap\partial\Omega)}
\right)^\alpha,
\end{equation*}
which, by \eqref{tau-V}, is \eqref{Est-local-at-the-boundary}.
\end{proof}

\begin{proof}[Proof of Theorem \ref{Thm-interpolation-from-interior}] Following the overall strategy of the proof of \cite[Theorem 5.11]{leellipticV1},  we adapt the proof of Lemma \ref{lemma-local-at-the-boundary}, keeping its notation but using Theorem \ref{Thm2-Carleman-General} instead of Theorem \ref{Thm-Carleman-General}. We only detail the necessary modifications: the weight, the cut-off functions, the choice of $\rho$ and the sets $K$, $\mathcal{A}_0$, $\mathcal{A}_1$ and $\mathcal{W}$, which are redefined below. For $\varsigma >0$, we define $\Omega_\varsigma= \{x\in \Omega; \quad\dist (x,\partial \Omega) >\varsigma\} $. We prove below two local estimates, a boundary estimate and an interior estimate. 
\begin{itemize}
    \item \textit{Boundary estimate.} For every $y \in \partial \Omega \cap \mathcal{U}_1,$ there exist $C_y > 0$, $\theta_1 \in (0,1)$ $\varrho_y>0$ and $\varsigma_y>0$ such that 
    \begin{equation}\label{boundary-estimate}
        \|u\|_{H^1(B_y(\varrho_y) \cap \Omega)} \leq C_y e^{C_y \|V\|_{L^q(\Omega)}^{\gamma_d(q)}}  \|u\|_{H^1(\Omega)}^{1-\theta_1} \left(\|f\|_{L^2(\Omega)} + \|u\|_{L^2(\Omega_{\varsigma_y})}\right)^{\theta_1}. 
    \end{equation}
    \item  \textit{Interior estimate.} There exist $s_0>0$, $c\in (0,\frac{1}{2})$, $\theta_2 \in (0,1)$, and for every $s\in (0,s_0]$, a constant $C_s > 0$ such that, for all $x_0 \in \Omega $ and $e\in \S^{d-1}$ with $B_{x_0} (3s)\subset \Omega$, 
      \begin{equation}\label{interior-estimate}
        \|u\|_{H^1(B_{x_0}(2 c s) )} \leq C_s e^{C_s \|V\|_{L^q(\Omega)}^{\gamma_d(q)}}  \|u\|_{H^1(\Omega)}^{1-\theta_2} \left(\|f\|_{L^2(\Omega)} + \|u\|_{L^2(B_{x_0}(s) \cap \{(z-x_0)\cdot e>cs\})}\right)^{\theta_2}. 
    \end{equation}
\end{itemize}
We also deduce from \eqref{interior-estimate} a three ball inequality: for every $r>0$, there exist $C_r>0$ and $\theta _r\in (0,1)$ such that, whenever $B_{x_0}(4r) \subset \Omega$, 
\begin{equation}\label{three-ball}
  \|f\|_{L^2(\Omega)}+   \|u\|_{H^1(B_{x_0}(2 r) )} \leq C_r e^{C_r \|V\|_{L^q(\Omega)}^{\gamma_d(q)}}  \|u\|_{H^1(\Omega)}^{1-\theta_r} \left(\|f\|_{L^2(\Omega)} + \|u\|_{L^2(B_{x_0}(r) )}\right)^{\theta_r}.
\end{equation}
These estimates imply \eqref{Thm-interpolation-from-interior} as in the proofs of \cite[Theorems 5.6 and 5.11]{leellipticV1}, which we recall. Here and below, $C>0$ denotes a constant that may change from line to line and does not depend on $u$, $f$ and $V$. If $\|f\|_{L^2(\Omega)} + \|u\|_{L^2(\omega)} \geq \|u\|_{H^1(\Omega)}$, then \eqref{Est-local-up-to-boundary} is obvious, so we assume that $\|f\|_{L^2(\Omega)} + \|u\|_{L^2(\omega)} \leq \|u\|_{H^1(\Omega)}$. 
\medskip

\noindent\textit{Interior propagation.} We first show that, for every $\varsigma >0$ there exist $C_\varsigma>0$ and $\theta_\varsigma \in (0,1) $ such that 
\begin{equation}
    \|u\|_{H^1(\Omega_\varsigma)} \leq C_\varsigma e^{C_\varsigma \|V\|_{L^q(\Omega)}^{\gamma_d(q)}} \|u\|_{H^1(\Omega )}^{1-\theta_\varsigma} \left(\|f\|_{L^2(\Omega)} + \|u\|_{L^2(\omega )}\right)^{\theta_\varsigma}. 
\end{equation}
Fix $x_\omega\in\omega$ and $r_\omega>0$ with $B_{x_\omega}(r_\omega)\subset\omega$, and let $y\in\overline{\Omega_\varsigma}$. Since $\Omega$ is connected, there exists a continuous path $\Gamma:[0,1]\to\Omega$ from $x_\omega$ to $y$, and $r_\Gamma:=\operatorname{dist}(\Gamma([0,1]),\partial\Omega)>0$. Set $r:=\min\{r_\omega,\frac{r_\Gamma}4\}$. Following $\Gamma$, we define points $x_0:=x_\omega,x_1,\dots,x_N:=y$ of $\Gamma([0,1])$: $x_{j+1}$ is the first point of $\Gamma$ after $x_j$ on the sphere $\{|z-x_j|=r\}$, or $y$ if there is none. The number $N$ is finite since $\Gamma$ is uniformly continuous. Then $|x_{j+1}-x_j|\leq r$, so that $B_{x_{j+1}}(r)\subset B_{x_j}(2r)$, and $B_{x_j}(4r)\subset\Omega$. Applying \eqref{three-ball} at $x_0,\dots,x_{N-1}$, and using $\|u\|_{L^2(B_{x_{j+1}}(r))}\leq\|u\|_{H^1(B_{x_j}(2r))}$ for $j\leq N-2$, $B_y(r)\subset B_{x_{N-1}}(2r)$ and $B_{x_0}(r)\subset\omega$, we obtain
\begin{equation*}
\|u\|_{H^1(B_y(r))}\leq Ce^{C\|V\|_{L^q(\Omega)}^{\gamma_d(q)}}\|u\|_{H^1(\Omega)}^{1-\theta_r^N}\Big(\|f\|_{L^2(\Omega)}+\|u\|_{L^2(\omega)}\Big)^{\theta_r^N}.
\end{equation*}
Covering the compact set $\overline{\Omega_\varsigma}$ by finitely many such balls gives \eqref{interior-estimate}.
\smallskip

\noindent\textit{Propagation up to the boundary.} If $\bar{\mathcal{U}} \cap \partial \Omega= \emptyset$, then $\mathcal{U} \subset \Omega_\varsigma$ for some $\varsigma>0$, and \eqref{Est-local-up-to-boundary} follows from \eqref{interior-estimate}. Otherwise, the compact set $\bar{\mathcal{U}}\cap \partial \Omega$ is contained in $\mathcal{U}_1$, so it is covered by finitely many balls $B_{y_\ell}(\varrho_{y_\ell})$ given  \eqref{boundary-estimate}.By compactness of $\mathcal{U}$, there exists $\varsigma'>0$ such that 
\begin{equation*}
    \tilde{\mathcal{U}}= \{ x \in \mathcal{U}, \quad \dist(x,\partial \Omega ) < \varsigma'\} \subset \bigcup_{\ell} B_{y_\ell} (\varrho _{y_\ell}).
\end{equation*}
Set $\varsigma'':=\min_\ell\varsigma_{y_\ell}$. Since $\Omega_{\varsigma_{y_\ell}}\subset\Omega_{\varsigma''}$, \eqref{boundary-estimate} gives
\begin{equation*}
\|u\|_{H^1(\widetilde{\mathcal U})}\leq Ce^{C\|V\|_{L^q(\Omega)}^{\gamma_d(q)}}\|u\|_{H^1(\Omega)}^{1-\theta_1}\Big(\|f\|_{L^2(\Omega)}+\|u\|_{L^2(\Omega_{\varsigma''})}\Big)^{\theta_1}.
\end{equation*}
Moreover, by \eqref{interior-estimate} with $\varsigma''$, and since $\|f\|_{L^2(\Omega)}\leq\|u\|_{H^1(\Omega)}^{1-\theta_{\varsigma''}}\big(\|f\|_{L^2(\Omega)}+\|u\|_{L^2(\omega)}\big)^{\theta_{\varsigma''}}$ by our assumption,
\begin{equation*}
\|f\|_{L^2(\Omega)}+\|u\|_{L^2(\Omega_{\varsigma''})}\leq Ce^{C\|V\|_{L^q(\Omega)}^{\gamma_d(q)}}\|u\|_{H^1(\Omega)}^{1-\theta_{\varsigma''}}\Big(\|f\|_{L^2(\Omega)}+\|u\|_{L^2(\omega)}\Big)^{\theta_{\varsigma''}}.
\end{equation*}
Combining these two estimates, we obtain
\begin{equation*}
\|u\|_{H^1(\widetilde{\mathcal U})}\leq Ce^{C\|V\|_{L^q(\Omega)}^{\gamma_d(q)}}\|u\|_{H^1(\Omega)}^{1-\theta_1\theta_{\varsigma''}}\Big(\|f\|_{L^2(\Omega)}+\|u\|_{L^2(\omega)}\Big)^{\theta_1\theta_{\varsigma''}}.
\end{equation*}
Finally, $\mathcal U\subset\widetilde{\mathcal U}\cup\Omega_{\varsigma'/2}$, and \eqref{interior-estimate} with $\varsigma'/2$, together with the last estimate, gives \eqref{Est-local-up-to-boundary}.
\medskip

\noindent \textit{Proof of \eqref{boundary-estimate}.} Let $y \in \partial\Omega \cap \mathcal{U}_1 $. We flatten the boundary near $y$ and rescale as in Steps 1 and 2 of the proof of Lemma \ref{lemma-local-at-the-boundary}, with three changes: the chart neighbourhood of $y$ is denoted by $\mathcal{U}_y$, we take $\mathcal{V}= \mathcal{U}_1$, and $\varepsilon_2$ and $\tau_0$ are given by Theorem \ref{Thm2-Carleman-General} for the function $\varphi$ below. Note that \eqref{equation-vepsilon} holds for every cut-off function $\chi$. Since $u=0$ on $\partial \Omega \cap \mathcal{U}_1$ and $\varepsilon \Gamma_{r_0} \subset \Gamma_{r_{\mathcal{V}}}$, we have $u_{\varepsilon| x_1=0}$ on $\Gamma_{r_0}$. We take the weight 
\begin{equation*}
    \Phi(x)= \varphi(x_1) - \frac{\varepsilon}{2}|x'|^2, \qquad \varphi(x_1) = x_1 + \frac{x_1^2}{2}. 
\end{equation*}
Set $\rho = \frac{\varepsilon r_0^2}{16}, $ let $\chi_1 \in C^\infty(\R; [0,1])$ satisfy \eqref{Def-weight-1} with this value of $\rho$, and let $\chi_2 \in C^\infty(\R; [0,1]) $ satisfy 
\begin{equation*} 
\chi_2(t) = \left\{
\begin{aligned}
&0 , 
&&  t \geq 2 \rho,\\
& 1 
&&  t < \rho,
\end{aligned}
\right. 
\end{equation*}
Then, we define $\chi(x)=\chi_1(\Phi(x)) \chi_2(x_1)$ and $v_\varepsilon=\chi u_\varepsilon$.

On $\operatorname{supp} \chi \cap \overline{B_{2r_0}^+}$, we have $x_1 \leq 2 \rho$ and $\Phi \geq -3 \rho,$ hence $\Phi \leq \varphi(2 \rho) \leq 3 \rho$ and $\frac{\varepsilon}{2}|x'|^2 \leq 6 \rho$. Therefore 
\begin{equation}
    \operatorname{supp}\chi \cap  \overline{B_{2r_0}^+} \subset K :=\Big\{|x'|\leq \tfrac{\sqrt  3 }{2} r_0 ,\ 0\leq x_1\leq 2 \rho\Big\} \Subset B_{r_0}^+. 
\end{equation}
and in particular $v_{\varepsilon|x_1=0}=0$ on $\Gamma$. Each term of $\nabla \chi$ and of $\div(A_\varepsilon \nabla \chi)$ contains a derivative of $\chi_1$ or of $\chi_2$. So these functions are supported in $\mathcal{A}_0 \cup \mathcal{A}_1$, where  
\begin{equation*}  
\mathcal A_0:=\big\{0 \leq x_1<2 \rho,\ -3 \rho\leq \Phi \leq - 2 \rho \big\},
\qquad
\mathcal A_1:=\big\{ \rho \leq x_1<2 \rho,\ \Phi \geq -3 \rho  \big\},
\end{equation*}
On $\mathcal A_0$, we have $\Phi \leq - 2 \rho =: C_1 $, and on $\mathcal A_1$, we have $\Phi \leq  \varphi(2 \rho) \leq 3 \rho =: C_3$ (see Figure \eqref{Geometry near the boundary-interior observation}).
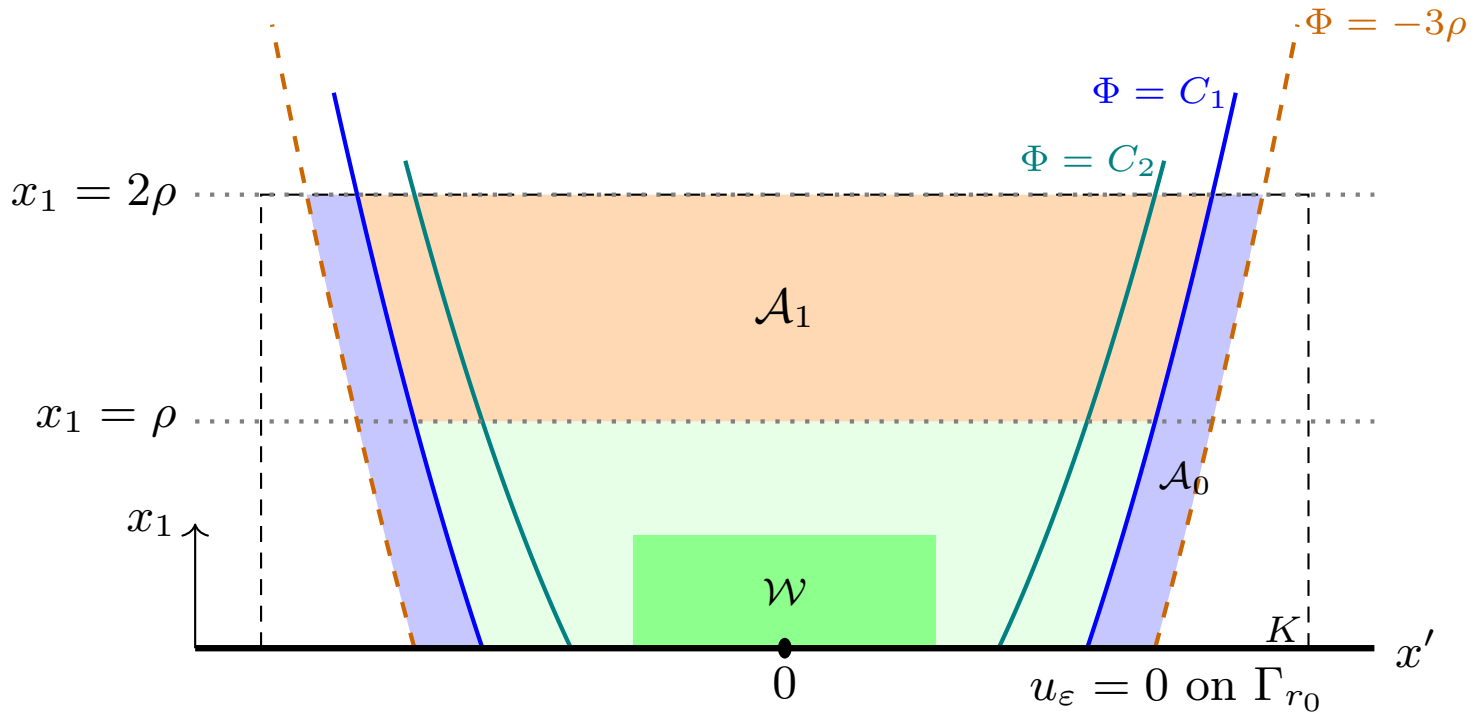
\begin{figure}[ht]\centering
\begin{tikzpicture}[xscale=1,yscale=1.5,
  lab/.style={font=\scriptsize,inner sep=1pt}]
\fill[green!10] (-2,0) -- (2,0) -- plot[domain=2:2.4495,samples=40] (\x,{0.5*\x*\x-2})
  -- plot[domain=-2.4495:-2,samples=40] (\x,{0.5*\x*\x-2}) -- cycle;
\fill[orange!30] (-2.8284,1) -- (2.8284,1)
  -- plot[domain=2.8284:3.1623,samples=40] (\x,{0.5*\x*\x-3})
  -- plot[domain=-3.1623:-2.8284,samples=40] (\x,{0.5*\x*\x-3}) -- cycle;
\foreach \s in {1,-1}
  \fill[blue!22] plot[domain=2:2.8284,samples=60] ({\s*\x},{0.5*\x*\x-2}) -- ({\s*3.1623},2)
    -- plot[domain=3.1623:2.4495,samples=60] ({\s*\x},{0.5*\x*\x-3}) -- cycle;
\fill[green!45] (-1,0) rectangle (1,0.5);
\foreach \s in {1,-1}{
  \draw[orange!80!black,thick,dashed] plot[domain=2.4495:3.3912,samples=60] ({\s*\x},{0.5*\x*\x-3});
  \draw[blue,thick] plot[domain=2:2.9833,samples=60] ({\s*\x},{0.5*\x*\x-2});
  \draw[teal,thick] plot[domain=1.4142:2.51,samples=60] ({\s*\x},{0.5*\x*\x-1});}
\draw[densely dashed] (-3.4641,0) rectangle (3.4641,2);
\draw[very thick] (-3.9,0) -- (3.9,0) node[right]{\small$x'$};
\draw[->] (-3.9,0) -- (-3.9,0.55) node[left]{\small$x_1$};
\draw[gray,dotted,thick] (-3.9,1) -- (3.9,1);
\draw[gray,dotted,thick] (-3.9,2) -- (3.9,2);
\node[left] at (-3.9,1) {\small$x_1=\rho$};
\node[left] at (-3.9,2) {\small$x_1=2\rho$};
\fill (0,0) circle (1.3pt) node[below]{\small$0$};
\node[below] at (2.6,0) {\small$u_\varepsilon=0$ on $\Gamma_{r_0}$};
\node at (0,1.5) {\small$\mathcal A_1$};
\node[lab] at (2.65,0.75) {$\mathcal A_0$};
\node[lab] at (0,0.25) {$\mathcal W$};
\node[above left,lab] at (3.4641,0) {$K$};
\node[right,lab,text=orange!80!black] at (3.3912,2.75) {$\Phi=-3\rho$};
\node[left,lab,text=blue] at (2.9833,2.45) {$\Phi=C_1$};
\node[left,lab,text=teal] at (2.51,2.15) {$\Phi=C_2$};
\end{tikzpicture}
\caption{Geometry of the proof of \eqref{boundary-estimate}, in the flattened and rescaled variables. The cut-off function $\chi=\chi_1(\Phi)\chi_2(x_1)$ is supported in $K$, equals $1$ in the light green region, which contains $\mathcal W$, and its gradient is supported in $\mathcal A_0$ (blue), where $\Phi\leq C_1$, and in $\mathcal A_1$ (orange), where $\Phi\leq C_3$. On $\mathcal W$ (green), $\Phi\geq C_2$. No boundary term appears since $u_\varepsilon=0$ on $\Gamma_{r_0}$. The observation region $\mathcal A_1$ lies in $\{x_1\geq\rho\}$.}
\label{Geometry near the boundary-interior observation}
\end{figure}

Since $v_{\varepsilon|x_1=0}=0$, the boundary terms in Theorem \ref{Thm2-Carleman-General} vanish. We write the commutator in \eqref{equation-vepsilon} in divergence form 
\begin{equation*}
    2 A_\varepsilon \nabla \chi \cdot \nabla u_\varepsilon + u_\varepsilon \div(A_\varepsilon \nabla \chi) = 2 \div( u_\varepsilon A _\varepsilon \nabla \chi) - u_\varepsilon \div(A_\varepsilon \nabla \chi). 
\end{equation*}
Accordingly, in the decompositions of Step~3, we replace $-2A_\varepsilon\nabla\chi\cdot\nabla u_\varepsilon-u_\varepsilon\div(A_\varepsilon\nabla\chi)$ by $u_\varepsilon\div(A_\varepsilon\nabla\chi)$ in $\tilde f_{2,\varepsilon}$, and we take $\widetilde F_\varepsilon:=-2u_\varepsilon A_\varepsilon\nabla\chi$ instead of $F=0$. Otherwise, Steps~3 and~6 of the proof of Lemma \ref{lemma-local-at-the-boundary} apply without change. Since $\Phi\leq C_3$ on $\operatorname{supp}\chi$, and since $u_\varepsilon\div(A_\varepsilon\nabla\chi)$ and $\widetilde F_\varepsilon$ are supported in $\mathcal A_0\cup\mathcal A_1$, we obtain, for all $\tau\geq\tau(V)$,
\begin{multline}\label{Carleman-combined-fixed-2}
\tau^{\frac{3}{2}}\|e^{\tau\Phi}v_\varepsilon\|_{L^2(B_{2r_0}^+)}+\tau^{\frac{1}{2}}\|e^{\tau\Phi}\nabla v_\varepsilon\|_{L^2(B_{2r_0}^+)}\\
\leq C\tau e^{C_1\tau}\|u_\varepsilon\|_{L^2(B_{2r_0}^+)}+C\tau e^{C_3\tau}\Big(\varepsilon^2\|f(\varepsilon\cdot)\|_{L^2(B_{2r_0}^+)}+\|u_\varepsilon\|_{L^2(\mathcal A_1)}\Big).
\end{multline}
On $\mathcal W:=\{0<x_1<\frac\rho2,\ |x'|<\sqrt{\rho/\varepsilon}\}$, we have $\chi=1$ and $\Phi\geq-\rho=:C_2$, so that the left-hand side of \eqref{Carleman-combined-fixed-2} is bounded from below by $e^{C_2\tau}\|u_\varepsilon\|_{H^1(\mathcal W)}$. Set $\eta:=\frac{C_2-C_1}2$, so that $C_1+\eta<C_2<C_3+\eta$. Since $\tau\leq C_\eta e^{\eta\tau}$. Then by optimization, there exist $C>0$ and $\theta_1 \in (0,1)$ such that 
\begin{equation*}
\|u_\varepsilon\|_{H^1(\mathcal W)}\leq Ce^{C\|V\|_{L^q(\Omega)}^{\gamma_d(q)}}\|u_\varepsilon\|_{H^1(B_{2r_0}^+)}^{1-\theta_1}\Big(\varepsilon^2\|f(\varepsilon\cdot)\|_{L^2(B_{2r_0}^+)}+\|u_\varepsilon\|_{L^2(\mathcal A_1)}\Big)^{\theta_1}.
\end{equation*}
Returning to the original variables as in Step~5 of the proof of Lemma \ref{lemma-local-at-the-boundary}, we obtain
\begin{equation*}
\|u\|_{H^1(B_y(\varrho_y)\cap\Omega)}\leq C_ye^{C_y\|V\|_{L^q(\Omega)}^{\gamma_d(q)}}\|u\|_{H^1(\Omega)}^{1-\theta_1}\Big(\|f\|_{L^2(\Omega)}+\|u\|_{L^2(K_y)}\Big)^{\theta_1},
\end{equation*}
where $K_y$ is the image of $\mathcal A_1$. Since $\mathcal A_1$ is a compact subset of $B_0(r_0)\cap\{x_1\geq\rho\}$, $K_y$ is a compact subset of $\Omega$. Hence \eqref{boundary-estimate} holds with $\varsigma_y:=\frac{1}{2}\operatorname{dist}(K_y,\partial\Omega)$, since then $K_y\subset\Omega_{\varsigma_y}$.
\medskip

\noindent\textit{Proof of \eqref{interior-estimate}} Let $\varphi(x_1)= x_1+ \frac{x_1^2}{2}$ be as above, fix $r_0= \frac{1}{4}$, let $\varepsilon_2$ and $\tau_0$ be given by Theorem \ref{Thm2-Carleman-General} with these $r_0$ and $\varphi$, and set 
\begin{equation*}
    \varepsilon=\varepsilon_1, \quad \Phi(x)= \varphi(x_1) - \frac{\varepsilon}{2} |x'|^2, \quad \rho = \frac{\varepsilon r_0^2}{16}, \quad x_*= (\frac{15}{16} \rho, 0) \in \R^d
\end{equation*}
$s_0= \varepsilon r_0$ and $c=\frac{\rho}{32 r_0}$. Note that $c < \frac{1}{2}$ and $\rho \leq \frac{r_0}{64}$. 
\smallskip

\emph{Rescaling.} Let $s\in(0,s_0]$, $x_0\in\Omega$ and $e\in\mathbb S^{d-1}$ with $B_{x_0}(3s)\subset\Omega$. After a translation and a rotation $R$, which replace the coefficients by $R^TA(x_0+R\,\cdot)R$ and preserve \eqref{Ellipticity-A}--\eqref{Lipschitz-condition-A} and the norms in \eqref{interior-estimate}, we may assume that $x_0=0$ and $e=e_1$. Set $\lambda:=\frac s{r_0}\in(0,\varepsilon]$ and
\begin{equation*}
u_\lambda(x):=u\big(\lambda(x-x_*)\big),\qquad x\in B_{2r_0}^+,
\end{equation*}
which is well defined since $\{\lambda(x-x_*), \quad x \in B_{2r_0}^+\}\subset B_0(3s)\subset\Omega$. The function $u_\lambda$ satisfies 
\begin{equation}\label{equation-ulambda}
-\div\big(\tilde A(\varepsilon x)\nabla u_\lambda\big)+\lambda^2V\big(\lambda(x-x_*)\big)u_\lambda=\lambda^2f\big(\lambda(x-x_*)\big)\qquad\text{in }B_{2r_0}^+,
\end{equation}
where $\tilde A(z):=A\big(\lambda(\varepsilon^{-1}z-x_*)\big)$, so that $\tilde A(\varepsilon x)=A\big(\lambda(x-x_*)\big)$. Since $\lambda\leq\varepsilon$, $\tilde A$ satisfies \eqref{Ellipticity-A}--\eqref{Lipschitz-condition-A}, its Lipschitz constant being $\frac\lambda\varepsilon\Lambda_1$. Since $\lambda\leq1$, we have $\lambda^2\|V_p(\lambda(\cdot-x_*))\|_{L^p(B_{2r_0}^+)}\leq\|V_p\|_{L^p(\Omega)}$ for $p\in\{\frac{d}{2},d,\infty\}$, so that the potential is absorbed as in the proof of Lemma \ref{lemma-local-at-the-boundary}.
\smallskip

\emph{Cut-off function and levels of the weight.} We let $\chi_1, \chi_2 \in C_0^\infty (\R; [0,1])$ satisfying  
\begin{equation*}  \chi_1(t) = \left\{
\begin{aligned}
&1 , 
&&  t \geq \varphi( \rho/2),\\
& 0 
&&  t < \varphi(\rho/4),
\end{aligned}
\right.  \qquad 
\chi_2(t) = \left\{
\begin{aligned}
&0 , 
&&  t \geq 2 \rho,\\
& 1 
&&  t < \rho,
\end{aligned}
\right. 
\end{equation*}
and set $\chi(x):=\chi_1(\Phi(x))\chi_2(x_1)$ and $v_\lambda:=\chi u_\lambda$. On $\operatorname{supp}\chi\cap\overline{B_{2r_0}^+}$, we have $\varphi(\frac\rho4)\leq\Phi\leq\varphi(x_1)$ and $x_1\leq2\rho$. Since $\varphi$ is increasing on $[0,2r_0]$, this gives $x_1\geq\frac\rho4$ and $\frac\varepsilon2|x'|^2\leq\varphi(2\rho)\leq3\rho$. Therefore
\begin{equation*}
\operatorname{supp}\chi\cap\overline{B_{2r_0}^+}\subset Q:=\Big\{\frac\rho4\leq x_1\leq2\rho,\ |x'|\leq\frac{\sqrt6}4r_0\Big\}\subset B_0(r_0)\cap B_{x_*}(r_0).
\end{equation*}
In particular, $v_\lambda$ vanishes near $\{x_1=0\}$, so that the boundary terms in Theorem \ref{Thm2-Carleman-General} vanish without any boundary condition. Moreover (see Figure \ref{Geometry far from the boundary-interior observation}): 
\begin{figure}[ht]\centering
\begin{tikzpicture}[xscale=1.35,yscale=1.8,
  lab/.style={font=\scriptsize,inner sep=1pt}]
\fill[green!10] plot[domain=-1:1,samples=60] (\x,{0.5+0.5*\x*\x}) -- cycle;
\fill[orange!30] (-1.2247,1) -- (1.2247,1)
  -- plot[domain=1.2247:1.8708,samples=40] (\x,{0.25+0.5*\x*\x})
  -- (-1.8708,2)
  -- plot[domain=-1.8708:-1.2247,samples=40] (\x,{0.25+0.5*\x*\x}) -- cycle;
\fill[blue!22] plot[domain=-1.8708:1.8708,samples=120] (\x,{0.25+0.5*\x*\x})
  -- (1.7321,2)
  -- plot[domain=1.7321:-1.7321,samples=120] (\x,{0.5+0.5*\x*\x}) -- cycle;
\fill[green!45] (0,0.9375) ellipse[x radius=0.0815,y radius=0.0611];
\draw[orange!80!black,thick,dashed] plot[domain=-2.1679:2.1679,samples=120] (\x,{0.25+0.5*\x*\x});
\draw[blue,thick] plot[domain=-1.8439:1.8439,samples=120] (\x,{0.5+0.5*\x*\x});
\draw[teal,thick] plot[domain=-1.4142:1.4142,samples=120] (\x,{0.75+0.5*\x*\x});
\draw[densely dashed] (-2.4495,0.25) rectangle (2.4495,2);
\draw[very thick] (-2.8,0) -- (2.8,0) node[right]{\small$x'$};
\node[below] at (0,0) {\small$\{x_1=0\}$};
\draw[->] (-2.8,0) -- (-2.8,0.55) node[left]{\small$x_1$};
\draw[gray,dotted,thick] (-2.8,1) -- (2.8,1);
\draw[gray,dotted,thick] (-2.8,2) -- (2.8,2);
\node[left] at (-2.8,1) {\small$x_1=\rho$};
\node[left] at (-2.8,2) {\small$x_1=2\rho$};
\node at (0,1.5) {\small$\mathcal A_1$};
\node[lab] at (0.85,0.736) {$\mathcal A_0$};
\node[below left,lab] at (2.4495,2) {$Q$};
\node[right,lab,text=orange!80!black] at (2.1679,2.6) {$\Phi=\varphi(\frac\rho4)$};
\node[left,lab,text=blue] at (1.8439,2.2) {$\Phi=C_1$};
\draw[teal] (-1.1,1.355) -- (-2.6,1.5) node[left,lab,fill=white,text=teal]{$\Phi=C_2$};
\fill (0,0.9375) circle (0.7pt);
\draw (0,0.9375) -- (0,0.74) node[below,inner sep=1pt]{\scriptsize$x_*$};
\draw (-0.08,0.94) -- (-1.55,0.62) node[left,inner sep=1pt]{\scriptsize$\mathcal W'$};
\end{tikzpicture}
\caption{Geometry of the proof of \eqref{interior-estimate}, in the variables of Theorem \ref{Thm2-Carleman-General}. The cut-off function $\chi=\chi_1(\Phi)\chi_2(x_1)$ is supported in $Q$, equals $1$ in the light green region, which contains $\mathcal W'=B_{x_*}(\frac\rho{16})$ (green), and its gradient is supported in the blue region, contained in $\mathcal A_0$, where $\Phi\leq C_1$, and in the orange region, contained in $\mathcal A_1$, where $\Phi\leq C_3$. On $\mathcal W'$, $\Phi\geq C_2$. Since $Q\subset\{x_1\geq\frac\rho4\}$, no boundary condition is needed on the artificial face $\{x_1=0\}$. In the original variables, $\mathcal W'$ becomes $B_{x_0}(2cs)$ and $\mathcal A_1$ lies in $B_{x_0}(s)\cap\{(z-x_0)\cdot e>cs\}$.}
\label{Geometry far from the boundary-interior observation}
\end{figure}
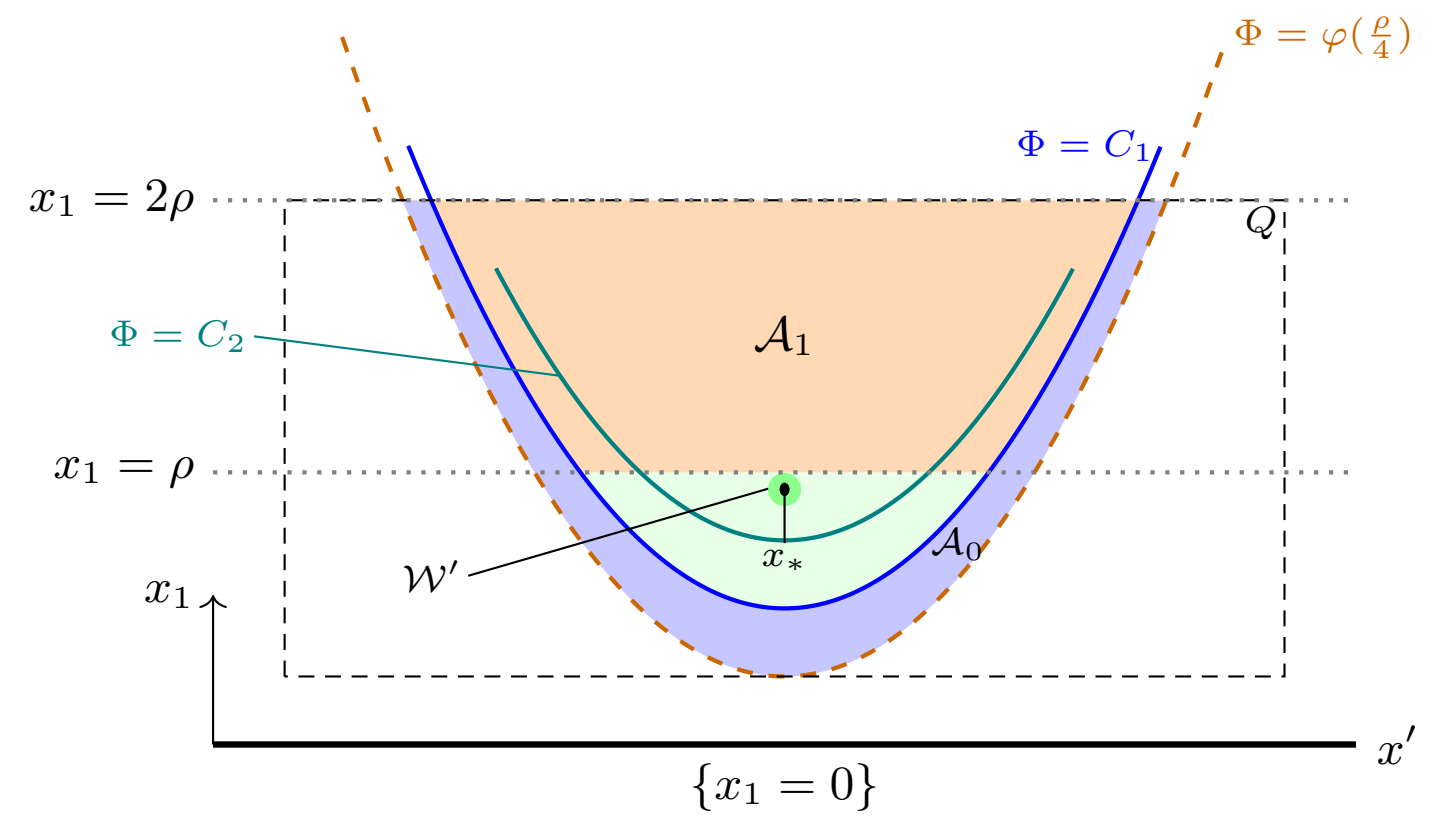

\begin{itemize}
\item in $\overline{B_{2r_0}^+}$, the derivatives of $\chi$ are supported in $\mathcal A_0\cup\mathcal A_1$, where $\mathcal A_0:=Q\cap\{\Phi\leq\varphi(\frac\rho2)\}$ and $\mathcal A_1:=Q\cap\{x_1\geq\rho\}$; on $\mathcal A_0$, we have $\Phi\leq C_1:=\varphi(\frac\rho2)$, and on $Q$, we have $\Phi\leq C_3:=\varphi(2\rho)$;
\item on the ball $\mathcal W':=B_{x_*}(\frac\rho{16})$, we have $\frac78\rho<x_1<\rho$ and $\frac\varepsilon2|x'|^2<\frac\rho8$; hence $\chi=1$ and, since $\varphi'\geq1$, $\Phi\geq\varphi(\frac78\rho)-\frac\rho8\geq\varphi(\frac{3}{4}\rho)=:C_2$.
\end{itemize}

\emph{Carleman estimate.} We now argue as in the proof of \eqref{boundary-estimate}, with $u_\lambda$, $\lambda^2f(\lambda(\cdot-x_*))$ and $\tilde A(\varepsilon\cdot)$ instead of $u_\varepsilon$, $\varepsilon^2f(\varepsilon\cdot)$ and $A_\varepsilon$. Since $C_1< C_2 <C_2$, by optimization there exist $C>0$ and $\theta_2\in (0,1)$ such that 
\begin{equation*}
\|u_\lambda\|_{H^1(\mathcal W')}\leq Ce^{C\|V\|_{L^q(\Omega)}^{\gamma_d(q)}}\|u_\lambda\|_{H^1(B_{2r_0}^+)}^{1-\theta_2}\Big(\lambda^2\|f(\lambda(\cdot-x_*))\|_{L^2(B_{2r_0}^+)}+\|u_\lambda\|_{L^2(\mathcal A_1)}\Big)^{\theta_2}. 
\end{equation*}

\emph{Back to the original variables.} In the variable $z=\lambda(x-x_*)$, the norms are multiplied by powers of $\lambda$, which are absorbed in $C_s$. The image of $\mathcal W'$ is $B_0(\frac{\lambda\rho}{16})=B_0(2cs)$. Since $\mathcal A_1\subset Q\cap\{x_1\geq\rho\}$ and $Q\subset B_{x_*}(r_0)$, the image of $\mathcal A_1$ is contained in $\{|z|<\lambda r_0,\ z_1\geq\frac{\lambda\rho}{16}\}=\{|z|<s,\ z_1\geq2cs\}$, hence in $B_0(s)\cap\{z_1>cs\}$. Undoing the translation and the rotation, we obtain \eqref{interior-estimate}.

\medskip
\noindent\textit{Proof of \eqref{three-ball}.} Let $B_{x_0}(4r)\subset\Omega$, and set $s:=\min\{s_0,cr\}$ and $h:=\frac{cs}2$. We first show that, for every $t\in[r,2r)$,
\begin{equation}\label{three-ball-step}
\|u\|_{H^1(B_{x_0}(t+h))}\leq Ce^{C\|V\|_{L^q(\Omega)}^{\gamma_d(q)}}\|u\|_{H^1(\Omega)}^{1-\theta_2}\Big(\|f\|_{L^2(\Omega)}+\|u\|_{L^2(B_{x_0}(t))}\Big)^{\theta_2}.
\end{equation}
Let $z\in\overline{B_{x_0}(t+h)}$. Set $e_z:=\frac{x_0-z}{|x_0-z|}$ if $z\neq x_0$, and let $e_z\in\mathbb S^{d-1}$ be arbitrary if $z=x_0$. We claim that
\begin{equation}\label{cap-inclusion}
B_z(s)\cap\{w:\ (w-z)\cdot e_z>cs\}\subset B_{x_0}(t).
\end{equation}
Indeed, let $w=z+v$ with $|v|<s$ and $v\cdot e_z>cs$, and set $a:=|z-x_0|\in[0,t+h]$. Since $z-x_0=-a\,e_z$,
\begin{equation*}
|w-x_0|^2=a^2-2a\,v\cdot e_z+|v|^2<a^2-2csa+s^2=:g(a).
\end{equation*}
Since $g$ is convex, $g(a)\leq\max\{g(0),g(t+h)\}$. Moreover, $g(0)=s^2<t^2$ and, since $h=\frac{cs}2$ and $s\leq ct$,
\begin{equation*}
g(t+h)=t^2+s^2-tcs-\tfrac34c^2s^2\leq t^2.
\end{equation*}
Hence $|w-x_0|<t$, which proves \eqref{cap-inclusion}. Since $c<\frac{1}{2}$, we have $|z-x_0|+3s\leq t+h+3s<4r$, hence $B_z(3s)\subset\Omega$. Therefore \eqref{interior-estimate} at $z$, with $e=e_z$, and \eqref{cap-inclusion} give
\begin{equation*}
\|u\|_{H^1(B_z(2cs))}\leq C_se^{C_s\|V\|_{L^q(\Omega)}^{\gamma_d(q)}}\|u\|_{H^1(\Omega)}^{1-\theta_2}\Big(\|f\|_{L^2(\Omega)}+\|u\|_{L^2(B_{x_0}(t))}\Big)^{\theta_2}.
\end{equation*}
Covering the compact set $\overline{B_{x_0}(t+h)}$ by finitely many balls $B_z(2cs)$ with $z\in\overline{B_{x_0}(t+h)}$, whose number depends only on $d$ and $\frac rs$, we obtain \eqref{three-ball-step} (see Figure \ref{Three-ball-step}).
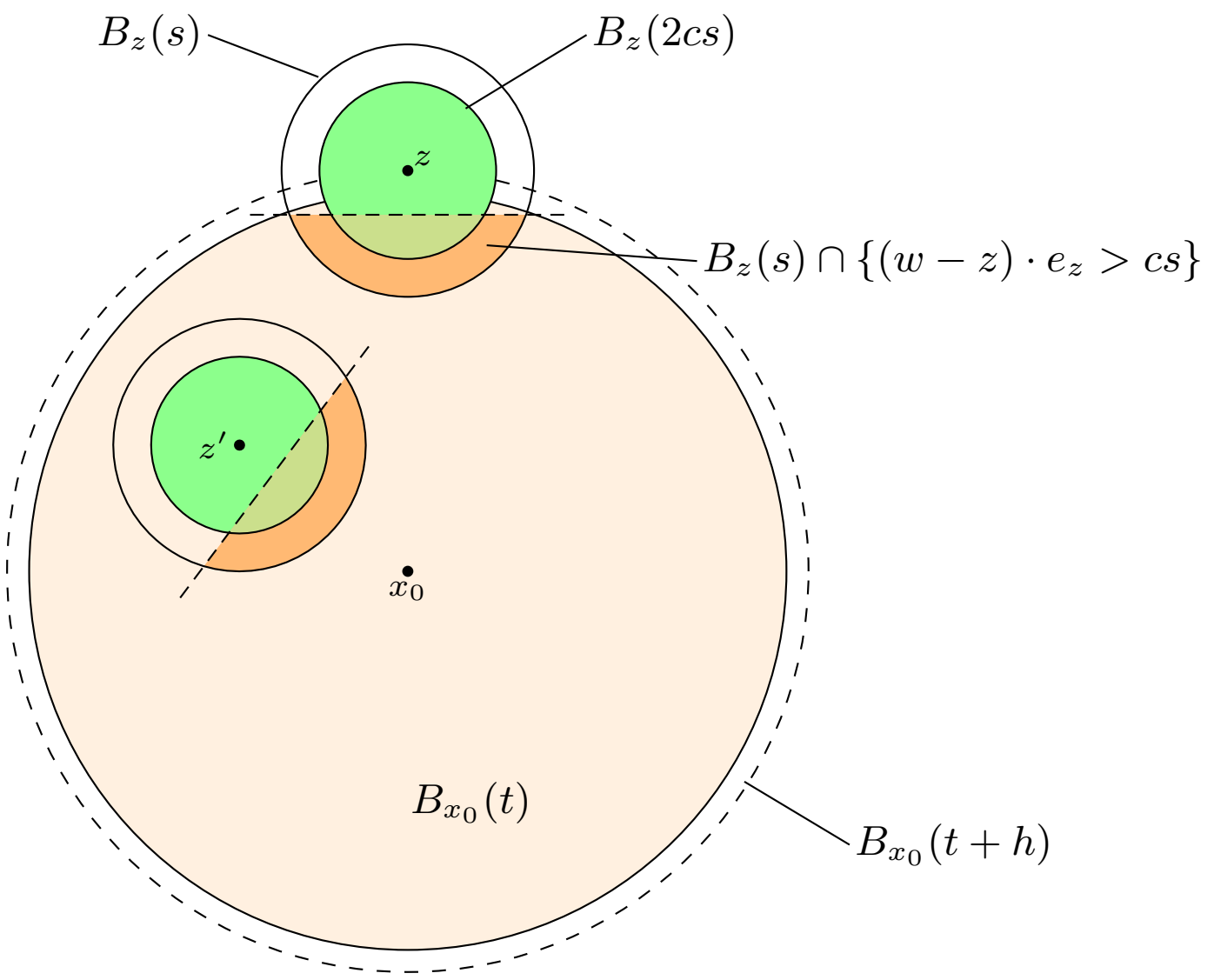
\begin{figure}[ht]\centering
\begin{tikzpicture}[scale=0.8]
\fill[orange!12] (0,0) circle (3.6);
\draw (0,0) circle (3.6);
\draw[dashed] (0,0) circle (3.81);
\begin{scope}[shift={(0,3.81)}]
  \fill[green!45] (0,0) circle (0.84);
  \fill[orange!55] (-1.1241,-0.42) arc[start angle=200.49,end angle=339.51,radius=1.2] -- cycle;
  \fill[orange!55!green!45] (-0.7272,-0.42) arc[start angle=210,end angle=330,radius=0.84] -- cycle;
  \draw (0,0) circle (1.2);
  \draw (0,0) circle (0.84);
  \draw[densely dashed] (-1.5,-0.42) -- (1.5,-0.42);
  \fill (0,0) circle (1.5pt);
  \node[above right,inner sep=1pt] at (0,0) {\scriptsize$z$};
\end{scope}
\begin{scope}[shift={(-1.6,1.2)},rotate=53.13]
  \fill[green!45] (0,0) circle (0.84);
  \fill[orange!55] (-1.1241,-0.42) arc[start angle=200.49,end angle=339.51,radius=1.2] -- cycle;
  \fill[orange!55!green!45] (-0.7272,-0.42) arc[start angle=210,end angle=330,radius=0.84] -- cycle;
  \draw (0,0) circle (1.2);
  \draw (0,0) circle (0.84);
  \draw[densely dashed] (-1.5,-0.42) -- (1.5,-0.42);
  \fill (0,0) circle (1.5pt);
\end{scope}
\node[left,inner sep=2pt] at (-1.6,1.2) {\scriptsize$z'$};
\fill (0,0) circle (1.5pt);
\node[below,inner sep=2pt] at (0,0) {\scriptsize$x_0$};
\node at (0.6,-2.2) {\small$B_{x_0}(t)$};
\draw (3.2,-2.0) -- (4.2,-2.6) node[right,inner sep=1pt]{\small$B_{x_0}(t+h)$};
\draw (-0.85,4.7) -- (-1.9,5.1) node[left,inner sep=1pt]{\small$B_z(s)$};
\draw (0.55,4.4) -- (1.7,5.1) node[right,inner sep=1pt]{\small$B_z(2cs)$};
\draw (0.75,3.1) -- (2.75,2.95) node[right,inner sep=1pt]{\small$B_z(s)\cap\{(w-z)\cdot e_z>cs\}$};
\end{tikzpicture}
\caption{One step of the proof of \eqref{three-ball}. For every $z\in\overline{B_{x_0}(t+h)}$, the cap $B_z(s)\cap\{(w-z)\cdot e_z>cs\}$ (dark orange), beyond the hyperplane at distance $cs$ from $z$ in the direction $e_z$ of $x_0$, lies in $B_{x_0}(t)$, where $u$ is controlled in $L^2$. Then \eqref{interior-estimate} controls $u$ in $H^1$ on $B_z(2cs)$ (green). This is shown for a point $z$ of the sphere of radius $t+h$ and for an interior point $z'$. The balls $B_z(2cs)$, $z\in\overline{B_{x_0}(t+h)}$, cover $\overline{B_{x_0}(t+h)}$.}
\label{Three-ball-step}
\end{figure}

We now iterate \eqref{three-ball-step}. Set $t_j:=r+jh$ and $N_r:=\lceil\frac rh\rceil$, so that $t_j<2r$ for $j<N_r$ and $t_{N_r}\geq2r$. Since $\|f\|_{L^2(\Omega)}+\|u\|_{L^2(B_{x_0}(t_j))}\leq2\|u\|_{H^1(\Omega)}$, we may add $\|f\|_{L^2(\Omega)}$ to the left-hand side of \eqref{three-ball-step}, and we obtain, for $0\leq j<N_r$,
\begin{equation*}
\|f\|_{L^2(\Omega)}+\|u\|_{H^1(B_{x_0}(t_{j+1}))}\leq Ce^{C\|V\|_{L^q(\Omega)}^{\gamma_d(q)}}\|u\|_{H^1(\Omega)}^{1-\theta_2}\Big(\|f\|_{L^2(\Omega)}+\|u\|_{L^2(B_{x_0}(t_j))}\Big)^{\theta_2}.
\end{equation*}
Since $\|u\|_{L^2(B_{x_0}(t_j))}\leq\|u\|_{H^1(B_{x_0}(t_j))}$, the right-hand side of each of these estimates is controlled by the left-hand side of the previous one. By induction, we obtain
\begin{equation*}
\|f\|_{L^2(\Omega)}+\|u\|_{H^1(B_{x_0}(t_{N_r}))}\leq Ce^{C\|V\|_{L^q(\Omega)}^{\gamma_d(q)}}\|u\|_{H^1(\Omega)}^{1-\theta_2^{N_r}}\Big(\|f\|_{L^2(\Omega)}+\|u\|_{L^2(B_{x_0}(r))}\Big)^{\theta_2^{N_r}}.
\end{equation*}
Since $B_{x_0}(2r)\subset B_{x_0}(t_{N_r})$, this gives \eqref{three-ball}, with $\theta_r:=\theta_2^{N_r}$.
\end{proof}

\medskip

We can now conclude the proof of the interpolation inequality.
\begin{proof}[Proof of the interpolation inequality \eqref{interpolation-inequality}]
If $\|f\|_{L^2(\Omega)}>\|u\|_{H^1(\Omega)}$, then $\|u\|_{H^1(\mathcal U)}\leq\|u\|_{H^1(\Omega)}^{1-\alpha}\|f\|_{L^2(\Omega)}^\alpha$ for every $\alpha\in(0,1)$. So we may assume that $\|f\|_{L^2(\Omega)}\leq\|u\|_{H^1(\Omega)}$. Set
\begin{equation*}
\mathcal B(u):=\|f\|_{L^2(\Omega)}+\|u_{|\partial\Omega}\|_{H^{\frac{1}{2}}(\mathcal V\cap\partial\Omega)}+\|\partial_{\nu_A}u\|_{H^{-\frac{1}{2}}(\mathcal V\cap\partial\Omega)}.
\end{equation*}
By Lemma \ref{lemma-local-at-the-boundary}, there exist a neighbourhood $W$ of $y$ in $\R^d$ and $\delta_1\in(0,1)$ such that
\begin{equation*}
\|f\|_{L^2(\Omega)}+\|u\|_{H^1(W\cap\Omega)}\leq Ce^{C\|V\|_{L^q(\Omega)}^{\gamma_d(q)}}\|u\|_{H^1(\Omega)}^{1-\delta_1}\mathcal B(u)^{\delta_1},
\end{equation*}
where we used $\|f\|_{L^2(\Omega)}\leq\|u\|_{H^1(\Omega)}^{1-\delta_1}\mathcal B(u)^{\delta_1}$. Since $u=0$ on $\partial\Omega\cap\mathcal U_1$, Theorem \ref{Thm-interpolation-from-interior} applies with $\omega:=W\cap\Omega$, which is a nonempty open subset of $\Omega$. It gives $\delta_2\in(0,1)$ such that
\begin{equation*}
\|u\|_{H^1(\mathcal U)}\leq Ce^{C\|V\|_{L^q(\Omega)}^{\gamma_d(q)}}\|u\|_{H^1(\Omega)}^{1-\delta_2}\Big(\|f\|_{L^2(\Omega)}+\|u\|_{L^2(W\cap\Omega)}\Big)^{\delta_2}.
\end{equation*}
Combining the two estimates gives \eqref{interpolation-inequality} with $\alpha:=\delta_1\delta_2$.
\end{proof}
	%
\section{Proof of spectral estimates}\label{Sec-Proof-spectral}
The goal of this section is to prove Theorem \ref{thm-spectral-inquality} and \ref{thm-null-controllability}. 
\begin{proof}[Proof of Theorem \ref{thm-spectral-inquality}]
Let $\mu\geq0$ and $w=\sum_{\mu_j\leq\mu}w_j\varphi_j$. For $\lambda\in\R$ and $s\geq0$, we set
\begin{equation*}
F_\lambda(s)=
\begin{cases}
\sinh(s\sqrt\lambda)/\sqrt\lambda, & \lambda>0,\\
s, & \lambda=0,\\
\sin(s\sqrt{-\lambda})/\sqrt{-\lambda}, & \lambda<0,
\end{cases}
\qquad\text{so that}\qquad
F_\lambda''=\lambda F_\lambda,\quad F_\lambda(0)=0,\quad F_\lambda'(0)=1,
\end{equation*}
$u(s,\cdot)=\sum_{\mu_j\leq\mu}w_jF_{\mu_j}(s)\varphi_j $. The sum being finite, $u\in H^2(Z)$, with 
\begin{equation*}
Qu=-\partial_s^2u+Pu=0\ \text{ in }Z,\qquad u|_{(0,S_0)\times\partial\Omega}=0,\qquad u|_{s=0}=0,\qquad\partial_su|_{s=0}=w .
\end{equation*}
Corollary \ref{corollary-interpolations-Q} therefore gives
\begin{equation}\label{Interpolation-u-spectral}
\|u\|_{H^1(Y)}\leq Ce^{C\|V\|_{L^q(\Omega)}^{\gamma_{d+1}(q)}}\|u\|^{1-\delta}_{H^1(Z)}\|w\|^\delta_{L^2(\omega)},
\end{equation}
where $C>0$ and $\delta\in(0,1)$ depend neither on $V$ nor on $\mu$.

We now bound $\|u\|_{H^1(Z)}$ from above. Recall that $\mu_j\geq-R_V$. For $-R_V\leq\lambda\leq\mu$ and $s\in[0,S_0]$, the explicit formulas give $|F_\lambda(s)|\leq S_0e^{S_0\sqrt\mu}$ and $|F_\lambda'(s)|\leq e^{S_0\sqrt\mu}$, so that, by orthonormality,
\begin{equation*}
\|u\|^2_{L^2(Z)}+\|\partial_su\|^2_{L^2(Z)}\leq Ce^{2S_0\sqrt\mu}\|w\|^2_{L^2(\Omega)} .
\end{equation*}
Moreover, by \eqref{quantitative-lower-bound-P} and since $-R_V\leq\mu_j\leq\mu$, for all $s\in[0,S_0]$,
\begin{align*}
\frac{\Lambda_0}2\|\nabla_xu(s)\|^2_{L^2(\Omega)}&\leq\langle Pu(s),u(s)\rangle+R_V\|u(s)\|^2_{L^2(\Omega)}\\
&=\sum_{\mu_j\leq\mu}(\mu_j+R_V)|w_j|^2F_{\mu_j}(s)^2\leq(\mu+R_V)\|u(s)\|^2_{L^2(\Omega)} .
\end{align*}
Hence
\begin{equation}\label{Upper-bound-u-H1-Z}
\|u\|^2_{H^1(Z)}\leq C(1+\mu+R_V)e^{2S_0\sqrt\mu}\|w\|^2_{L^2(\Omega)} .
\end{equation}

We then bound $\|u\|_{L^2(Y)}$ from below. For $\lambda\geq0$, we have $F_\lambda(s)\geq s$ for $s\geq0$. For $\lambda=-k^2<0$, the function $k\mapsto(1+k^2)\int_\beta^{S_0-\beta}k^{-2}\sin^2(ks)\,ds$ is continuous and positive on $(0,+\infty)$, and tends to $\int_\beta^{S_0-\beta}s^2\,ds$ as $k\to0$ and to $\frac{S_0-2\beta}2$ as $k\to+\infty$; it is therefore bounded from below by a constant $c>0$ depending only on $S_0$ and $\beta$. Since $-\mu_j\leq R_V$, the orthonormality of $(\varphi_j)$ gives
\begin{equation}\label{Lower-bound-u-Y}
\frac c{1+R_V}\|w\|^2_{L^2(\Omega)}\leq\|u\|^2_{L^2(Y)}\leq\|u\|^2_{H^1(Y)} .
\end{equation}

Combining \eqref{Interpolation-u-spectral}, \eqref{Upper-bound-u-H1-Z} and \eqref{Lower-bound-u-Y}, we obtain
\begin{equation*}
\Big(\frac c{1+R_V}\Big)^{\frac{1}{2}}\|w\|_{L^2(\Omega)}\leq Ce^{C\|V\|_{L^q(\Omega)}^{\gamma_{d+1}(q)}}\Big(C(1+\mu+R_V)^{\frac{1}{2}}e^{S_0\sqrt\mu}\|w\|_{L^2(\Omega)}\Big)^{1-\delta}\|w\|^\delta_{L^2(\omega)} .
\end{equation*}
If $w=0$, there is nothing to prove. Otherwise, we divide by $\|w\|^{1-\delta}_{L^2(\Omega)}$ and raise to the power $\frac1\delta$. Since $\delta$ depends neither on $V$ nor on $\mu$, since $(1+\mu)^{\frac1{2\delta}}\leq Ce^{\sqrt\mu}$, and since $\log(1+R_V)\leq C\big(1+\|V\|_{L^q(\Omega)}^{\gamma_{d+1}(q)}\big)$, $R_V$ being bounded by a power of $1+\|V\|_{L^q(\Omega)}$, all the factors are bounded by $Ce^{C(\|V\|_{L^q(\Omega)}^{\gamma_{d+1}(q)}+\sqrt\mu)}$, which proves \eqref{Spectral-inequality}.
\end{proof}

\begin{proof}[Proof of Theorem \ref{thm-null-controllability}]
Set $\widetilde P:=P+R_V$, which is self-adjoint and nonnegative by \eqref{quantitative-lower-bound-P}, with eigenfunctions $\varphi_j$ and eigenvalues $\mu_j+R_V\geq0$. For $\mu\geq0$, the range of $\mathbf 1_{(-\infty,\mu]}(\widetilde P)$ is spanned by the $\varphi_j$ with $\mu_j\leq\mu-R_V$, hence contained in $\operatorname{span}\{\varphi_j;\ \mu_j\leq\mu\}$, since $R_V\geq0$. Theorem \ref{thm-spectral-inquality} therefore gives, for all $z\in L^2(\Omega)$,
\begin{equation*}
\big\|\mathbf 1_{(-\infty,\mu]}(\widetilde P)z\big\|_{L^2(\Omega)}\leq d_0e^{d_1\sqrt\mu}\big\|\mathbf 1_{(-\infty,\mu]}(\widetilde P)z\big\|_{L^2(\omega)},\qquad d_0=Ce^{C\|V\|_{L^q(\Omega)}^{\gamma_{d+1}(q)}},\quad d_1=C .
\end{equation*}
By the abstract form of the Lebeau--Robbiano method with explicit constants, see~\cite{EgidiNakicEtAl2020}, this spectral inequality, with exponent $\frac{1}{2}$, implies a final-state observability estimate for the semigroup $(e^{-t\widetilde P})_{t\geq0}$ with the observation operator $z\mapsto z|_\omega$. Bounding its constant in terms of $d_0$, $d_1$ and $T$, and absorbing the powers of $T^{-1}$ into $e^{C/T}$, we obtain, for all $T>0$ and $z_0\in L^2(\Omega)$,
\begin{equation*}
\|e^{-T\widetilde P}z_0\|_{L^2(\Omega)}\leq Ce^{C\left(1+\frac1T+\|V\|_{L^q(\Omega)}^{\gamma_{d+1}(q)}\right)}\|e^{-t\widetilde P}z_0\|_{L^2((0,T)\times\omega)} .
\end{equation*}
By the standard Hilbert Uniqueness Method, this estimate is equivalent to the null-controllability of $\partial_tz+\widetilde Pz=g\mathbf 1_\omega$ in $(0,T)\times\Omega$, $z=0$ on $(0,T)\times\partial\Omega$, $z(0,\cdot)=y_0$, with a control $g$ satisfying
\begin{equation*}
\|g\|_{L^2((0,T)\times\omega)}\leq Ce^{C\left(1+\frac1T+\|V\|_{L^q(\Omega)}^{\gamma_{d+1}(q)}\right)}\|y_0\|_{L^2(\Omega)} .
\end{equation*}
Finally, $y(t)=e^{R_Vt}z(t)$ and $v(t)=e^{R_Vt}g(t)$ satisfy \eqref{Heat-equation} and $y(T,\cdot)=0$, with
\begin{equation*}
\|v\|_{L^2((0,T)\times\omega)}\leq e^{TR_V}\|g\|_{L^2((0,T)\times\omega)},
\end{equation*}
which gives \eqref{control-cost}.
\end{proof}
	%
\section{Local Carleman estimates  for the Laplacian}\label{Sec-Carleman-Laplacian}
The main result of this section is the following local Carleman estimates for the Laplacian. 
\begin{theorem}
\label{Thm-Local-Carl-x-0}
Let $d\geq 3$, $r_0\in(0,\tfrac{1}{2}]$, and let $\Omega=B^+_{2r_0}$ and $ \Gamma$  be as in Theorem \ref{Thm-Carleman-General}. 
Let $\varphi\in C^2([0,2r_0])$ satisfy \eqref{Weight-function-varphi}
and set $ \Phi(x)=\varphi(x_1)-\frac{\varepsilon}{2}\,|x'|^2 .$

There exist $C>0$, $\varepsilon_0\in(0,1]$ and $\tau_0 \geq 1$, depending only on
$d,\alpha_1, \alpha_2, \beta,r_0$ such that for all $0 < \varepsilon \leq \varepsilon_0$, all $\tau\geq \tau_0$, all
$x_0 \in\overline{B^+_{r_0}}$, all $w\in\overline{C_c^\infty}^{\,\Omega}\big(B_{x_0}((\varepsilon\tau)^{-1/2})\cap B_0(r_0)\big)$
 satisfying 
\begin{equation}
    \label{Laplacian-conjugated-w}
    e^{\tau \Phi }\Delta \left(e^{-\tau \Phi }  w\right)   
= f_{2} + f_{2_*'} + \div F \quad\text{in }\Omega\cap B_{x_0}((\varepsilon\tau)^{-1/2}),
\end{equation} 
with 
$ f_{2}\in L^2(\Omega),\quad f_{2_*'}\in L^{\frac{2d}{d+2}}(\Omega),\quad
F\in L^2(\Omega;\C^d)$ vanishing outside $B_{x_0}((\varepsilon\tau)^{-1/2})\cap B_0(\frac{3}{2}r_0)$,  the following estimates hold    
 \begin{multline}
		\label{Carleman-Local-w-x0-L2}
		\tau^{\frac{3}{2}} \| w  \|_{L^2(\Omega)}
		+  
		\tau^{\frac{1}{2}} \| \nabla w  \|_{L^2(\Omega)}
        + 
         \tau^{\frac{3}{4}- \frac{1}{2d}} \| w \|_{L^{\frac{2d}{d-2}}(\Omega)} 
		\leq
		C\Big(
			\| f_{2 } \|_{L^2(\Omega)}
			+ 
			\tau^{\frac{3}{4}- \frac{1}{2d}} \| f_{2_*' } \|_{L^{\frac{2d}{d+2}}(\Omega)} 
			+
			\tau \| F  \|_{L^2(\Omega)} \\ 
			+\tau^{\frac{3}{4}}
\|w\|_{H^{\frac{1}{2}}(\Gamma_{x_0})}+\tau^{\frac{5}{4}}\|w\|_{L^2(\Gamma_{x_0})}
+\tau^{\frac{1}{4}-\theta} \|(\nabla w-F)\cdot\nu_{|\Gamma}\|_{H^{\theta}\Gamma_{x_0})} \Big),
\end{multline}
and
\begin{multline}
		\label{Carleman-Local-w-x0-L2_*} 
		\tau^{\frac{3}{4}+ \frac{1}{2d}} \| w \|_{L^{\frac{2d}{d-2}}(\Omega)} 
         + 
         \tau^{\frac{3}{2}} \| w  \|_{L^2(\Omega)}
		+  
		\tau^{\frac{1}{2}} \| \nabla w  \|_{L^2(\Omega)} 
		\leq
		C\left(
			\| f_{2 } \|_{L^2(\Omega)} 
			+ 
			\tau^{\frac{3}{4}+ \frac{1}{2d}} \| f_{2_*' } \|_{L^{\frac{2d}{d+2}}(\Omega)} 
			+
			\tau \| F  \|_{L^2(\Omega)} \right. \\
+\tau^{\frac{5}{4}}\|w\|_{L^{2}(\Gamma_{x_0})}
+\tau^{\frac{3}{4}+\frac{1}{2d}}\Big(\|w\|_{H^{\frac{1}{2}}(\Gamma_{x_0})}
+\tau^{-\frac{1}{2}-\theta}\|(\nabla w-F)\cdot\nu_{|\Gamma}\|_{H^{\theta}(\Gamma_{x_0})}\Big),
\end{multline}
where  $\Gamma_{x_0}=\Gamma\cap B_{x_0}(\varepsilon\tau)^{-1/2})$ and $\theta \in \{0,-1/2\}$. 
\end{theorem}

\medskip
To prove these estimates, we first establish a Carleman estimate with the weight function
$x\mapsto\varphi(x_1)+\ell\cdot x'$, where $\ell\in\R^{d-1}$ is a vector with $|\ell|$ small. As explained earlier, the feature of this weight, which is linear with respect to $x'$, is that the
conjugated Laplacian
\begin{equation}\label{Conjugated-Operator-strip}
e^{\tau(\varphi+\ell\cdot x')}\,\Delta\,\big(e^{-\tau(\varphi+\ell\cdot x')}\,\cdot\,\big)
=\Delta-2\tau\varphi'(x_1)\,\partial_1-2\tau\,\ell\cdot\nabla'
+\tau^2\varphi'(x_1)^2+\tau^2|\ell|^2-\tau\varphi''(x_1)
\end{equation}
has coefficients independent of the tangential variables, so that it can be factorized after a partial
Fourier transform in $x'$, allowing us to build a parametrix for this operator. The general case is then
recovered by localizing the scale $(\varepsilon\tau)^{-1/2}$ and linearizing the tangential part of the weight at the center of
each ball; see Section \ref{Sec-Proof-Carl-x-0}.
\medskip

Therefore, we focus on the proof of
a Carleman estimate for the Laplacian in the vertical strip
\begin{equation}\label{Def-Strip}
\mathcal S=(0,L)\times\R^{d-1},\qquad 0<L\leq 1,
\end{equation}
with the weight $x\mapsto\varphi(x_1)+\ell\cdot x'$.
	%
	%
\subsection{A Carleman estimate in a strip} \label{Sec-Carleman-Strip}
Throughout this section, $\mathcal{S}$ denotes a strip of the form \eqref{Def-Strip}. The purpose of this subsection is to establish the following Carleman estimate, which extends \cite[Théorème 2.5.10 and Théorème 3.2.1]{Thabouti2024}, corresponding to $\ell=0$, to the presence of the linear tangential term $\ell \cdot x'$ in the weight. 
\begin{proposition}\label{Prop-Carleman-Strip-2} Let $d\geq 3$, let $\mathcal S$ be as in \eqref{Def-Strip} with $0<L\leq1$, let $\varphi\in C^2([0,L])$ satisfy 
\eqref{Weight-function-varphi} on $[0,L]$, and let $\ell\in\R^{d-1}$ with
\begin{equation}\label{Def-ell-star}
|\ell|\leq\ell_*=\frac{\alpha_1}{4} .
\end{equation}
Then there exist $C>0$ and $\tau_0\geq1$, depending only on $d,\alpha_1,\alpha_2,\beta$, and in particular independent of $L \in (0,1]$ and of $\ell$, such that for all $\tau \geq \tau_0$, all 
\begin{equation}\label{Reg-Source-Terms-Strip}
f_2\in L^2(\mathcal S),\quad f_{2_*'}\in L^{\frac{2d}{d+2}}(\mathcal S),\quad
F\in L^2(\mathcal S;\C^d),\quad g\in H^{\frac{1}{2}}(\R^{d-1}),
\end{equation}
and all $w$ solving
\begin{equation}\label{Elliptic-Strip-w}
\left\{
\begin{array}{ll}
e^{\tau(\varphi+\ell\cdot x')}\,\Delta\big(e^{-\tau(\varphi+\ell\cdot x')}w\big)
=f_2+f_{2_*'}+\div F & \text{in }\mathcal S,\\[1mm]
w(0,x')=g(x'), & x'\in\R^{d-1},\\[1mm]
w(L,x')=0, & x'\in\R^{d-1},
\end{array}
\right.
\end{equation}
we have 
\begin{multline}\label{Carleman-Strip-w-1}
\tau^{\frac{3}{2}}\|w\|_{L^2(\mathcal S)}+\tau^{\frac{1}{2}}\|\nabla w\|_{L^2(\mathcal S)}
\leq C\Big(\|f_2\|_{L^2(\mathcal S)}
+\tau^{\frac{3}{4}-\frac{1}{2d}}\|f_{2_*'}\|_{L^{\frac{2d}{d+2}}(\mathcal S)}
+\tau\|F\|_{L^2(\mathcal S)}\\
+\tau^{\frac{3}{4}}\|g\|_{H^{\frac{1}{2}}(\R^{d-1})}
+\tau^{\frac{5}{4}}\|g\|_{L^{2}(\R^{d-1})}
+\tau^{\frac{1}{4}-\theta}\|\left(\partial_1w-F_1\right)(0,\cdot)\|_{H^{\theta}(\R^{d-1})}\Big),
\end{multline}
and
\begin{multline}\label{Carleman-Strip-w-2}
\tau^{\frac{3}{4}+\frac{1}{2d}}\|w\|_{L^{\frac{2d}{d-2}}(\mathcal S)}
\leq C\Big(\|f_2\|_{L^2(\mathcal S)}
+\tau^{\frac{3}{4}+\frac{1}{2d}}\|f_{2_*'}\|_{L^{\frac{2d}{d+2}}(\mathcal S)}
+\tau\|F\|_{L^2(\mathcal S)}\\
+\tau^{\frac{3}{4}+\frac{1}{2d}}\|g\|_{H^{\frac{1}{2}}(\R^{d-1})}
+\tau^{\frac{1}{4}+\frac{1}{2d}-\theta}\|\left(\partial_1w-F_1\right)(0,\cdot)\|_{H^{\theta}(\R^{d-1})}\Big),
\end{multline}
where $\theta \in \{0,-1/2\}$. 
\end{proposition}
\begin{remark}\label{Rem-smallness-ell}
The smallness condition $|\ell|\leq\ell_*<\alpha_1$ keeps the frequencies at which the analysis of Section \ref{Sec-restriction-frequency} degenerates, namely those for which $|\xi'|\lesssim\tau|\ell|$, away from the set where the conjugated operator ceases to be elliptic. The corresponding estimates are then
obtained directly, see \eqref{low-below-turning} and Subsubsection \ref{Sec-low-frequencies}. Conditions of
this type are standard for weights carrying a tangential dependence, see \cite[Section 4E]{LeRousseau2013interface} and \cite{Vessella2017Jumps}.
\end{remark}
The proof of Proposition \ref{Prop-Carleman-Strip-2} consists of the following steps:
\begin{enumerate}
\item We construct an explicit parametrix, which expresses the solution $w$ of \eqref{Elliptic-Strip-w} in terms
of the source terms $(f_2,f_{2_*'},F)$ and of the boundary terms
$\big(w(0,\cdot),(\partial_1w-F_1)(0,\cdot)\big)$, see Subsection \ref{Sec-Construction-parametrix}.
\item We estimate the operators involved in the parametrix. These estimates rest on a version of the the Stein--Tomas restriction theorem, uniform on the family of ellipsoids $\{\operatorname{Re}\psi = c\}$, where $\psi $ is defined in \eqref{Def-Psi-ell} below, and on the corresponding estimates for Fourier multipliers, see Subsection \ref{Sec-restriction-frequency}. Since these ellipsoids degenerate at low frequencies when
$\ell\neq0$, we split the tangential frequencies into the two regions \eqref{Def-regions}. On $\mathcal R_{\mathrm{low}}$, the operators are estimated directly, by a Schur test together with Bernstein's inequality, see Subsubsection \ref{Sec-low-frequencies}. On $\mathcal R_{\mathrm{high}}$, we use the  estimates for Fourier multiplier to bound the interior operators $K_{\tau,j}$ and the boundary operators $G_\tau$ and $H_\tau$, see Subsubsection \ref{Subsubsec-The high-frequency region}. 
\item We combine these estimates, see Subsubsection \ref{Subsubsec-Proof-Prop-Strip}.
\end{enumerate}

	%
	%
\subsection{Construction of the paramertrix}\label{Sec-Construction-parametrix}
We adapt the construction of the parametrix in \cite[Section 2.2.4]{Thabouti2024} (see also \cite[Section 3]{Dehman-Ervedoza-Thabouti-2023}) to the presence of the small linear tangential term $x'\mapsto \ell \cdot x'$ in the Carleman weight. Conjugating by $e^{\tau \ell \cdot x'}$ and taking the partial Fourier transform in $x' \mapsto \xi'$ amounts to the substitution $\xi_j \mapsto \xi_j+\tau \ell_j $, $j\in\{2,\cdots, d\}$, in the tangential symbol. In order to factorize the operators, we define the square root for any $z=a+ib\in\C$ by
\begin{equation}\label{Def-sqrt}
\sqrt{z}=\frac{1}{\sqrt2}\left(\big(\sqrt{a^2+b^2}+a\big)^{\frac{1}{2}}
+i\,\operatorname{sgn}(b)\,\big(\sqrt{a^2+b^2}-a\big)^{\frac{1}{2}}\right), 
\end{equation}
with the convention $\operatorname{sgn}(0)=1$. Then $(\sqrt{z})^2=z$ and  $\operatorname{Re}\big(\sqrt{z}\big)\geq0$ for all $z \in \C.$. We therefore introduce the function $\psi:\R^{d-1}\to\C$ defined by
\begin{equation}\label{Def-Psi-ell}
\psi(\xi')=\Big(\sum_{j=2}^{d}\big(\xi_j+i\tau\ell_j\big)^2\Big)^{\frac{1}{2}}
=\big(|\xi'|^2-\tau^2|\ell|^2+2i\tau\,\ell\cdot\xi'\big)^{\frac{1}{2}}.
\end{equation}
so that, writing
\begin{equation}\label{Psi-ell-R-and-I}
\psi^2=a+ib,\qquad a(\xi')=|\xi'|^2-\tau^2|\ell|^2,\qquad b(\xi')=2\tau\,\ell\cdot\xi',
\end{equation}
we have
\begin{equation}\label{Re-Im-Psi-ell}
\operatorname{Re}\psi=\frac{1}{\sqrt2}\big(\sqrt{a^2+b^2}+a\big)^{\frac{1}{2}}\geq0,
\qquad
\operatorname{Im}\psi=\frac{\operatorname{sgn}(b)}{\sqrt2}\big(\sqrt{a^2+b^2}-a\big)^{\frac{1}{2}} .
\end{equation}
Note that $\psi$ is smooth outside the set $\{\ell \cdot \xi'=0, \, |\xi'|\leq \tau |\ell|\},$ on which $\psi^2 \in (-\infty, 0]$, and that when $\ell=0$, $\psi=|\xi'|$ and the construction below reduces to that of \cite[Section 2.2.4]{Thabouti2024}.

\begin{proposition}\label{Prop-parametrix-Strip-1}
Under the setting of Proposition \ref{Prop-Carleman-Strip-2}, let $\tau\geq 1$, let $(f_2,f_{2_*'},F,g)$ be as in \eqref{Reg-Source-Terms-Strip}, and let $w \in H^1(\mathcal{S)} $ solve \eqref{Elliptic-Strip-w}. Then, formally
\begin{equation}\label{parametrix}
w=K_{\tau,0}\big(f_2+f_{2_*'}\big)+K_{\tau,1}(F_1)+\sum_{j=2}^dK_{\tau,j}(F_j)
+G_\tau(g)+H_\tau\big((\partial_1w-F_1)(0,\cdot)\big),
\end{equation}
where, using the partial Fourier transform, the operators $K_{\tau,j}$, $j\in\{0,\dots,d\}$, are defined for $f$ depending on $(x_1,x')\in\mathcal S$ by
\begin{equation}\label{Def-K-tau-j-Fourier}
\widehat{K_{\tau,j}f}(x_1,\xi')=\int_0^Lk_{\tau,j}(x_1,y_1,\xi')\,\widehat f(y_1,\xi')\,dy_1,
\qquad(x_1,\xi')\in\mathcal S,
\end{equation}
with kernels given, for $(x_1,y_1,\xi')\in[0,L]^2\times\R^{d-1}$, by
\begin{align}
k_{\tau,0}(x_1,y_1,\xi')
&=\mathbf 1_{\operatorname{Re}\psi(\xi')<-\tau\varphi'(x_1)}\,\mathbf 1_{y_1<x_1}\,
e^{\tau(\varphi(x_1)-\varphi(y_1))}\,\frac{\sinh\big(\psi(\xi')(x_1-y_1)\big)}{\psi(\xi')}
\nonumber\\
&\quad-\mathbf 1_{\operatorname{Re}\psi(\xi')\geq-\tau\varphi'(x_1)}\,
e^{\tau(\varphi(x_1)-\varphi(y_1))-\psi(\xi')|x_1-y_1|}
\int_0^{L-\max\{x_1,y_1\}}e^{-2\psi(\xi')\tilde x_1}\,d\tilde x_1,
\label{Def-k-tau-0-x1-y1-xi'}\\[1mm]
k_{\tau,1}(x_1,y_1,\xi')
&=\mathbf 1_{\operatorname{Re}\psi(\xi')<-\tau\varphi'(x_1)}\,\mathbf 1_{y_1<x_1}\,
e^{\tau(\varphi(x_1)-\varphi(y_1))+(x_1-y_1)\psi(\xi')}
\nonumber\\
&\quad-\mathbf 1_{\operatorname{Re}\psi(\xi')\geq-\tau\varphi'(x_1)}\,\mathbf 1_{y_1>x_1}\,
e^{\tau(\varphi(x_1)-\varphi(y_1))+(x_1-y_1)\psi(\xi')} -\big(\psi(\xi')-\tau\varphi'(y_1)\big)\,k_{\tau,0}(x_1,y_1,\xi'),
\label{Def-k-tau-1-x1-y1-xi'}\\[1mm]
k_{\tau,j}(x_1,y_1,\xi')
&=i\,\xi_j\,k_{\tau,0}(x_1,y_1,\xi'),\qquad j\in\{2,\dots,d\}.
\label{Def-k-tau-j-x1-y1-xi'}
\end{align}
The operators $G_\tau$ and $H_\tau$ are given in Fourier, for $g_0\in\mathscr S(\R^{d-1})$, by
\begin{align}
\widehat{G_\tau g_0}(x_1,\xi')&=g_\tau(x_1,\xi')\,\widehat g_0(\xi'),
\qquad\ \ (x_1,\xi')\in\mathcal S,
\label{Def-G-tau-Fourier}\\
\widehat{H_\tau g_0}(x_1,\xi')&=k_{\tau,0}(x_1,0,\xi')\,\widehat g_0(\xi'),
\qquad(x_1,\xi')\in\mathcal S,
\label{Def-H-tau-Fourier}
\end{align}
where $g_\tau$ is given, for $(x_1,\xi')\in[0,L]\times\R^{d-1}$, by
\begin{equation}\label{Def-g-tau}
g_\tau(x_1,\xi')
=\mathbf 1_{\operatorname{Re}\psi(\xi')<-\tau\varphi'(x_1)}\,e^{\tau\varphi(x_1)+x_1\psi(\xi')}
-\big(\tau\varphi'(0)+\psi(\xi')\big)\,k_{\tau,0}(x_1,0,\xi').
\end{equation}
\end{proposition}
\begin{remark}\label{Rem-parametrix-formal}
We emphasize that Proposition \ref{Prop-parametrix-Strip-1} is formal. It will be proved in
Subsections \ref{Sec-low-frequencies} and \ref{Subsubsec-The high-frequency region} that the operators
$K_{\tau,0}$, $(K_{\tau,j})_{j\in\{1,\dots,d\}}$, $G_\tau$ and $H_\tau$ respectively belong to
$\mathcal L\big(L^2(\mathcal S)+L^{\frac{2d}{d+2}}(\mathcal S);L^2(\mathcal S)\big)$,
$\big(\mathcal L(L^2(\mathcal S);L^2(\mathcal S))\big)^d$,
$\mathcal L\big(H^{\frac{1}{2}}(\R^{d-1});L^2(\mathcal S)\big)$ and
$\mathcal L\big(H^{-\frac{1}{2}}(\R^{d-1});L^2(\mathcal S)\big)$, after which a simple density argument justifies \eqref{parametrix} rigorously.
\end{remark}
\begin{proof}
\textit{Step 1: factorization.} Taking the partial Fourier transform $x'\mapsto\xi'$ in
\eqref{Elliptic-Strip-w}$_{(1)}$, and using \eqref{Conjugated-Operator-strip}, the tangential terms combine
into $-|\xi'|^2-2i\tau\,\ell\cdot\xi'+\tau^2|\ell|^2=-\psi(\xi')^2$ and the normal terms into
$(\partial_1-\tau\varphi')^2=\partial_1^2-2\tau\varphi'\partial_1+\tau^2\varphi'^2-\tau\varphi''$, so that
\begin{equation}\label{Conjugated-operator-in-Fourier}
\Big((\partial_1-\tau\varphi'(x_1))^2-\psi(\xi')^2\Big)\widehat w=\widehat H+\partial_1\widehat F_1,
\qquad
\widehat H= \widehat f_2+\widehat f_{2_*'}+\i\sum_{j\geq 2}\xi_j\widehat F_j .
\end{equation}
Since $\psi$ is independent of $x_1$, it commutes with $\partial_1$, and the operator factorizes as
\begin{equation}\label{factorisation}
(\partial_1-\tau\varphi')^2-\psi^2
=\big(\partial_1-\tau\varphi'+\psi\big)\big(\partial_1-\tau\varphi'-\psi\big)
\end{equation}
 Accordingly, introducing the additional unknown $\widehat z=(\partial_1-\tau\varphi'-\psi)\widehat w$, we obtain the
first-order system
\begin{equation}\label{Eq-System-w-z}
\left\{
\begin{array}{ll}
(\partial_1-\tau\varphi'(x_1)-\psi(\xi'))\,\widehat w=\widehat z & \text{in }\mathcal S,\\[1mm]
(\partial_1-\tau\varphi'(x_1)+\psi(\xi'))\,\widehat z=\widehat H+\partial_1\widehat F_1 & \text{in }\mathcal S,\\[1mm]
\widehat w(0,\xi')=\widehat g(\xi') & \text{on }\R^{d-1},\\[1mm]
\widehat w(L,\xi')=0 & \text{on }\R^{d-1}.
\end{array}
\right.
\end{equation}

\textit{Step 2: resolution of \eqref{Eq-System-w-z}$_{(2)}$.} Let $\xi'\in\R^{d-1}$ be fixed. The coefficient of \eqref{Eq-System-w-z}$_{(2)}$ is $\tau\varphi'-\psi$, whose real part satisfies $\operatorname{Re}\big(\tau\varphi'(x_1)-\psi(\xi')\big)\leq-\alpha_1\tau<0$, since $\operatorname{Re}\psi\geq 0$.
We can therefore solve \eqref{Eq-System-w-z}$_{(2)}$ by integrating it forward from $x_1=0$. Working on
$(\widehat z-\widehat F_1)(\cdot,\xi')$, which satisfies
$(\partial_1-\tau\varphi'+\psi)(\widehat z-\widehat F_1)=\widehat H-(\psi-\tau\varphi'(x_1))\widehat F_1$,
Duhamel's formula gives, for $x_1\in(0,L)$ (recall $\varphi(0)=0$),
\begin{multline}\label{Explicit-z}
\widehat z(x_1,\xi')=\widehat F_1(x_1,\xi')
+e^{\tau\varphi(x_1)-\psi(\xi')x_1}
\Big((\partial_1\widehat w-\widehat F_1)(0,\xi')-\big(\psi(\xi')+\tau\varphi'(0)\big)\widehat g(\xi')\Big)\\
+\int_0^{x_1}e^{\tau(\varphi(x_1)-\varphi(y_1))-\psi(\xi')(x_1-y_1)}
\Big(\widehat H(y_1,\xi')-\big(\psi(\xi')-\tau\varphi'(y_1)\big)\widehat F_1(y_1,\xi')\Big)\,dy_1,
\end{multline}
where the boundary value is computed from \eqref{Eq-System-w-z}$_{(1)}$ at $x_1=0$ together with
\eqref{Eq-System-w-z}$_{(3)}$:
$\widehat z(0,\xi')=(\partial_1\widehat w)(0,\xi')-\big(\tau\varphi'(0)+\psi(\xi')\big)\widehat g(\xi')$.

\textit{Step 3: resolution of \eqref{Eq-System-w-z}$_{(1,3,4)}$.} We then focus on the equations
\eqref{Eq-System-w-z}$_{(1,3,4)}$, giving $\widehat w(\cdot,\xi')$ in terms of $\widehat z(\cdot,\xi')$. Since \eqref{Eq-System-w-z}$_{(3,4)}$ provide two boundary conditions for the
first-order equation \eqref{Eq-System-w-z}$_{(1)}$, we have to make a choice when solving it, based on the
requirement that the resulting formulas involve only exponentials with nonpositive real part. To this end,
we analyse the sign of the real part of the coefficient of $\widehat w$ in \eqref{Eq-System-w-z}$_{(1)}$,
that is, of the function $x_1\mapsto\operatorname{Re}\psi(\xi')+\tau\varphi'(x_1)$. Since
$\varphi''\geq \beta>0$, this function is strictly increasing on $[0,L]$, with derivative bounded below by
$\beta\tau$. It therefore vanishes at most once, at some point $x_{\tau,\xi'}\in[0,L]$, and is negative on
$[0,x_{\tau,\xi'})$ and positive on $(x_{\tau,\xi'},L]$, with the conventions $x_{\tau,\xi'}=L$ if it is
negative on all of $[0,L]$, and $x_{\tau,\xi'}=0$ if it is nonnegative at $x_1=0$.

Accordingly, for $x_1\in[0,L]$ such that $\operatorname{Re}\psi(\xi')<-\tau\varphi'(x_1)$, we integrate
forward from $x_1=0$, using \eqref{Eq-System-w-z}$_{(1,3)}$,
\begin{equation}\label{w-forward}
\widehat w(x_1,\xi')=e^{\tau\varphi(x_1)+x_1\psi(\xi')}\,\widehat g(\xi')
+\int_0^{x_1}e^{\tau(\varphi(x_1)-\varphi(y_1))+(x_1-y_1)\psi(\xi')}\,\widehat z(y_1,\xi')\,dy_1,
\end{equation}
while for $x_1\in[0,L]$ such that $\operatorname{Re}\psi(\xi')\geq -\tau\varphi'(x_1)$, we integrate backward
from $x_1=L$, using \eqref{Eq-System-w-z}$_{(1,4)}$,
\begin{equation}\label{w-backward}
\widehat w(x_1,\xi')=-\int_{x_1}^{L}e^{\tau(\varphi(x_1)-\varphi(y_1))+(x_1-y_1)\psi(\xi')}\,
\widehat z(y_1,\xi')\,dy_1 .
\end{equation}
In both cases, the exponentials involved have nonpositive real part, since
$\big|e^{\tau(\varphi(x_1)-\varphi(y_1))+(x_1-y_1)\psi(\xi')}\big|
=e^{\int_{y_1}^{x_1}(\operatorname{Re}\psi(\xi')+\tau\varphi'(r))\,dr}$ and the integrand has a constant
sign on each of the two regions. Combining \eqref{w-forward}--\eqref{w-backward}, we obtain
\begin{multline}\label{Explicit-w}
\widehat w(x_1,\xi')
=\mathbf 1_{\operatorname{Re}\psi(\xi')<-\tau\varphi'(x_1)}
\left(e^{\tau\varphi(x_1)+x_1\psi(\xi')}\,\widehat g(\xi')
+\int_0^{x_1}e^{\tau(\varphi(x_1)-\varphi(y_1))+(x_1-y_1)\psi(\xi')}\,\widehat z(y_1,\xi')\,dy_1\right)\\
-\mathbf 1_{\operatorname{Re}\psi(\xi')\geq-\tau\varphi'(x_1)}
\int_{x_1}^{L}e^{\tau(\varphi(x_1)-\varphi(y_1))+(x_1-y_1)\psi(\xi')}\,\widehat z(y_1,\xi')\,dy_1 .
\end{multline}
The parametrix formula \eqref{parametrix} then follows by inserting \eqref{Explicit-z} into
\eqref{Explicit-w} and collecting the resulting terms. Details are left to the reader.
\end{proof}
	%
	%
 \subsection{The restriction estimates and the frequency decomposition}\label{Sec-restriction-frequency}
 The moduli of the kernels of Proposition \ref{Prop-parametrix-Strip-1} depend on $\xi'$ essentially through $\operatorname{Re}\psi(\xi')$ alone: the exponentials satisfy $|e^{\pm \psi (\xi') (x_1-y_1) }|=|e^{\pm \operatorname{Re} \psi (\xi') (x_1-y_1) }| $, and the branch indicators compare $\operatorname{Re}\psi (\xi')$ with $-\tau \varphi'(x_1)$. We shall therefore use restriction estimates on the level sets of $\xi'\mapsto\operatorname{Re}\psi(\xi')$. For $\lambda>0$, we set
\begin{equation}\label{Def-Sigma-lambda}
\Sigma_\lambda:=\big\{\xi'\in\R^{d-1}:\ \operatorname{Re}\psi(\xi')=\lambda\big\}
=\Big\{\xi'\in\R^{d-1}:\ |\xi'|^2+\frac{\tau^2(\ell\cdot\xi')^2}{\lambda^2}=\lambda^2+\tau^2|\ell|^2\Big\} .
\end{equation}
Indeed, writing $\psi^2=a+ib$ as in \eqref{Psi-ell-R-and-I}, the condition $\operatorname{Re}\psi=\lambda$
reads, by \eqref{Re-Im-Psi-ell}, $\sqrt{a^2+b^2}=2\lambda^2-a$. By squaring, it is equivalent to $b^2=4\lambda^2(\lambda^2-a)$, since the latter gives $2 \lambda^2 - a \geq \lambda >0$, and this is the second equality in \eqref{Def-Sigma-lambda}.

Accordingly, we use the Fourier restriction on the following family of ellipsoids: for $s\geq0$ and a unit vector $e\in\mathbb S^{d-2}$, writing $\eta=u\,e+\eta_\perp$ with $u=e\cdot\eta$ and $\eta_\perp\perp e$, we set
\begin{equation}\label{Sigma-ellipsoid}
\Sigma_{s,e}:=\Big\{\eta\in\R^{d-1}:\ u^2+\frac{|\eta_\perp|^2}{1+s^2}=1\Big\},
\end{equation}
that is, $\Sigma_{s,e}$ is the ellipsoid of revolution about the axis $\R e$, with semi-axis $1$ in the direction $e$ and semi-axes $\sqrt{1+s^2}$ in the $d-2$ orthogonal directions. Setting $\eta =\xi'/\lambda$ and dividing  \eqref{Def-Sigma-lambda} by $\lambda^2(1+s^2)$, we obtain  
\begin{equation}\label{Sigma-lambda-dilation}
\Sigma_\lambda=\lambda\,\Sigma_{s,e_\ell},
\qquad\text{with}\quad s=\frac{\tau|\ell|}{\lambda}, \qquad e_\ell=\frac{\ell}{|\ell|} \text{ if } \ell \neq 0, 
\end{equation}
and $e_\ell $ any unit vector if $\ell =0$. After a rotation mapping $e$ to the coordinate vector $e_2$, the ellipsoid $\Sigma_{s,e}$ is of the form $\big\{\sum_{j=2}^{d}\mu_j\eta_j^2=1\big\}$ with $\mu_2=1,\qquad \mu_3=\cdots=\mu_{d}=\frac{1}{1+s^2}.$ Hence, for all $s_*>0$, $s\in [0,s_*]$ and $\eta \in \R^{d-1}$, 
\begin{equation}\label{Coercivity-Restriction}
\frac{1}{c_0}|\eta|^2\leq\sum_{j=2}^{d}\mu_j\eta_j^2\leq c_0|\eta|^2
\quad\text{with}\quad c_0=1+s^2 .
\end{equation}
This is precisely the setting of \cite[Theorem 5.2]{Dehman-Ervedoza-Thabouti-2023}, with $n=d-1$ and the coefficients indexed by $j\in \{2,\cdots,d\}$. Indeed, the coefficients in \eqref{Coercivity-Restriction} are of the form $\mu_j=1-a\lambda_j$, with $\lambda_2=0$, $\lambda_3=\cdots=\lambda_d=1$ and $a=\frac{s^2}{1+s^2 } $, so that $s \in [0,s_*]$ corresponds to $a \in [0,a_*]$, with $a_*=\frac{s_*^2}{1+s_*^2 }<1$ and \eqref{Coercivity-Restriction} is the coercivity condition \cite[(5.3)]{Dehman-Ervedoza-Thabouti-2023} on $[X_0,X_1]= [0,a_*]$. Since the constant in
\cite[Theorem 5.2]{Dehman-Ervedoza-Thabouti-2023} depends only on the dimension and on the coercivity constant, and since all the norms involved are invariant under rotations, we obtain the following result.
\begin{theorem}\label{Thm-Fourier-Rest-Sigma-s}
Let $d\geq3$ and $s_*>0$. There exists a constant $C>0$, depending only on $d$ and $s_*$, such that for all $s\in[0,s_*]$, all $e\in\mathbb S^{d-2}$, and all $f\in L^{\frac{2d}{d+2}}(\R^{d-1})$,
\begin{equation}\label{Uniform-Restriction-Thm}
\big\|\widehat f\big\|_{L^2(\Sigma_{s,e})}\leq C\,\|f\|_{L^{\frac{2d}{d+2}}(\R^{d-1})} .
\end{equation}
\end{theorem}
The uniformity in $s$ is the point of this Theorem: in \eqref{Uniform-Restriction-Thm}, the parameter $s =\tau |\ell |/\lambda$ depends on the level $\lambda$, whereas the constant in \eqref{Uniform-Restriction-Thm} does not. Undoing the dilation \eqref{Uniform-Restriction-Thm} gives the following form, on the level set of $\operatorname{Re} \psi $ themselves. 
\begin{corollary} \label{Cor-Restriction-Sigma-lambda}
Let $s_*>0$ and let $\psi$ be as in \eqref{Def-Psi-ell}. There exists $C>0$, depending only on $d$ and
$s_*$, such that for all $\tau\geq1$, all $\ell\in\R^{d-1}$, all $\lambda\geq\tau|\ell|/s_*$, and all $f\in L^{\frac{2d}{d+2}}(\R^{d-1})$,
\begin{equation}
\label{Uniform-Restriction-lambda}
\|\widehat f\|_{L^2(\Sigma_\lambda)}
\leq C\,\lambda^{\frac{d-2}{2d}}\,\|f\|_{L^{\frac{2d}{d+2}}(\R^{d-1})} .
\end{equation}
\end{corollary}

\begin{proof}
By \eqref{Sigma-lambda-dilation}, $\Sigma_\lambda=\lambda\,\Sigma_{s,e_\ell}$ with
$s=\tau|\ell|/\lambda\leq s_*$. Setting $g(x)=f(x/\lambda)$, so that
$\widehat g(\eta)=\lambda^{d-1}\widehat f(\lambda\eta)$ and $\|g\|_{L^p}=\lambda^{(d-1)/p}\|f\|_{L^p}$, with
$p=\tfrac{2d}{d+2}$, and
using that the surface measure of $\Sigma_\lambda=\lambda\,\Sigma_{s,e_\ell}$ satisfies
$dS_{\Sigma_\lambda}(\lambda\eta)=\lambda^{d-2}\,dS_{s,e_\ell}(\eta)$, we obtain
\begin{equation*}
    \|\widehat f\|_{L^2(\Sigma_\lambda)}^2
=\lambda^{d-2}\int_{\Sigma_{s,e_\ell}}\big|\widehat f(\lambda\eta)\big|^2\,dS_{s,e_\ell}(\eta)
=\lambda^{-d}\,\|\widehat g\|_{L^2(\Sigma_{s,e_\ell})}^2 .
\end{equation*}
Applying Theorem \ref{Thm-Fourier-Rest-Sigma-s} to $g$ yields
$\|\widehat f\|_{L^2(\Sigma_\lambda)}\leq C\lambda^{-\frac{d}{2}+\frac{d-1}{p}}\|f\|_{L^p}$, which is \eqref{Uniform-Restriction-lambda}, since $-\frac{d}{2}+\frac{d-1}{p}= \frac{d-2}{2d}$.
\end{proof}
\medskip
The condition $s\leq s_*$ is related to the curvature of $\Sigma_{s,e}$. Near its poles $\pm e$, $\Sigma_{s,e}$ is the graph of $\eta_\perp\mapsto\pm\big(1-|\eta_\perp|^2/(1+s^2)\big)^{1/2}$, so that its principal curvatures at $\pm e$ all equal $(1+s^2)^{-1}$, and its Gaussian curvature takes there the value $(1+s^2)^{-(d-2)}$, which is its minimal value on $\Sigma_{s,e}$. Hence, for $s\in[0,s_*]$, the Gaussian curvature of $\Sigma_{s,e}$ is bounded from below by $(1+s_*^2)^{-(d-2)}=c_0^{-(d-2)}$, whereas, as $s\to+\infty$, $\Sigma_{s,e}$ flattens near its poles. Since restriction estimates of Stein--Tomas type rely on the nonvanishing of the Gaussian curvature, and fail on flat hypersurfaces, the constant in \eqref{Uniform-Restriction-Thm} cannot be expected to be uniform as $s\to+\infty$. This is why the restriction estimates will only be used on the region $\{\operatorname{Re}\psi\geq2\tau|\ell|\}$, on which $s\leq\frac{1}{2}$, see \eqref{Def-regions} below.

we shall also need the following elementary bounds.  
\begin{lemma}\label{Lem-psi-bounds}
For all $\xi'\in\R^{d-1}$,
\begin{equation}\label{psi-bounds}
 (\operatorname{Re}\psi(\xi'))^2\ \leq\ |\xi'|^2\ \leq\ (\operatorname{Re}\psi(\xi'))^2+\tau^2|\ell|^2, \qquad \operatorname{Im}\psi(\xi')\ \leq \tau|\ell|, 
\end{equation}
so that $\operatorname{Re}\psi(\xi')\ \leq\ |\xi'|\ \leq\ \operatorname{Re}\psi(\xi')+\tau|\ell|,
$ and $
\operatorname{Re}\psi(\xi')\ \leq\ |\psi(\xi')|\ \leq\ \operatorname{Re}\psi(\xi')+\tau|\ell|$. Moreover, $\psi$ is smooth on $\{\operatorname{Re}\psi>0\}$ and for all $\lambda>0$ and all $\xi'\in\Sigma_\lambda$, with $s=\tau|\ell|/\lambda$,
\begin{equation}\label{gradient-bounds}
\frac{1}{\sqrt{1+s^2}}\ \leq\ \big|\nabla_{\xi'}\!\operatorname{Re}\psi(\xi')\big|\ \leq\ \sqrt{1+2 s^2} .
\end{equation}
\end{lemma}
\begin{proof}
\textit{Proof of \eqref{psi-bounds}.} With $a=|\xi'|^2- \tau^2 |\ell|^2$ as in \eqref{Psi-ell-R-and-I}, we have $(\operatorname{Re}\psi)^2= \frac{1}{2} (|\psi|^2+a)$ and $(\operatorname{Im}\psi)^2= \frac{1}{2} (|\psi|^2-a)$, while $ |\psi|^2 \geq a$  and 
\begin{equation*}
  |\psi|^2 =  |   \sum_{j=2}^{d}(\xi_j+i \tau \ell_j )|^2  \leq    |\xi'|^2+\tau^2|\ell|^2. 
\end{equation*}
Hence $a\leq (\operatorname{Re}\psi)^2 \leq |\xi'|^2$ and $ (\operatorname{Im}\psi)^2 \leq \tau^2|\ell|^2,$ which is \eqref{psi-bounds}. The other bounds follow, using $\sqrt{x+y} \leq \sqrt{x} + \sqrt{y}$ and $|\psi| \leq \operatorname{Re}\psi+ |\operatorname{Im}\psi|$. 
\smallskip

\textit{Proof of \eqref{gradient-bounds}.} On $\{\operatorname{Re}\psi>0\}$, we have $\psi^2\notin(-\infty,0]$, so that $\psi$ is smooth there, and differentiating $\psi^2=|\xi'|^2-\tau^2|\ell|^2+2i\tau\,\ell\cdot\xi'$ gives $\nabla_{\xi'}\psi=(\xi'+i\tau\ell)/\psi$. Let $\lambda>0$ and $\xi'\in\Sigma_\lambda$, and write $\psi=\lambda+i\iota$. By \eqref{psi-bounds}, $|\xi'|^2\leq\lambda^2+\tau^2|\ell|^2$ and $\iota^2\leq\tau^2|\ell|^2$. The upper bound follows from
\begin{equation*}
\big|\nabla_{\xi'}\operatorname{Re}\psi\big|\leq\big|\nabla_{\xi'}\psi\big|=\frac{\big(|\xi'|^2+\tau^2|\ell|^2\big)^{1/2}}{|\psi|}\leq\frac{\big(\lambda^2+2\tau^2|\ell|^2\big)^{1/2}}{\lambda}.
\end{equation*}
For the lower bound, the imaginary and real parts of $\psi^2=(\lambda+i\iota)^2$ give $\tau\,\ell\cdot\xi'=\lambda\iota$ and $|\xi'|^2+\iota^2=\lambda^2+\tau^2|\ell|^2$, so that
\begin{equation*}
\xi'\cdot\nabla_{\xi'}\operatorname{Re}\psi
=\operatorname{Re}\frac{\big(|\xi'|^2+i\tau\,\ell\cdot\xi'\big)(\lambda-i\iota)}{|\psi|^2}
=\frac{\lambda\big(|\xi'|^2+\iota^2\big)}{\lambda^2+\iota^2}
=\frac{\lambda\big(\lambda^2+\tau^2|\ell|^2\big)}{\lambda^2+\iota^2}\geq\lambda .
\end{equation*}
Hence $\big|\nabla_{\xi'}\operatorname{Re}\psi\big|\geq\lambda/|\xi'|\geq\lambda/(\lambda^2+\tau^2|\ell|^2)^{1/2}$, which is the lower bound in \eqref{gradient-bounds}.
\end{proof}

Now, for $g\in\mathscr S(\R^{d-1})$, we decompose $g$ along the level sets of $\operatorname{Re}\psi$. Since $\operatorname{Re}\psi$ is smooth with nonvanishing gradient on $\{\operatorname{Re}\psi>0\}$, by Lemma \ref{Lem-psi-bounds} and since $\{\operatorname{Re} \psi =0\}= \{\ell \cdot \xi'=0, \, |\xi'| \leq \tau |\ell|\}$ is a null set, the coarea formula gives 
\begin{multline}\label{Def-E-lambda}
g=\int_0^{\infty}E_\lambda[g]\,d\lambda,
\quad  \text{where} \quad 
E_\lambda[g](x')=\mathcal F^{-1}\Big(\delta\big(\operatorname{Re}\psi(\xi')-\lambda\big)\,\widehat g(\xi')\Big)(x')
\\ =\frac{1}{(2\pi)^{d-1}}\int_{\Sigma_\lambda}e^{ix'\cdot\xi'}\,\widehat g(\xi')\,
\frac{dS_\lambda(\xi')}{\big|\nabla_{\xi'}\!\operatorname{Re}\psi(\xi')\big|}, 
\end{multline}
with $dS_\lambda(\xi')$ is the surface measure on $\Sigma_\lambda$.  When $\ell=0$, \eqref{Def-E-lambda} is the spherical decomposition of
\cite[Section 2.5.4]{Thabouti2024}. Denoting by $R_\lambda g=\widehat g\,\big|_{\Sigma_\lambda}$ the restriction operator and by
$R_\lambda^\star$ its adjoint, we have $E_\lambda= R_\lambda^\star\circ\big(|\nabla_{\xi'}\!\operatorname{Re}\psi|^{-1}\,\cdot\,\big)\circ R_\lambda$. Hence, Corollary \ref{Cor-Restriction-Sigma-lambda} with $s_*=\frac{1}{2},$ its dual and  \eqref{gradient-bounds} imply, for all $\lambda >0$ with $\lambda \geq 2 \tau |\ell |$, 
\begin{equation}\label{E-lambda-bound}
\big\|E_\lambda[g]\big\|_{L^{\frac{2d}{d-2}}(\R^{d-1})}
\leq C\,\lambda^{\frac{d-2}{d}}\,\|g\|_{L^{\frac{2d}{d+2}}(\R^{d-1})}, 
\end{equation}
where $C$ depending only on $d$. Formula \eqref{Def-E-lambda}, together with \eqref{Uniform-Restriction-lambda}--\eqref{E-lambda-bound}, is the tool by which the operators of Proposition \ref{Prop-parametrix-Strip-1} will be estimated on the high frequencies.

 Accordingly, we split the tangential frequencies into tow regions
\begin{equation}\label{Def-regions}
\mathcal R_{\mathrm{high}}=\big\{\xi'\in\R^{d-1}:\ \operatorname{Re}\psi(\xi')\geq2\tau|\ell|\big\},
\qquad
\mathcal R_{\mathrm{low}}=\big\{\xi'\in\R^{d-1}:\ \operatorname{Re}\psi(\xi')<2\tau|\ell|\big\},
\end{equation}
and we denote by $\Pi_{\mathrm{high}}$ and $\Pi_{\mathrm{low}}$ the corresponding Fourier projections, $\widehat{\Pi_{\mathrm{high}}g}:=\mathbf 1_{\mathcal R_{\mathrm{high}}}\widehat g$ and $\widehat{\Pi_{\mathrm{low}}g}:=\mathbf 1_{\mathcal R_{\mathrm{low}}}\widehat g$. Note that $\mathcal R_{\mathrm{low}}=\emptyset$ when $\ell=0$. On $\mathcal R_{\mathrm{high}}$, the level parameter satisfies $s=\tau|\ell|/\operatorname{Re}\psi\leq\frac{1}{2}$, so that Corollary \ref{Cor-Restriction-Sigma-lambda} applies with $s_*=\frac{1}{2}$, \eqref{E-lambda-bound} holds and, by \eqref{gradient-bounds},
\begin{equation}\label{gradient-high}
\frac{2}{\sqrt5}\leq\big|\nabla_{\xi'}\operatorname{Re}\psi\big|\leq\sqrt{\frac{3}{2}}
\qquad\text{on }\Sigma_\lambda,\quad\lambda\geq2\tau|\ell| .
\end{equation}
On $\mathcal R_{\mathrm{low}}$, by \eqref{psi-bounds} and $|\ell|\leq\ell_*=\frac{\alpha_1}{4}$, we have, for all $x_1\in[0,L]$,
\begin{equation}\label{low-below-turning}
|\xi'|<3\tau|\ell|\leq\frac{3\alpha_1}{4}\,\tau
\qquad\text{and}\qquad
\operatorname{Re}\psi(\xi')<2\tau|\ell|\leq\frac{\alpha_1}{2}\,\tau\leq-\tau\varphi'(x_1)-\frac{\alpha_1}{2}\,\tau .
\end{equation}
Thus, on $\mathcal R_{\mathrm{low}}$, all the frequencies lie below the turning point, for all $x_1\in[0,L]$: only the first terms of the kernels of Proposition \ref{Prop-parametrix-Strip-1} occur, they involve only exponentials decaying in $x_1-y_1\geq0$ at least at the rate $\frac{\alpha_1}{2}\tau$, and $|\xi'|\lesssim\tau$. The corresponding operators are then estimated directly, by a Schur test together with Bernstein's type inequality, see Proposition \ref{Prop-low-frequencies}. 

	%
    \subsection{$L^p-L^q$ estimates for Fourier multipliers} \label{Subsec-Fourier-Multipliers}
    In this subsection, we prove the $L^p-L^q$ estimates for Fourier multipliers by which the operators of Proposition \ref{Prop-parametrix-Strip-1} will be estimated, following \cite[Section 5.2]{Dehman-Ervedoza-Thabouti-2023}.
    
    Recall from \eqref{Sigma-lambda-dilation} that, for $\lambda>0$,
$\Sigma_\lambda=\lambda\,\Sigma_{s,e_\ell}$ with $s=\tau|\ell|/\lambda$. Accordingly, we parametrize  $\{\operatorname{Re}\psi >0\}$, whose complement in $\R^{d-1}$ is a null set, by the map
\begin{equation}\label{Def-polar-coordinates}
\xi'\ \longmapsto\ \big(\lambda,\omega'\big):=\Big(\operatorname{Re}\psi(\xi'),\
\frac{\xi'}{\operatorname{Re}\psi(\xi')}\Big), \qquad \lambda \in\ \R_+^*, \quad \omega ' \in \Sigma_{s,e_\ell},
\quad s=\frac{\tau|\ell|}{\lambda},
\end{equation}
whose inverse is $(\lambda, \omega ') \mapsto \lambda \omega'$, and in which $\lambda$ plays the role of a radius and $\omega'$ that of an angular variable. For $\ell=0$, \eqref{Def-polar-coordinates} coincides with the classical polar coordinates of $\R^{d-1}$. Contrary to the
situation of \cite[Section 5.2]{Dehman-Ervedoza-Thabouti-2023}, where the function $\psi(a,\cdot)$ is
homogeneous of degree one and the surface $\Sigma_a$ is therefore independent of the radius, the surface
$\Sigma_{s,e_\ell}$ in \eqref{Def-polar-coordinates} depends on $\lambda$ through the parameter $s$. This
is compensated by the fact that the restriction estimate of Theorem \ref{Thm-Fourier-Rest-Sigma-s} is
uniform with respect to $s\in[0,s_*]$. 

For a measurable function $k$ defined on
$\big\{(\lambda,\omega'):\, \lambda>0, \, \lambda\geq2\tau|\ell|,\, \omega'\in\Sigma_{\tau|\ell|/\lambda,\,e_\ell}\big\}$, we
consider the Fourier multiplier operator $K_k$ defined, for $f\in \mathcal{S}(\R^{d-1})$, by 
\begin{equation}\label{Def-Fourier-mulptilier-K}
\widehat{K_kf}(\xi')=\mathbf 1_{\mathcal R_{\mathrm{high}}}(\xi')\,
k\Big(\operatorname{Re}\psi(\xi'),\frac{\xi'}{\operatorname{Re}\psi(\xi')}\Big)\,\widehat f(\xi'),
\qquad \xi'\in\R^{d-1},
\end{equation}
and we write,  
\begin{equation*}
  \big\|k(\lambda,\cdot)\big\|_{L^\infty}:=\big\|k(\lambda,\cdot)\big\|_{L^\infty(\Sigma_{\tau|\ell|/\lambda,\,e_\ell})} ,  
\end{equation*}
which is the supremum of the multiplier of $K_k$ on the level set $\Sigma_\lambda$.

\begin{proposition}  \label{Lemma-E-lambda-multiplier}
Let $d\geq3$. There exists a constant $C>0$, depending only on $d$, such that for all $\tau\geq1$, all
$\ell\in\R^{d-1}$ and all $k$ as above:
\begin{itemize}
\item if $k$ is bounded, that is, if $\sup_{\lambda\geq 2 \tau |\ell|} \|k(\lambda,\cdot)\|_{L^\infty} <\infty$, then $K_k$ extends as a bounded operator on $L^2(\R^{d-1})$, and
\begin{equation}\label{Est-K-k-L2}
\big\|K_k\big\|_{\mathcal L(L^2(\R^{d-1}))}\ \leq\
\sup_{\lambda\geq2\tau|\ell|}\big\|k(\lambda,\cdot)\big\|_{L^\infty},
\end{equation}
\item if $\displaystyle\int_{2\tau|\ell|}^{\infty}\big\|k(\lambda,\cdot)\big\|_{L^\infty}\,
\lambda^{1-\frac{2}{d}}\,d\lambda<\infty$, then
$K_k$ extends as a bounded operator from  $L^{\frac{2d}{d+2}}(\R^{d-1})$ to $ L^{\frac{2d}{d-2}}(\R^{d-1})$,  and
\begin{equation}\label{Est-K-k-Lp-Lp'}
\big\|K_k\big\|_{\mathcal L\left(L^{\frac{2d}{d+2}},\,L^{\frac{2d}{d-2}}\right)}
\ \leq\ C\int_{2\tau|\ell|}^{\infty}\big\|k(\lambda,\cdot)\big\|_{L^\infty}\,\lambda^{1-\frac{2}{d}}\,d\lambda ,
\end{equation}
\item if $\displaystyle\int_{2\tau|\ell|}^{\infty}\big\|k(\lambda,\cdot)\big\|^2_{L^\infty}\,
\lambda^{1-\frac{2}{d}}\,d\lambda<\infty$, then
$K_k$ extends as a bounded operator from  $L^{\frac{2d}{d+2}}(\R^{d-1})$ to $L^2(\R^{d-1})$ and from $L^2(\R^{d-1})$ to $ L^{\frac{2d}{d-2}}(\R^{d-1})$,  and 
\begin{equation}\label{Est-K-k-mixed}
\big\|K_k\big\|_{\mathcal L\left(L^{\frac{2d}{d+2}},\,L^{2}\right)}
+\big\|K_k\big\|_{\mathcal L\left(L^{2},\,L^{\frac{2d}{d-2}}\right)}
\ \leq\ C\left(\int_{2\tau|\ell|}^{\infty}\big\|k(\lambda,\cdot)\big\|^2_{L^\infty}\,
\lambda^{1-\frac{2}{d}}\,d\lambda\right)^{\frac{1}{2}} .
\end{equation}
\end{itemize}
\end{proposition}

\begin{proof}
Denote by $R_\lambda g:=\widehat g\,\big|_{\Sigma_\lambda}$ the restriction operator, by $R_\lambda^*$ its
adjoint, and by $w_\lambda:=\big|\nabla_{\xi'}\!\operatorname{Re}\psi\big|^{-1}$ on $\Sigma_\lambda$, so
that \eqref{Def-E-lambda} reads $E_\lambda[g]= (2\pi) ^{-(d-1)}R_\lambda^*\big[w_\lambda R_\lambda g\big]$.  Since, Since $\mathbf 1_{\mathcal R_{\mathrm{high}}}=1$ on $\Sigma_\lambda$ if $\lambda\geq2\tau|\ell|$, and $\mathbf 1_{\mathcal R_{\mathrm{high}}}=0$ on $\Sigma_\lambda$ otherwise, \eqref{Def-E-lambda} and \eqref{Def-Fourier-mulptilier-K} give, for $f\in\mathscr S(\R^{d-1})$,
\begin{equation}\label{K-k-decomposition}
K_kf=\int_{2\tau|\ell|}^{\infty}E_\lambda\big[K_kf\big]\,d\lambda
=(2\pi)^{1-d}\int_{2\tau|\ell|}^{\infty}R_\lambda^\star\big[w_\lambda\,k_\lambda\,R_\lambda f\big]\,d\lambda,
\end{equation}
where $k_\lambda(\xi'):=k(\lambda,\xi'/\lambda)$ for $\xi'\in\Sigma_\lambda$, so that $\|k_\lambda\|_{L^\infty(\Sigma_\lambda)}=\|k(\lambda,\cdot)\|_{L^\infty}$. By \eqref{gradient-high}, we have $\frac{1}{2}\leq w_\lambda\leq2$ on $\Sigma_\lambda$ for $\lambda\geq2\tau|\ell|$. We shall use two consequences of Theorem \ref{Thm-Fourier-Rest-Sigma-s}, both uniform in $\lambda\geq2\tau|\ell|$:
\begin{equation}\label{R-Rstar-bounds}
\big\|R_\lambda\big\|_{\mathcal L\left(L^{\frac{2d}{d+2}}(\R^{d-1}),\,L^2(\Sigma_\lambda)\right)}
\leq C\lambda^{\frac{d-2}{2d}},
\qquad
\big\|R_\lambda^\star\big\|_{\mathcal L\left(L^2(\Sigma_\lambda),\,L^{\frac{2d}{d-2}}(\R^{d-1})\right)}
\leq C\lambda^{\frac{d-2}{2d}},
\end{equation}
the first one being Corollary \ref{Cor-Restriction-Sigma-lambda} with $s_*=\frac{1}{2}$ and the second one its dual, and
\begin{equation}\label{Plancherel-layers}
\int_{2\tau|\ell|}^{\infty}\big\|R_\lambda g\big\|^2_{L^2(\Sigma_\lambda)}\,d\lambda
\leq2\int_0^{\infty}\!\!\int_{\Sigma_\lambda}\big|\widehat g\big|^2\,w_\lambda\,dS_\lambda\,d\lambda
=2\,\big\|\widehat g\big\|^2_{L^2(\R^{d-1})}\leq C\,\|g\|^2_{L^2(\R^{d-1})},
\end{equation}
which follows from $w_\lambda\geq\frac{1}{2}$ and the coarea formula.

\medskip

\textit{Proof of \eqref{Est-K-k-L2}.} The multiplier of $K_k$ is bounded, in modulus, by $\sup_{\lambda\geq2\tau|\ell|}\|k(\lambda,\cdot)\|_{L^\infty}$, so that \eqref{Est-K-k-L2} follows from Plancherel's identity.

\medskip

\textit{Proof of \eqref{Est-K-k-Lp-Lp'}.} By Minkowski's integral inequality in $\lambda$ applied to \eqref{K-k-decomposition}, then by \eqref{R-Rstar-bounds} and $w_\lambda\leq2$, we have, for $f\in\mathscr S(\R^{d-1})$,
\begin{equation}\label{K-k-Minkowski}
\big\|K_kf\big\|_{L^{\frac{2d}{d-2}}}
\leq C\int_{2\tau|\ell|}^{\infty}\Big\|R_\lambda^\star\big[w_\lambda\,k_\lambda\,R_\lambda f\big]\Big\|_{L^{\frac{2d}{d-2}}}d\lambda
\leq C\int_{2\tau|\ell|}^{\infty}\lambda^{\frac{d-2}{2d}}\,\big\|k(\lambda,\cdot)\big\|_{L^\infty}\,
\big\|R_\lambda f\big\|_{L^2(\Sigma_\lambda)}\,d\lambda .
\end{equation}
Bounding $\|R_\lambda f\|_{L^2(\Sigma_\lambda)}\leq C\lambda^{\frac{d-2}{2d}}\|f\|_{L^{\frac{2d}{d+2}}}$ by \eqref{R-Rstar-bounds} then yields \eqref{Est-K-k-Lp-Lp'}, since $\frac{d-2}{2d}+\frac{d-2}{2d}=1-\frac{2}{d}$, and we conclude by density of $\mathscr S(\R^{d-1})$ in $L^{\frac{2d}{d+2}}(\R^{d-1})$.

\medskip

\textit{Proof of \eqref{Est-K-k-mixed}.} Applying instead the Cauchy--Schwarz inequality in $\lambda$ to the right-hand side of \eqref{K-k-Minkowski}, together with \eqref{Plancherel-layers}, we get, for $f\in\mathscr S(\R^{d-1})$,
\begin{multline*}
\big\|K_kf\big\|_{L^{\frac{2d}{d-2}}}
\leq C\left(\int_{2\tau|\ell|}^{\infty}\big\|k(\lambda,\cdot)\big\|^2_{L^\infty}\lambda^{1-\frac{2}{d}}\,d\lambda\right)^{\frac{1}{2}}
\left(\int_{2\tau|\ell|}^{\infty}\big\|R_\lambda f\big\|^2_{L^2(\Sigma_\lambda)}\,d\lambda\right)^{\frac{1}{2}}\\
\leq C\left(\int_{2\tau|\ell|}^{\infty}\big\|k(\lambda,\cdot)\big\|^2_{L^\infty}\lambda^{1-\frac{2}{d}}\,d\lambda\right)^{\frac{1}{2}}\|f\|_{L^2},
\end{multline*}
which, by density, is the second half of \eqref{Est-K-k-mixed}. The first half follows by duality: the adjoint of $K_k$ is $K_{\overline k}$, and $\|\overline k(\lambda,\cdot)\|_{L^\infty}=\|k(\lambda,\cdot)\|_{L^\infty}$, so that $\|K_k\|_{\mathcal L(L^{\frac{2d}{d+2}},L^2)}=\|K_{\overline k}\|_{\mathcal L(L^2,L^{\frac{2d}{d-2}})}$.
\end{proof} 
	%
    \subsection{$L^p$-estimates on the parametrix} \label{Sec-Est-Parametrix} 
    This subsection is devoted to the estimates of the operators of Proposition \ref{Prop-parametrix-Strip-1}. 

    \subsubsection{The low-frequency region} \label{Sec-low-frequencies} 
    On $\mathcal{R}_{\mathrm{low}}$, we have the following result: 
\begin{proposition}
\label{Prop-low-frequencies}
There exists $C>0$, depending only on $d$, $\alpha_1$ and $\alpha_2$, such that for all $\tau\geq1$ and all
$|\ell|\leq\ell_*$, the following estimates hold. For all $f\in L^{\frac{2d}{d+2}}(\mathcal S)$,
\begin{equation}
\label{Est-K-tau-0-f-L-2*'-low}
\big\|K_{\tau,0}\Pi_{\mathrm{low}}f\big\|_{L^{\frac{2d}{d-2}}(\mathcal S)}
+\tau \big\|K_{\tau,0}\Pi_{\mathrm{low}}f\big\|_{L^2(\mathcal S)}
+ \big\|\nabla K_{\tau,0}\Pi_{\mathrm{low}}f\big\|_{L^2(\mathcal S)}
\leq C\,\|f\|_{L^{\frac{2d}{d+2}}(\mathcal S)},
\end{equation}
for all $f\in L^2(\mathcal S)$,
\begin{equation}
\label{Est-K-tau-0-f-L-2-low}
\tau\big\|K_{\tau,0}\Pi_{\mathrm{low}}f\big\|_{L^{\frac{2d}{d-2}}(\mathcal S)}
+\tau^2\big\|K_{\tau,0}\Pi_{\mathrm{low}}f\big\|_{L^2(\mathcal S)}
+\tau\big\|\nabla K_{\tau,0}\Pi_{\mathrm{low}}f\big\|_{L^2(\mathcal S)}
\leq C\,\|f\|_{L^2(\mathcal S)},
\end{equation}
for all $j\in\{1,\dots,d\}$ and all $f\in L^2(\mathcal S)$,
\begin{equation}
\label{Est-K-tau-j-f-L-2-low}
\big\|K_{\tau,j}\Pi_{\mathrm{low}}f\big\|_{L^{\frac{2d}{d-2}}(\mathcal S)}
+\tau\big\|K_{\tau,j}\Pi_{\mathrm{low}}f\big\|_{L^2(\mathcal S)}
+\big\|\nabla K_{\tau,j}\Pi_{\mathrm{low}}f\big\|_{L^2(\mathcal S)}
\leq C\,\|f\|_{L^2(\mathcal S)},
\end{equation}
for all $g\in H^{\frac{1}{2}}(\R^{d-1})$,
\begin{equation}
\label{Est-G-tau-g-H-1/2-L-2_*-low}
\begin{gathered}
\tau^{\frac{1}{2}}\big\|G_\tau\Pi_{\mathrm{low}}g\big\|_{L^2(\mathcal S)}
+\tau^{-\frac{1}{2}}\big\|\nabla G_\tau\Pi_{\mathrm{low}}g\big\|_{L^2(\mathcal S)}
\leq C\,\|g\|_{L^2(\R^{d-1})},\\
\big\|G_\tau\Pi_{\mathrm{low}}g\big\|_{L^{\frac{2d}{d-2}}(\mathcal S)}
\leq C\,\|g\|_{H^{\frac{1}{2}}(\R^{d-1})},
\end{gathered}
\end{equation}
and, for all $\theta\in\{0,-\frac{1}{2}\}$ and all $h\in H^{\theta}(\R^{d-1})$,
\begin{equation}\label{Est-H-tau-g-H-1/2-L-2_*-low}
\tau^{\frac{3}{2}}\big\|H_\tau\Pi_{\mathrm{low}}h\big\|_{L^2(\mathcal S)}
+\tau^{\frac{1}{2}}\big\|\nabla H_\tau\Pi_{\mathrm{low}}h\big\|_{L^2(\mathcal S)}
+\tau^{\frac{1}{2}}\big\|H_\tau\Pi_{\mathrm{low}}h\big\|_{L^{\frac{2d}{d-2}}(\mathcal S)}
\leq C\,\tau^{-\theta}\,\|h\|_{H^{\theta}(\R^{d-1})} .
\end{equation}
\end{proposition}
\begin{remark}
These estimates imply those of Proposition \ref{Prop-Carleman-Strip-2} for the low-frequency part of \eqref{parametrix}. In $L^{\frac{2d}{d-2}}(\mathcal{S})$, the bounds of $G_\tau \Pi_{\mathrm{low}}$ and $H_\tau \Pi_{\mathrm{low}}$, and that of $K_{\tau,0}\Pi_{\mathrm{low}}$ on $L^{\frac{2d}{d+2}}(\mathcal{S})$, are scale invariant; all the other ones are better than required by at least a factor $\tau^{-\frac{1}{4}+\frac{1}{2d}}$. 
\end{remark}
\medskip 

Before proving this proposition, let us prove the following useful estimates, for all $x_1, y_1 \in [0,L]$ and $\xi' \in \R^{d-1}$ 
\begin{align}\label{sinh-bounds}
\Big|\frac{\sinh\big(\psi(\xi')(x_1-y_1)\big)}{\psi(\xi')}\Big|
&\leq 2 \min\Big\{|x_1-y_1|,\frac{1}{2\lambda}\Big\}e^{\lambda|x_1-y_1|},
\\ 
\big|\cosh\big(\psi(\xi')(x_1-y_1)\big)\big|&\leq e^{\lambda|x_1-y_1|},\label{cosh-bounds}
\end{align}
where $\lambda= \operatorname{Re} \psi(\xi')$. 

To prove that, set $\lambda=\operatorname{Re}\psi(\xi')$ and $\iota=\operatorname{Im}\psi(\xi')$, and writing $u:=\lambda|x_1-y_1|$, the identities
$|\sinh(a+ib)|^2=\sinh^2a+\sin^2b$ and $|\cosh(a+ib)|^2=\cosh^2a-\sin^2b$ give, together with
$|\iota|\leq|\psi(\xi')|$, $\sinh u\geq u$ and $\sinh(u)e^{-u}=\tfrac{1}{2}(1-e^{-2u})\leq\min\{u,\tfrac{1}{2}\}$,
\begin{align*}
\Big|\frac{\sinh\big(\psi(\xi')(x_1-y_1)\big)}{\psi(\xi')}\Big|
&\leq\frac{\sinh u}{|\psi(\xi')|}+\frac{|\iota||x_1-y_1|}{|\psi(\xi')|} \\& \leq2\,\frac{\sinh u}{\lambda} \nonumber
\leq2\min\Big\{|x_1-y_1|,\frac{1}{2\lambda}\Big\}e^{\lambda|x_1-y_1|},
\\ 
\big|\cosh\big(\psi(\xi')(x_1-y_1)\big)\big|&\leq e^{\lambda|x_1-y_1|} . \nonumber
\end{align*}
\begin{proof}
    \textit{Pointwise bounds.} Let $\xi' \in \mathcal{R}_{\mathrm{low}}$ be fixed and $\lambda=\operatorname{Re}\psi(\xi')$. By \eqref{low-below-turning}, $\lambda  <- \tau \varphi'(x_1)-\frac{\alpha_1}{2}\tau$ for all $x_\in [0,L]$. Consequently, in
\eqref{Def-k-tau-0-x1-y1-xi'} only the first line contributes, in \eqref{Def-k-tau-1-x1-y1-xi'} only the
first line and the last term, while in \eqref{Def-g-tau} both terms are present, the indicator in the first one being satisfied for every $x_1\in[0,L]$. In particular, all the kernels $k_{\tau,j}(x_1, y_1, \xi')$ vanish for $y_1>x_1$. Moreover, by \eqref{psi-bounds} and \eqref{low-below-turning},  $|\xi'|\leq\frac{3\alpha_1}{4}\tau$ and $|\psi(\xi')|\leq\alpha_1\tau$, so that $|\tau\varphi'(x)|+|\psi(\xi')|\leq2\alpha_2\tau$ for all $x\in[0,L]$. Using $\varphi'\leq-\alpha_1$, we have, for $y_1<x_1$,
$$
\tau\big(\varphi(x_1)-\varphi(y_1)\big)\leq-\alpha_1\tau\,(x_1-y_1),
$$
so that, by \eqref{sinh-bounds} and $\lambda<\frac{\alpha_1}{2}\tau$, for $y_1<x_1$,
\begin{multline}
    \label{low-pointwise}
\big|k_{\tau,0}(x_1,y_1,\xi')\big|\leq2\,(x_1-y_1)\,e^{-(\alpha_1\tau-\lambda)(x_1-y_1)}
\leq2\,(x_1-y_1)\,e^{-\frac{\alpha_1}{2}\tau(x_1-y_1)},
\\
\big|e^{\tau(\varphi(x_1)-\varphi(y_1))+(x_1-y_1)\psi(\xi')}\big|\leq e^{-\frac{\alpha_1}{2}\tau(x_1-y_1)} .
\end{multline}
Moreover, we have, for $x_1\neq y_1$,
\begin{equation}\label{low-d1}
\partial_{x_1}k_{\tau,0}(x_1,y_1,\xi')=\tau\varphi'(x_1)\,k_{\tau,0}(x_1,y_1,\xi')
+\mathbf 1_{y_1<x_1}\,e^{\tau(\varphi(x_1)-\varphi(y_1))}\cosh\big(\psi(\xi')(x_1-y_1)\big),
\end{equation}
where, by \eqref{cosh-bounds}, for $y_1<x_1$,
$$
\big|e^{\tau(\varphi(x_1)-\varphi(y_1))}\cosh\big(\psi(\xi')(x_1-y_1)\big)\big|
\leq e^{-\alpha_1\tau(x_1-y_1)}\cosh\big(\lambda(x_1-y_1)\big)
\leq e^{-\frac{\alpha_1}{2}\tau(x_1-y_1)} .
$$
Since $k_{\tau,0}(\cdot,y_1,x_1)$ is continuous and vanishes at $x_1= y_1$, the operator $\partial_1 K_{\tau,0} \Pi_\mathrm{low}$ has the kernel \eqref{low-d1}.  Let
\begin{equation}\label{Def-kappa-low}
\kappa(t):=(1+t)\,e^{-\frac{\alpha_1}{2}t}\,\mathbf 1_{t>0},\qquad t\in\R .
\end{equation}
Since $|k_{\tau,j}|=|\xi_j|\,|k_{\tau,0}|$ for $j\in\{2,\dots,d\}$, and since the last term of $k_{\tau,1}$ is bounded by $|\psi(\xi')-\tau\varphi'(y_1)|\,|k_{\tau,0}|\leq2\alpha_2\tau\,|k_{\tau,0}|$, the above bounds give, for $x_1\neq y_1$ and $j\in\{1,\dots,d\}$,
\begin{equation}\label{low-pointwise-kappa}
\begin{gathered}
\big|k_{\tau,0}(x_1,y_1,\xi')\big|\leq C\tau^{-1}\kappa\big(\tau(x_1-y_1)\big),\\
\big|\partial_{x_1}k_{\tau,0}(x_1,y_1,\xi')\big|+|\xi'|\,\big|k_{\tau,0}(x_1,y_1,\xi')\big|+\big|k_{\tau,j}(x_1,y_1,\xi')\big|\leq C\kappa\big(\tau(x_1-y_1)\big),\\
|\xi'|\,\big|k_{\tau,j}(x_1,y_1,\xi')\big|\leq C\tau\,\kappa\big(\tau(x_1-y_1)\big).
\end{gathered}
\end{equation}
\medskip

\textit{Reduction.} Let $m(x_1,y_1,\xi')$ be a kernel satisfying $|m(x_1,y_1,\xi')|\leq C\tau^b\kappa\big(\tau(x_1-y_1)\big)$ for all $\xi'\in\mathcal R_{\mathrm{low}}$ and $x_1\neq y_1$, and let $T_m$ be defined by $\widehat{T_mf}(x_1,\xi')=\mathbf 1_{\mathcal R_{\mathrm{low}}}(\xi')\int_0^Lm(x_1,y_1,\xi')\,\widehat f(y_1,\xi')\,dy_1$. Then
$$
\sup_{x_1\in[0,L]}\int_0^L\big|m(x_1,y_1,\xi')\big|\,dy_1
+\sup_{y_1\in[0,L]}\int_0^L\big|m(x_1,y_1,\xi')\big|\,dx_1
\leq2C\tau^b\int_\R\kappa(\tau z)\,dz=2C\,\|\kappa\|_{L^1(\R)}\,\tau^{b-1},
$$
so that Schur's test in the variable $x_1$, applied for each fixed $\xi'\in\mathcal R_{\mathrm{low}}$, together with Plancherel's theorem in $x'$, yields
\begin{equation}\label{low-L2-reduction}
\big\|T_mf\big\|_{L^2(\mathcal S)}\leq C\tau^{b-1}\,\|f\|_{L^2(\mathcal S)} .
\end{equation}
\medskip

\textit{$L^2$-estimates for $K_{\tau,0}$.} By \eqref{low-pointwise-kappa}, \eqref{low-L2-reduction} applies to $K_{\tau,0}\Pi_{\mathrm{low}}$ with $b=-1$, and, with $b=0$, to $\partial_1K_{\tau,0}\Pi_{\mathrm{low}}$ and to the tangential derivatives, which amount to multiplication by $i\xi'$. Hence
\begin{align}
\tau^{2}\big\|K_{\tau,0}\Pi_{\mathrm{low}}f\big\|_{L^2(\mathcal S)}&\leq C\,\|f\|_{L^2(\mathcal S)},\label{low-K0-L2}\\
\tau\big\|\nabla K_{\tau,0}\Pi_{\mathrm{low}}f\big\|_{L^2(\mathcal S)}&\leq C\,\|f\|_{L^2(\mathcal S)} .\label{low-K0-grad}
\end{align}
\medskip

\textit{$L^2$-estimates for $K_{\tau,j}$.} By \eqref{low-pointwise-kappa}, \eqref{low-L2-reduction} applies to $K_{\tau,j}\Pi_{\mathrm{low}}$ with $b=0$, and to its tangential derivatives with $b=1$. For $j\in\{2,\dots,d\}$, $\partial_1K_{\tau,j}\Pi_{\mathrm{low}}$ has the kernel $i\xi_j\,\partial_{x_1}k_{\tau,0}(x_1,y_1,\xi')$, bounded by $C\tau\kappa\big(\tau(x_1-y_1)\big)$. For $j=1$, since $k_{\tau,1}(x_1,y_1,\xi')$ vanishes for $y_1>x_1$ and tends to $1$ as $y_1\to x_1^-$, differentiating $\widehat{K_{\tau,1}\Pi_{\mathrm{low}}f}(x_1,\xi')=\mathbf 1_{\mathcal R_{\mathrm{low}}}(\xi')\int_0^{x_1}k_{\tau,1}(x_1,y_1,\xi')\,\widehat f(y_1,\xi')\,dy_1$ with respect to $x_1$ gives
\begin{equation}\label{low-d1-K1}
\partial_1K_{\tau,1}\Pi_{\mathrm{low}}f=\Pi_{\mathrm{low}}f+K_{\tau,1,\partial_1}\Pi_{\mathrm{low}}f ,
\end{equation}
where the kernel of $K_{\tau,1,\partial_1}$, namely
$$
\mathbf 1_{y_1<x_1}\big(\tau\varphi'(x_1)+\psi(\xi')\big)e^{\tau(\varphi(x_1)-\varphi(y_1))+(x_1-y_1)\psi(\xi')}-\big(\psi(\xi')-\tau\varphi'(y_1)\big)\,\partial_{x_1}k_{\tau,0}(x_1,y_1,\xi'),
$$
is bounded by $C\tau\kappa\big(\tau(x_1-y_1)\big)$, by \eqref{low-pointwise} and \eqref{low-pointwise-kappa}. Since $\|\Pi_{\mathrm{low}}f\|_{L^2(\mathcal S)}\leq\|f\|_{L^2(\mathcal S)}$, \eqref{low-L2-reduction} gives
\begin{equation}\label{low-Kj-L2}
\tau\big\|K_{\tau,j}\Pi_{\mathrm{low}}f\big\|_{L^2(\mathcal S)}
+\big\|\nabla K_{\tau,j}\Pi_{\mathrm{low}}f\big\|_{L^2(\mathcal S)}\leq C\,\|f\|_{L^2(\mathcal S)} .
\end{equation}

\textit{$L^2$-estimates for $G_\tau$ and $H_\tau$.} Both are Fourier multipliers in $x'$, with kernels $g_\tau(x_1,\xi')$ and $k_{\tau,0}(x_1,0,\xi')$ given by \eqref{Def-g-tau} and \eqref{Def-k-tau-0-x1-y1-xi'} respectively. Since $\varphi(0)=0$, \eqref{low-pointwise} with $y_1=0$ gives $|e^{\tau\varphi(x_1)+x_1\psi(\xi')}|\leq e^{-\frac{\alpha_1}{2}\tau x_1}$. Hence, by \eqref{Def-g-tau}, \eqref{low-pointwise-kappa} with $y_1=0$, $|\tau\varphi'(x)|+|\psi(\xi')|\leq2\alpha_2\tau$ and
$$
\partial_{x_1}g_\tau(x_1,\xi')=\big(\tau\varphi'(x_1)+\psi(\xi')\big)e^{\tau\varphi(x_1)+x_1\psi(\xi')}-\big(\tau\varphi'(0)+\psi(\xi')\big)\,\partial_{x_1}k_{\tau,0}(x_1,0,\xi'),
$$
we have, for all $x_1\in(0,L]$,
\begin{equation}\label{low-pointwise-boundary}
\begin{gathered}
\big|k_{\tau,0}(x_1,0,\xi')\big|\leq C\tau^{-1}\kappa(\tau x_1),\\\ 
\big|\partial_{x_1}k_{\tau,0}(x_1,0,\xi')\big|+|\xi'|\,\big|k_{\tau,0}(x_1,0,\xi')\big|+\big|g_\tau(x_1,\xi')\big|\leq C\kappa(\tau x_1),\\
\big|\partial_{x_1}g_\tau(x_1,\xi')\big|+|\xi'|\,\big|g_\tau(x_1,\xi')\big|\leq C\tau\,\kappa(\tau x_1).
\end{gathered}
\end{equation}
Let $s(x_1,\xi')$ be a function satisfying $|s(x_1,\xi')|\leq C\tau^b\kappa(\tau x_1)$ for all $x_1\in(0,L]$ and $\xi'\in\mathcal R_{\mathrm{low}}$, and let $S_s$ be defined by $\widehat{S_su}(x_1,\xi')=\mathbf 1_{\mathcal R_{\mathrm{low}}}(\xi')\,s(x_1,\xi')\,\widehat u(\xi')$. Plancherel's theorem in $x'$ gives
\begin{equation}\label{low-L2-reduction-GH}
\big\|S_su\big\|_{L^2(\mathcal S)}\leq C\tau^b\,\big\|\kappa(\tau\,\cdot)\big\|_{L^2(\R)}\,\big\|\Pi_{\mathrm{low}}u\big\|_{L^2(\R^{d-1})}=C\tau^{b-\frac{1}{2}}\,\big\|\Pi_{\mathrm{low}}u\big\|_{L^2(\R^{d-1})} .
\end{equation}
By \eqref{low-pointwise-boundary}, \eqref{low-L2-reduction-GH} applies to $H_\tau\Pi_{\mathrm{low}}$ with $b=-1$, to $G_\tau\Pi_{\mathrm{low}}$, $\partial_1H_\tau\Pi_{\mathrm{low}}$ and the tangential derivatives of $H_\tau\Pi_{\mathrm{low}}$ with $b=0$, and to $\partial_1G_\tau\Pi_{\mathrm{low}}$ and the tangential derivatives of $G_\tau\Pi_{\mathrm{low}}$ with $b=1$. Hence
\begin{align}
\big\|H_\tau\Pi_{\mathrm{low}}h\big\|_{L^2(\mathcal S)}&\leq C\tau^{-\frac{3}{2}}\|\Pi_{\mathrm{low}}h\|_{L^2},
&
\big\|G_\tau\Pi_{\mathrm{low}}g\big\|_{L^2(\mathcal S)}&\leq C\tau^{-\frac{1}{2}}\|g\|_{L^2},\label{low-GH-L2}\\
\big\|\nabla H_\tau\Pi_{\mathrm{low}}h\big\|_{L^2(\mathcal S)}&\leq C\tau^{-\frac{1}{2}}\|\Pi_{\mathrm{low}}h\|_{L^2},
&
\big\|\nabla G_\tau\Pi_{\mathrm{low}}g\big\|_{L^2(\mathcal S)}&\leq C\tau^{\frac{1}{2}}\|g\|_{L^2} .\label{low-GH-grad}
\end{align}

\textit{ $L^\frac{2d}{d-2}$-estimates.} By \eqref{low-below-turning}, all the functions above
have, for each $x_1$, tangential Fourier support contained in $B\big(0,\frac{3\alpha_1}{4}\tau\big)$. For $u\in\mathscr S(\R^{d-1})$ with $\operatorname{supp}\widehat u\subset B(0,\rho)$, the Hausdorff--Young inequality $\|u\|_{L^{\frac{2d}{d-2}}}\leq\|\widehat u\|_{L^{\frac{2d}{d+2}}}$ and Hölder's inequality with $\frac{d+2}{2d}=\frac1d+\frac{1}{2}$ give Bernstein's type inequality in $\R^{d-1}$ and its $H^{\frac{1}{2}}$ version,
\begin{equation}\label{low-Bernstein}
\|u\|_{L^{\frac{2d}{d-2}}(\R^{d-1})}\leq C\rho^{\frac{d-1}{d}}\|u\|_{L^2(\R^{d-1})},
\qquad
\|u\|_{L^{\frac{2d}{d-2}}(\R^{d-1})}\leq C\rho^{\frac{d-2}{2d}}\|u\|_{H^{\frac{1}{2}}(\R^{d-1})},
\end{equation}
the second one since $\int_{B(0,\rho)}(1+|\xi'|^2)^{-\frac d4}\,d\xi'\leq\int_{B(0,\rho)}|\xi'|^{-\frac{d}{2}}\,d\xi'=C\rho^{\frac{d-2}{2}}$ for $d\geq3$, as well as the dual form
$\|\Pi_{\mathrm{low}}f\|_{L^2(\R^{d-1})}\leq C\tau^{\frac{d-1}{d}}\|f\|_{L^{\frac{2d}{d+2}}(\R^{d-1})}$. Now, let $m$ and $T_m$ be as above. Plancherel's theorem in $x'$ and Minkowski's inequality in $y_1$ give, for all $x_1\in(0,L)$,
\begin{equation}\label{low-reduction}
\begin{aligned}
\big\|T_mf(x_1,\cdot)\big\|_{L^2_{x'}(\R^{d-1})}
&\leq\int_0^L\sup_{\xi'\in\mathcal R_{\mathrm{low}}}\big|m(x_1,y_1,\xi')\big|\,\big\|\Pi_{\mathrm{low}}f(y_1,\cdot)\big\|_{L^2_{y'}(\R^{d-1})}\,dy_1\\
&\leq C\tau^b\int_0^L\kappa\big(\tau(x_1-y_1)\big)\,\big\|\Pi_{\mathrm{low}}f(y_1,\cdot)\big\|_{L^2_{y'}(\R^{d-1})}\,dy_1 .
\end{aligned}
\end{equation}
Let $q_1\in\{\frac{2d}{d+2},2\}$, $q_2\in\{2,\frac{2d}{d-2}\}$ and $\frac1r=1+\frac1{q_2}-\frac1{q_1}$. Using \eqref{low-Bernstein} in $x'$ if $q_2=\frac{2d}{d-2}$, then \eqref{low-reduction} and Young's inequality in $x_1$, with $\|\kappa(\tau\,\cdot)\|_{L^r(\R)}=\tau^{-\frac1r}\|\kappa\|_{L^r(\R)}$, and finally the dual form of \eqref{low-Bernstein} if $q_1=\frac{2d}{d+2}$, we obtain
\begin{equation}\label{low-transfer}
\begin{aligned}
\|T_mf\|_{L^{q_2}(\mathcal S)}
&\leq C\tau^{(d-1)(\frac{1}{2}-\frac1{q_2})}\Big\|\,\big\|T_mf(x_1,\cdot)\big\|_{L^2_{x'}(\R^{d-1})}\Big\|_{L^{q_2}_{x_1}(0,L)}\\
&\leq C\tau^{(d-1)(\frac{1}{2}-\frac1{q_2})+b-\frac1r}\Big\|\,\big\|\Pi_{\mathrm{low}}f(y_1,\cdot)\big\|_{L^2_{y'}(\R^{d-1})}\Big\|_{L^{q_1}_{y_1}(0,L)}\\
&\leq C\tau^{\,b-1+d(\frac1{q_1}-\frac1{q_2})}\,\|f\|_{L^{q_1}(\mathcal S)} ,
\end{aligned}
\end{equation}
since $(d-1)(\frac1{q_1}-\frac1{q_2})-\frac1r=d(\frac1{q_1}-\frac1{q_2})-1$. Here $d(\frac1{q_1}-\frac1{q_2})$ equals $2$ for $(q_1,q_2)=(\frac{2d}{d+2},\frac{2d}{d-2})$, $1$ if exactly one of $q_1$ and $q_2$ equals $2$, and $0$ for $(q_1,q_2)=(2,2)$, in which case \eqref{low-transfer} is \eqref{low-L2-reduction}. By \eqref{low-pointwise-kappa}, \eqref{low-transfer} applies to $K_{\tau,0}\Pi_{\mathrm{low}}$ with $b=-1$, and, with $b=0$, to $\partial_1K_{\tau,0}\Pi_{\mathrm{low}}$, to the tangential derivatives of $K_{\tau,0}\Pi_{\mathrm{low}}$ and to $K_{\tau,j}\Pi_{\mathrm{low}}$. We get
\begin{gather*}
\big\|K_{\tau,0}\Pi_{\mathrm{low}}f\big\|_{L^{\frac{2d}{d-2}}(\mathcal S)}+\tau\big\|K_{\tau,0}\Pi_{\mathrm{low}}f\big\|_{L^2(\mathcal S)}+\big\|\nabla K_{\tau,0}\Pi_{\mathrm{low}}f\big\|_{L^2(\mathcal S)}\leq C\|f\|_{L^{\frac{2d}{d+2}}(\mathcal S)},\\
\tau\big\|K_{\tau,0}\Pi_{\mathrm{low}}f\big\|_{L^{\frac{2d}{d-2}}(\mathcal S)}+\big\|K_{\tau,j}\Pi_{\mathrm{low}}f\big\|_{L^{\frac{2d}{d-2}}(\mathcal S)}\leq C\|f\|_{L^2(\mathcal S)},
\end{gather*}
which, together with \eqref{low-K0-L2}, \eqref{low-K0-grad} and \eqref{low-Kj-L2}, gives \eqref{Est-K-tau-0-f-L-2*'-low}, \eqref{Est-K-tau-0-f-L-2-low} and \eqref{Est-K-tau-j-f-L-2-low}.

For $G_\tau$ and $H_\tau$, there is no integral in $y_1$: by \eqref{low-pointwise-boundary} and \eqref{low-Bernstein}, for all $x_1\in(0,L]$,
\begin{equation*}
\begin{aligned}
\big\|G_\tau\Pi_{\mathrm{low}}g(x_1,\cdot)\big\|_{L^{\frac{2d}{d-2}}_{x'}(\R^{d-1})}&\leq C\tau^{\frac{d-2}{2d}}\kappa(\tau x_1)\,\|g\|_{H^{\frac{1}{2}}(\R^{d-1})},\\
\big\|H_\tau\Pi_{\mathrm{low}}h(x_1,\cdot)\big\|_{L^{\frac{2d}{d-2}}_{x'}(\R^{d-1})}&\leq C\tau^{\frac{d-1}{d}-1}\kappa(\tau x_1)\,\|\Pi_{\mathrm{low}}h\|_{L^2(\R^{d-1})},
\end{aligned}
\end{equation*}
so that, taking the $L^{\frac{2d}{d-2}}$ norm in $x_1$, with $\|\kappa(\tau\,\cdot)\|_{L^{\frac{2d}{d-2}}(\R)}=C\tau^{-\frac{d-2}{2d}}$,
\begin{equation}\label{low-GH-2*}
\big\|G_\tau\Pi_{\mathrm{low}}g\big\|_{L^{\frac{2d}{d-2}}(\mathcal S)}\leq C\|g\|_{H^{\frac{1}{2}}(\R^{d-1})},
\qquad
\big\|H_\tau\Pi_{\mathrm{low}}h\big\|_{L^{\frac{2d}{d-2}}(\mathcal S)}\leq C\tau^{-\frac{1}{2}}\|\Pi_{\mathrm{low}}h\|_{L^2} .
\end{equation}

\textit{Passage to the Sobolev norms.} For $G_\tau$, \eqref{low-GH-L2}--\eqref{low-GH-grad} and \eqref{low-GH-2*} are exactly \eqref{Est-G-tau-g-H-1/2-L-2_*-low}. For $H_\tau$, they give \eqref{Est-H-tau-g-H-1/2-L-2_*-low} for $\theta=0$. For $\theta=-\frac{1}{2}$, we use that $(1+|\xi'|^2)^{\frac{1}{2}}\leq C\tau$ on $\mathcal R_{\mathrm{low}}$ by \eqref{low-below-turning}, whence
\begin{equation}\label{low-L2-vs-H-1/2}
\big\|\Pi_{\mathrm{low}}h\big\|^2_{L^2(\R^{d-1})}
=\int_{\mathcal R_{\mathrm{low}}}\big(1+|\xi'|^2\big)^{\frac{1}{2}}\big(1+|\xi'|^2\big)^{-\frac{1}{2}}
\big|\widehat h(\xi')\big|^2\,d\xi'
\leq C\,\tau\,\|h\|^2_{H^{-\frac{1}{2}}(\R^{d-1})} ,
\end{equation}
so that $\|\Pi_{\mathrm{low}}h\|_{L^2}\leq C\tau^{\frac{1}{2}}\|h\|_{H^{-\frac{1}{2}}}$, which, combined with the case $\theta=0$, gives \eqref{Est-H-tau-g-H-1/2-L-2_*-low} for $\theta=-\frac{1}{2}$.
\end{proof} 

	%
	%
	
	 \subsubsection{The high-frequency region}\label{Subsubsec-The high-frequency region} 
We start this subsection by proving the boundedness of the interior operators $K_{\tau,j}$. 
\begin{proposition}\label{Prop-K-high}
There exists $C>0$, depending only on $d,\alpha_1,\alpha_2,\beta$, such that for all $\tau\geq1$ and all
$|\ell|\leq\ell_*$, we have, for all $f\in L^2(\mathcal S)$,
\begin{equation}\label{K-high-1}
\tau^{\frac{3}{4}+\frac1{2d}}\big\|K_{\tau,0}\Pi_{\mathrm{high}}f\big\|_{L^{\frac{2d}{d-2}}(\mathcal S)}
+\tau^{\frac{3}{2}}\big\|K_{\tau,0}\Pi_{\mathrm{high}}f\big\|_{L^2(\mathcal S)}
+\tau^{\frac{1}{2}}\big\|\nabla K_{\tau,0}\Pi_{\mathrm{high}}f\big\|_{L^2(\mathcal S)}
\leq C\,\|f\|_{L^2(\mathcal S)},
\end{equation}
for all $f\in L^{\frac{2d}{d+2}}(\mathcal S)$,
\begin{equation}\label{K-high-2}
\big\|K_{\tau,0}\Pi_{\mathrm{high}}f\big\|_{L^{\frac{2d}{d-2}}(\mathcal S)}
+\tau^{\frac{3}{4}+\frac1{2d}}\big\|K_{\tau,0}\Pi_{\mathrm{high}}f\big\|_{L^2(\mathcal S)} + \tau^{-\frac{1}{4}+\frac1{2d}} \big\|\nabla K_{\tau,0}\Pi_{\mathrm{high}}f\big\|_{L^2(\mathcal S)}
\leq C\,\|f\|_{L^{\frac{2d}{d+2}}(\mathcal S)},
\end{equation}
 and for all $j\in\{1,\dots,d\}$ and all $f\in L^2(\mathcal S)$,
\begin{equation}\label{K-high-3}
\tau^{-\frac{1}{4}+\frac1{2d}}\big\|K_{\tau,j}\Pi_{\mathrm{high}}f\big\|_{L^{\frac{2d}{d-2}}(\mathcal S)}
+\tau^{\frac{1}{2}}\big\|K_{\tau,j}\Pi_{\mathrm{high}}f\big\|_{L^2(\mathcal S)}
+\tau^{-\frac{1}{2}}\big\|\nabla K_{\tau,j}\Pi_{\mathrm{high}}f\big\|_{L^2(\mathcal S)}
\leq C\,\|f\|_{L^2(\mathcal S)}.
\end{equation}
\end{proposition}

Before estimating the operators, we prove the following pointwise bounds on the kernels of Proposition \ref{Prop-parametrix-Strip-1}. Following \eqref{Def-polar-coordinates}, we identify with the pair $(\lambda, \omega') \in \R^*_+\times \Sigma_{s,e_\ell}$, $s=\tau |\ell|/\lambda$, through $\lambda= \operatorname{Re}\psi(\xi')$ and $\omega'= \xi'/\lambda$, and we write, with a slight abuse of notation, $k_{\tau,j}(x_1,y_1,\lambda,\omega')=k_{\tau,j}(x_1,y_1,\lambda\omega')$. We set, for $(x_1,y_1)\in[0,L]^2$ and $\lambda>0$,
\begin{equation}\label{Def-Theta-E}
\begin{gathered}
\Theta_\tau(x_1,y_1,\lambda):=
\begin{cases}
\min\Big\{|x_1-y_1|,\dfrac1{2\lambda}\Big\}, & \lambda<-\tau\varphi'(x_1),\\[3mm]
\dfrac1{2\lambda}, & \lambda\geq-\tau\varphi'(x_1),
\end{cases}\\[2mm]
\mathcal E_\tau(x_1,y_1,\lambda):=e^{-\frac{1}{2}|\lambda+\tau\varphi'(x_1)||x_1-y_1|-\frac{\beta}{2}\tau(x_1-y_1)^2} .
\end{gathered}
\end{equation}

\begin{lemma}\label{Lemma-kernels-bounds}
There exists $C>0$, depending only on $\alpha_1,\alpha_2,$ and $\beta$, such that for all $\tau\geq1$, all
$\ell\in\R^{d-1}$ with $|\ell|\leq\ell_*$, all $(x_1,y_1)\in[0,L]^2$ and all $\lambda>0$ with $\lambda \geq 2 \tau |\ell|$, writing
$\Sigma:=\Sigma_{\tau|\ell|/\lambda,\,e_\ell}$, so that $\big\|k (x_1,y_1,\lambda,\cdot)\big\|_{L^\infty(\Sigma)}= \sup_{\xi' \in \Sigma_\lambda}|k(x_1,y_1,\xi')|$, which we also denote by $\big\|k (x_1,y_1,\lambda,\cdot)\big\|_{L^\infty(\Sigma)}$
\begin{align}
\big\|k_{\tau,0}(x_1,y_1,\lambda,\cdot)\big\|_{L^\infty(\Sigma)}
&\leq C\,\Theta_\tau(x_1,y_1,\lambda)\,\mathcal E_\tau(x_1,y_1,\lambda),
\label{pw-k0}\\
\big\|k_{\tau,j}(x_1,y_1,\lambda,\cdot)\big\|_{L^\infty(\Sigma)}
&\leq C\,\lambda\,\Theta_\tau(x_1,y_1,\lambda)\,\mathcal E_\tau(x_1,y_1,\lambda),
\qquad j\in\{2,\dots,d\},
\label{pw-kj}\\
\big\|k_{\tau,1}(x_1,y_1,\lambda,\cdot)\big\|_{L^\infty(\Sigma)}
&\leq C\Big(1+(\lambda+\tau)\,\Theta_\tau(x_1,y_1,\lambda)\Big)\,\mathcal E_\tau(x_1,y_1,\lambda),
\label{pw-k1}\\
\big\|\partial_1k_{\tau,0}(x_1,y_1,\lambda,\cdot)\big\|_{L^\infty(\Sigma)}
&\leq C\,\mathcal E_\tau(x_1,y_1,\lambda),
\label{pw-d1k0}\\
\big\|\partial_1k_{\tau,1}(x_1,y_1,\lambda,\cdot)\big\|_{L^\infty(\Sigma)}
&\leq C\,(\lambda+\tau)\,\mathcal E_\tau(x_1,y_1,\lambda),
\label{pw-d1k1}
\end{align}
the last estimate being valid for $\lambda\neq-\tau\varphi'(x_1)$ and $x_1\neq y_1$. Moreover, the multipliers of the boundary
operators satisfy
\begin{multline}\label{pw-boundary}
\big\|k_{\tau,0}(x_1,0,\lambda,\cdot)\big\|_{L^\infty(\Sigma)}
\leq C\,\Theta_\tau(x_1,0,\lambda)\,\mathcal E_\tau(x_1,0,\lambda),
\\
\big\|g_\tau(x_1,\lambda,\cdot)\big\|_{L^\infty(\Sigma)}
\leq C\Big(1+(\lambda+\tau)\Theta_\tau(x_1,0,\lambda)\Big)\mathcal E_\tau(x_1,0,\lambda).
\end{multline}
\end{lemma}

\begin{proof} The right-hand sides depend on $\xi'$ on through $\lambda$, so that it suffices to prove the corresponding bounds pointwise for $\xi' \in \Sigma_\lambda$. 

We begin with some elementary bounds. First from \eqref{psi-bounds}, $\lambda \geq 2 \tau |\ell |$ and $|\ell| \leq \frac{\alpha_1}{4}$, one can get  
\begin{equation}\label{psi-modulus}
\begin{gathered}
\lambda\ \leq\ |\psi(\xi')|\ \leq\ |\xi'|+\tau|\ell|\ \leq\ \lambda+2\tau|\ell|\ \leq\ \lambda+\tfrac{\alpha_1}2\tau, \\ 
 |\operatorname{Im} \psi(\xi')| \leq \tau |\ell| \leq \frac{\alpha_1}{4} \tau, \quad\text{and}\quad
|\xi'|\leq\lambda+\tau|\ell|\leq\tfrac32\lambda\  
\end{gathered}
\end{equation}
On the other hand, the convexity of $\varphi$ implies 
\begin{equation}\label{convexity-ineq}
\varphi(x_1)-\varphi(y_1)\ \leq\ \varphi'(x_1)(x_1-y_1) -\frac{\beta}{2}\,(x_1-y_1)^2,
\end{equation}
Using this inequality,  for $\lambda\leq-\tau\varphi'(x_1)$ and all $y_1<x_1$, we get
\begin{equation}\label{exponent-forward}
\tau\big(\varphi(x_1)-\varphi(y_1)\big)+\lambda(x_1-y_1)
\ \leq\ \big(\lambda+\tau\varphi'(x_1)\big)(x_1-y_1)-\frac{\beta}{2}\tau(x_1-y_1)^2 ,
\end{equation}
the first term of the right-hand side equals $-\big|\lambda+\tau\varphi'(x_1)\big|\,|x_1-y_1|$. Similarly, for
$\lambda\geq-\tau\varphi'(x_1)$ and all $(x_1,y_1)\in[0,L]^2$,
\begin{equation}\label{exponent-backward}
\tau\big(\varphi(x_1)-\varphi(y_1)\big)-\lambda|x_1-y_1|
\ \leq\ -\big|\lambda+\tau\varphi'(x_1)\big|\,|x_1-y_1|-\frac{\beta}{2}\tau(x_1-y_1)^2, 
\end{equation}
indeed, if $y_1<x_1$ the left-hand side is bounded by
$\big(\tau\varphi'(x_1)-\lambda\big)(x_1-y_1)-\tfrac{\beta}{2}\tau(x_1-y_1)^2$ and
$\lambda-\tau\varphi'(x_1)\geq\lambda+\tau\varphi'(x_1)=|\lambda+\tau\varphi'(x_1)|$ since $\varphi'\leq0$;
and if $y_1>x_1$ it is bounded by
$-\big(\lambda+\tau\varphi'(x_1)\big)(y_1-x_1)-\tfrac{\beta}{2}\tau(x_1-y_1)^2$, and
$\lambda+\tau\varphi'(x_1)=|\lambda+\tau\varphi'(x_1)|$ on this branch.
\medskip 

Now, we prove our estimates: 

\textit{Proof of \eqref{pw-k0} and \eqref{pw-boundary}.} If
$\lambda<-\tau\varphi'(x_1)$ and $y_1<x_1$, then by \eqref{Def-k-tau-0-x1-y1-xi'}, \eqref{sinh-bounds} and \eqref{exponent-forward},
\begin{multline*}
  \big|k_{\tau,0}(x_1,y_1,\xi')\big|
=e^{\tau(\varphi(x_1)-\varphi(y_1))}\Big|\frac{\sinh(\psi(\xi') (x_1-y_1))}{\psi(\xi')}\Big| \\ 
\leq2\min\Big\{|x_1-y_1|,\frac1{2\lambda}\Big\}\,
e^{\tau(\varphi(x_1)-\varphi(y_1))+\lambda(x_1-y_1)} \leq 2 \Theta_\tau\,\mathcal E_\tau  . 
\end{multline*}

If $\lambda\geq-\tau\varphi'(x_1)$, then, bounding
$\big|\int_0^{L-\max\{x_1,y_1\}}e^{-2\psi\tilde x_1}d\tilde x_1\big|
\leq\int_0^{\infty}e^{-2\lambda\tilde x_1}d\tilde x_1=\tfrac{1}{2\lambda}$ and using
\eqref{exponent-backward},
$$
\big|k_{\tau,0}(x_1,y_1,\xi')\big|\leq\frac{1}{2\lambda}\,
e^{\tau(\varphi(x_1)-\varphi(y_1))-\lambda|x_1-y_1|}\leq\Theta_\tau\,\mathcal E_\tau ,
$$
which proves \eqref{pw-k0}, and \eqref{pw-boundary}$_1$ by taking $y_1=0$. 
\medskip

\textit{Proof of \eqref{pw-kj}.} Next, for all $j\in \{2,\cdots,d\}$, we have  $|k_{\tau,j} (x_1,y_1,\xi')|=|\xi_j|\,|k_{\tau,0}(x_1,y_1,\xi')|$ with $|\xi'|\leq\lambda+\tau|\ell|\leq2\lambda$ on $\mathcal R_{\mathrm{high}}$, whence \eqref{pw-kj}.
\medskip

\textit{Proof of  \eqref{pw-k1} and \eqref{pw-boundary}.}  For $k_{\tau,1}$, the first two terms of
\eqref{Def-k-tau-1-x1-y1-xi'} have modulus
$e^{\tau(\varphi(x_1)-\varphi(y_1))+\operatorname{Re}(\psi)(x_1-y_1)}\leq\mathcal E_\tau$, by
\eqref{exponent-forward} on the first branch and \eqref{exponent-backward} on the second, while, by
\eqref{psi-modulus} and $|\varphi'|\leq\alpha_2$,
$|\psi-\tau\varphi'(y_1)|\leq\lambda+\tfrac32\alpha_2\tau\leq C(\lambda+\tau)$, combined with
\eqref{pw-k0}, this gives \eqref{pw-k1}, and \eqref{pw-boundary}$_2$ follows in the same way from
\eqref{Def-g-tau}, since $|\tau\varphi'(0)+\psi(\xi')|\leq C(\lambda+\tau)$ and $\varphi(0)=0$.

\medskip

\textit{Proof of  \eqref{pw-d1k0}.} Assume first $\lambda<-\tau\varphi'(x_1)$ and $y_1<x_1$, since $k_{\tau, 0}=0$ for $y_1>x_1$ on this branch. Differentiating $k_{\tau,0}(x_1,y_1,\xi')=e^{\tau(\varphi(x_1)-\varphi(y_1))}\sinh\big(\psi(\xi') (x_1-y_1)\big)/\psi$ with respect to $x_1$ and using $\cosh z=\sinh z+e^{-z}$,
\begin{equation}\label{d1-k0-forward}
\partial_1k_{\tau,0}(x_1,y_1,\xi')
=\big(\psi(\xi')+\tau\varphi'(x_1)\big)\,k_{\tau,0}(x_1,y_1,\xi')
+e^{\tau(\varphi(x_1)-\varphi(y_1))-\psi(\xi')(x_1-y_1)} .
\end{equation}
The second term has modulus $\leq e^{\tau(\varphi(x_1)-\varphi(y_1))-\lambda(x_1-y_1)}\leq\mathcal E_\tau$,
by \eqref{exponent-backward}, whose proof for $y_1<x_1$ does not use the branch condition. As for the
first, \eqref{psi-modulus} gives
$|\psi+\tau\varphi'(x_1)|\leq|\lambda+\tau\varphi'(x_1)|+\tfrac{\alpha_1}2\tau$, and we distinguish two
cases. 

If $\lambda\geq\tfrac{1}{2}\big(-\tau\varphi'(x_1)\big)$, then
$\Theta_\tau\leq\tfrac{1}{2\lambda}\leq\tfrac{1}{\alpha_1\tau}$ and
$|\psi(\xi')+\tau\varphi'(x_1)|\leq(\alpha_2+\tfrac{\alpha_1}2)\tau$, so that
$|\psi(\xi')+\tau\varphi'(x_1)|\,\Theta_\tau\leq C$. 

If $\lambda<\tfrac{1}{2}\big(-\tau\varphi'(x_1)\big)$, then
$|\lambda+\tau\varphi'(x_1)|\geq\tfrac{1}{2}\big(-\tau\varphi'(x_1)\big)\geq\tfrac{\alpha_1}2\tau$, so that
$|\psi+\tau\varphi'(x_1)|\leq C|\lambda+\tau\varphi'(x_1)|$ and, by $\Theta_\tau\leq|x_1-y_1|$ and
$ue^{-\frac{u}{2}}\leq2e^{-1}e^{-\frac{u}{4}}$,
$$
|\psi+\tau\varphi'(x_1)|\,\Theta_\tau\,\mathcal E_\tau
\leq C\,\big|\lambda+\tau\varphi'(x_1)\big|\,|x_1-y_1|\,
e^{-\frac{1}{2}|\lambda+\tau\varphi'(x_1)||x_1-y_1|}e^{-\frac{\beta}{2}\tau(x_1-y_1)^2}\leq C\,\mathcal E_\tau .
$$
In both cases $\big|(\psi+\tau\varphi'(x_1))k_{\tau,0}\big|\leq C\mathcal E_\tau$, which gives
\eqref{pw-d1k0} on this branch.

If now $\lambda>-\tau\varphi'(x_1)$, differentiating the second line of
\eqref{Def-k-tau-0-x1-y1-xi'} gives
\begin{multline*}
    \partial_1k_{\tau,0}(x_1,y_1,\xi')=\Big(\tau\varphi'(x_1)-\psi(\xi')\,\mathrm{sgn}(x_1-y_1)\Big)k_{\tau,0}(x_1,y_1,\xi') \\ 
+\mathbf 1_{y_1<x_1}\,e^{\tau(\varphi(x_1)-\varphi(y_1))-\psi(\xi')(x_1-y_1)-2\psi(\xi')(L-x_1)} ,
\end{multline*}
the modulus of the last term is again
$\leq\mathcal E_\tau$ by \eqref{exponent-backward}, while, by \eqref{psi-modulus} and
$-\tau\varphi'(x_1)\leq\lambda$,
$\big|\tau\varphi'(x_1)-\psi\,\mathrm{sgn}(x_1-y_1)\big|\leq2\lambda+\tfrac{\alpha_1}2\tau\leq C\lambda$,
so that, since $\Theta_\tau=\tfrac1{2\lambda}$ on this branch,
$\big|\tau\varphi'(x_1)-\psi\,\mathrm{sgn}(x_1-y_1)\big|\,\Theta_\tau\leq C$. This proves
\eqref{pw-d1k0}.

\textit{Proof of \eqref{pw-d1k1}.} For $y_1\neq x_1$, differentiating \eqref{Def-k-tau-1-x1-y1-xi'} with respect to $x_1$ multiplies its first two terms by $\psi(\xi')+\tau\varphi'(x_1)$, of modulus at most $C(\lambda+\tau)$ by \eqref{psi-modulus}, and turns the last one into $-\big(\psi(\xi')-\tau\varphi'(y_1)\big)\partial_1k_{\tau,0}(x_1,y_1,\xi')$. Hence \eqref{pw-d1k1} follows from the proof of \eqref{pw-k1} and from \eqref{pw-d1k0}.
\end{proof}
\begin{remark}\label{Rem-comparison-DET}
The bounds of Lemma \ref{Lemma-kernels-bounds} have the same form as those of \cite[Lemma 6.4]{Dehman-Ervedoza-Thabouti-2023}, the turning level $\tau$ there being replaced by $-\tau\varphi'(x_1)\in[\alpha_1\tau,\alpha_2\tau]$, and the Gaussian factor $e^{-\lambda(x_1-y_1)^2/C_1}$ by $e^{-\frac\beta2\tau(x_1-y_1)^2}$, which is comparable to it near the turning level. There is one difference: in \cite[(6.7)--(6.8)]{Dehman-Ervedoza-Thabouti-2023}, the kernel of $\partial_1K_{\tau,0}$ carries a factor $|\lambda-\tau|$, which vanishes at the turning level, whereas here the corresponding factor $\psi+\tau\varphi'(x_1)$ equals $i\operatorname{Im}\psi=i\tau\,\ell\cdot\xi'/\lambda$ at the turning level, whose supremum on $\Sigma_\lambda$ is $\tau|\ell|$. This is why \eqref{pw-d1k0} carries no gain, and why $\partial_1K_{\tau,j}$ is estimated together with $\nabla'K_{\tau,j}$ in Proposition \ref{Prop-K-high}, instead of separately as in \cite[Propositions 6.2 and 6.7]{Dehman-Ervedoza-Thabouti-2023}.
\end{remark}
\begin{proof}[Proof of Proposition \ref{Prop-K-high}] Our starting point is that each of the operators is of the form
\begin{equation}\label{K-as-integral}
K_{\tau,j}\Pi_{\mathrm{high}}f(x_1,\cdot)
=\int_0^L k_{\tau,j}(x_1,y_1,D')\,\Pi_{\mathrm{high}}f(y_1,\cdot)\,dy_1,
\qquad j\in\{0,\dots,d\},
\end{equation}
where, for each fixed pair $(x_1,y_1)\in[0,L]^2$, the operator
$k_{\tau,j}(x_1,y_1,D')\,\Pi_{\mathrm{high}}$ is precisely a Fourier multiplier of the class
\eqref{Def-Fourier-mulptilier-K}, associated with the function
\begin{equation*}
    k(\lambda,\omega'):=k_{\tau,j}\big(x_1,y_1,\lambda\omega'\big),
\qquad \lambda\geq2\tau|\ell|,\quad\omega'\in\Sigma_{\tau|\ell|/\lambda,\,e_\ell},
\end{equation*}
in the parametrization \eqref{Def-polar-coordinates}. Equivalently, by \eqref{Def-E-lambda},
\begin{equation}\label{K-as-layers}
k_{\tau,j}(x_1,y_1,D')\,\Pi_{\mathrm{high}}\,g
=\int_{2\tau|\ell|}^{\infty}E_\lambda\big[k_{\tau,j}(x_1,y_1,D')g\big]\,d\lambda,
\end{equation}
so that only the values of the kernel on the layers $\Sigma_\lambda$ with $\lambda\geq2\tau|\ell|$ are
involved, and the relevant quantity is
$\big\|k(\lambda,\cdot)\big\|_{L^\infty}=\big\|k_{\tau,j}(x_1,y_1,\lambda,\cdot)\big\|_{L^\infty(\Sigma_\lambda)}$,
which is exactly what Lemma \ref{Lemma-kernels-bounds} estimates. Applying Minkowski's integral inequality in the variable $y_1$ and then Proposition \ref{Lemma-E-lambda-multiplier} for each fixed pair $(x_1,y_1)$, we obtain
\begin{equation}\label{reduction-1d}
\big\|K_{\tau,j}\Pi_{\mathrm{high}}f(x_1,\cdot)\big\|_{L^q(\R^{d-1})}
\leq C\int_0^L\mathcal N_j\big(x_1,y_1\big)\,\big\|f(y_1,\cdot)\big\|_{L^p(\R^{d-1})}\,dy_1,
\end{equation}
where, according to the pair $(p,q)$ under consideration,
\begin{align*}
 \mathcal N_j=&\sup_{\lambda\geq2\tau|\ell|}\big\|k_{\tau,j}(x_1,y_1,\lambda,\cdot)\big\|_{L^\infty(\Sigma_\lambda)}
\quad\text{if }(p,q)=(2,2),
\\ 
\mathcal N_j=&\Big(\int_{2\tau|\ell|}^{\infty}\big\|k_{\tau,j}(x_1,y_1,\lambda,\cdot)\big\|^r_{L^\infty(\Sigma_\lambda)}
\lambda^{1-\frac{2}{d}}\,d\lambda\Big)^{\frac1r},   
\end{align*}
with $r=1$ if $(p,q)=\big(\tfrac{2d}{d+2},\tfrac{2d}{d-2}\big)$ and $r=2$ in the two mixed cases. 

For $(p,q)=(2,2)$, \eqref{reduction-1d} is too lossy for $\nabla K_{\tau,j}$, $j\geq1$, since the supremum in $\lambda$ is taken before integrating in $y_1$: the supremum in $\lambda$ of $\lambda^2\Theta_\tau\mathcal E_\tau$ is of order $|x_1-y_1|^{-1}$ for $|x_1-y_1|\leq\tau^{-1}$. We shall therefore use it only for $K_{\tau,0}$, and, for the other operators, Plancherel's identity in $x'$ and Young's inequality in $x_1$ for each fixed $\xi'\in\mathcal R_{\mathrm{high}}$: if the kernel $k$ of an operator $T$ satisfies $|k(x_1,y_1,\xi')|\leq m_\lambda(x_1-y_1)$ for all $(x_1,y_1)\in[0,L]^2$, with $\lambda=\operatorname{Re}\psi(\xi')$, then
\begin{equation}\label{reduction-L2}
\big\|T\,\Pi_{\mathrm{high}}f\big\|_{L^2(\mathcal S)}
\leq\sup_{\lambda\geq2\tau|\ell|}\big\|m_\lambda\big\|_{L^1(\R)}\,\|f\|_{L^2(\mathcal S)}.
\end{equation}

The same reductions apply to the derivatives, the tangential ones amount to replacing $k_{\tau,j}$ by $i\xi'k_{\tau,j}$, whose $L^\infty(\Sigma_\lambda)$ norm is bounded by
$\tfrac32\lambda\|k_{\tau,j}(x_1,y_1,\lambda,\cdot)\|_{L^\infty(\Sigma_\lambda)}$ by \eqref{psi-modulus}, while the normal one amounts to replacing $k_{\tau,j}$ by $\partial_1k_{\tau,j}$, estimated in
\eqref{pw-d1k0}, \eqref{pw-d1k1} and, for $j\geq2$, by $\partial_1k_{\tau,j}=i\xi_j\partial_1k_{\tau,0}$. Here, as in \cite[(4.6) and Remark 6.5]{Dehman-Ervedoza-Thabouti-2023}, $\partial_1k_{\tau,j}$ is computed for $\lambda\neq-\tau\varphi'(x_1)$, the set $\{(x_1,\xi'):\operatorname{Re}\psi(\xi')=-\tau\varphi'(x_1)\}$ being negligible; this is enough for our purpose since, $w$ being in $H^1(\mathcal S)$, the jumps of the terms of \eqref{parametrix} across this set cancel. Moreover, for $j=1$, since $k_{\tau,1}(x_1,x_1^-,\xi')-k_{\tau,1}(x_1,x_1^+,\xi')=1$ on both branches, differentiating \eqref{Def-K-tau-j-Fourier} with respect to $x_1$ gives $\partial_1K_{\tau,1}\Pi_{\mathrm{high}}f=\Pi_{\mathrm{high}}f+K_{\tau,1,\partial_1}\Pi_{\mathrm{high}}f$, where $K_{\tau,1,\partial_1}$ has kernel $\partial_1k_{\tau,1}$ computed for $y_1\neq x_1$.
\medskip 

The next step is to estimate the quantities $\mathcal N_j$. By Lemma \ref{Lemma-kernels-bounds} and
Remark \ref{Rem-comparison-DET}, the quantities $\mathcal N_j$ obey the same bounds as their counterparts in \cite[Lemma 6.6]{Dehman-Ervedoza-Thabouti-2023}: for all $(x_1,y_1)\in[0,L]^2$,
\begin{align}
&\sup_{\lambda\geq2\tau|\ell|}\big\|k_{\tau,0}(x_1,y_1,\lambda,\cdot)\big\|_{L^\infty(\Sigma_\lambda)}
\leq\frac{C}{\tau}\,e^{-\frac{\beta}{2}\tau(x_1-y_1)^2},
\label{N-sup}\\
&\int_{2\tau|\ell|}^{\infty}\big\|k_{\tau,0}(x_1,y_1,\lambda,\cdot)\big\|_{L^\infty(\Sigma_\lambda)}
\lambda^{1-\frac{2}{d}}\,d\lambda\ \leq\ \frac{C}{|x_1-y_1|^{1-\frac{2}{d}}},
\label{N-L1}\\
&\Big(\int_{2\tau|\ell|}^{\infty}\big\|k_{\tau,0}(x_1,y_1,\lambda,\cdot)\big\|^2_{L^\infty(\Sigma_\lambda)}
\lambda^{1-\frac{2}{d}}\,d\lambda\Big)^{\frac{1}{2}}\leq\widetilde k_{\tau,0}(x_1-y_1),\label{N-L2}\\
			& \qquad \text{ with $\widetilde k_{\tau,0} \in L^{\frac{d}{d-1}} (\R)$ and }
\big\|\widetilde k_{\tau,0}\big\|_{L^{\frac d{d-1}}(\R)}\leq C\,\tau^{-\frac{3}{4}-\frac1{2d}}, \notag 
\\
& \Big(\int_{2\tau|\ell|}^{\infty}\Big(\big\|\partial_1k_{\tau,0}\big\|_{L^\infty(\Sigma_\lambda)}
+\big\|\lambda\,k_{\tau,0}\big\|_{L^\infty(\Sigma_\lambda)}\Big)^2\lambda^{1-\frac{2}{d}}\,d\lambda\Big)^{\frac{1}{2}}
\leq\frac{C}{|x_1-y_1|^{1-\frac1d}}+\check k_{\tau,0}(x_1-y_1),\label{N-lambda-L2}
\\
			& \qquad \text{ with $\check k_{\tau,0} \in L^{\frac{d}{d-1}} (\R)$ and }
			\|\check k_{\tau,0}\|_{L^{\frac{d}{d-1}} (\R)} \leq C \tau^{\frac{1}{4}-\frac{1}{2d}},
			\notag
\end{align}
\eqref{N-sup} following directly from \eqref{pw-k0}, and the other proofs being those of \cite[Lemma 6.6]{Dehman-Ervedoza-Thabouti-2023}, which use only the form of the bounds of Lemma \ref{Lemma-kernels-bounds}, with $\tau$ replaced by $-\tau\varphi'(x_1)\in[\alpha_1\tau,\alpha_2\tau]$, so that the resulting bounds are uniform in $x_1\in[0,L]$. We therefore do not reproduce them. By
\eqref{pw-kj} and \eqref{pw-k1}, and since $\tau\Theta_\tau\leq\frac1{\alpha_1}$ for $\lambda>\frac{\alpha_1}2\tau$ while $\tau\Theta_\tau\mathcal E_\tau\leq Ce^{-\frac{\alpha_1}8\tau|x_1-y_1|}$ for $\lambda\leq\frac{\alpha_1}2\tau$, the kernels $k_{\tau,j}$, $j\in\{1,\dots,d\}$, satisfy
\eqref{N-lambda-L2} as well.

For \eqref{reduction-L2}, we distinguish three ranges of $\lambda$, independent of $x_1$: if $\lambda\leq\frac{\alpha_1}2\tau$, then $\lambda<-\tau\varphi'(x_1)$ and $|\lambda+\tau\varphi'(x_1)|\geq\frac{\alpha_1}2\tau$; if $\frac{\alpha_1}2\tau<\lambda<2\alpha_2\tau$, then $\Theta_\tau\leq\frac1{2\lambda}<\frac1{\alpha_1\tau}$; and if $\lambda\geq2\alpha_2\tau$, then $\Theta_\tau=\frac1{2\lambda}$ and $|\lambda+\tau\varphi'(x_1)|\geq\frac\lambda2$, for all $x_1\in[0,L]$. Hence
$\Theta_\tau\mathcal E_\tau\leq\kappa_\lambda(x_1-y_1)$ and $\mathcal E_\tau\leq e_\lambda(x_1-y_1)$, where $\kappa_\lambda(z)$ is
$|z|e^{-\frac{\alpha_1}4\tau|z|}$, $(\alpha_1\tau)^{-1}e^{-\frac{\beta}2\tau z^2}$ and $(2\lambda)^{-1}e^{-\frac\lambda4|z|}$, and $e_\lambda(z)$ is
$e^{-\frac{\alpha_1}4\tau|z|}$, $e^{-\frac{\beta}2\tau z^2}$ and $e^{-\frac\lambda4|z|}$, in these three ranges respectively. Computing their $L^1(\R)$ norms gives, for $a\in\{1,2\}$ and $b\in\{0,1\}$,
\begin{equation}\label{kappa-e-L1}
\sup_{\lambda\geq2\tau|\ell|}\big\|(\lambda+\tau)^a\kappa_\lambda\big\|_{L^1(\R)}\leq C\tau^{a-\frac{3}{2}},
\qquad
\sup_{\lambda\geq2\tau|\ell|}\big\|(\lambda+\tau)^b e_\lambda\big\|_{L^1(\R)}\leq C\tau^{b-\frac{1}{2}}.
\end{equation}
\medskip

Finally, it remains to insert \eqref{N-sup}--\eqref{kappa-e-L1} into
\eqref{reduction-1d} and \eqref{reduction-L2}, and to apply, in the variable $x_1$, Young's inequality or the Hardy--Littlewood--Sobolev inequality, exactly as in
\cite[Propositions 6.2, 6.7 and 6.10]{Dehman-Ervedoza-Thabouti-2023}. The proof is relatively standard, but for
the sake of completeness we outline it, leaving the details to the reader:  
\begin{itemize}
    \item \eqref{N-sup} together with
$\|\tau^{-1}e^{-\frac{\beta}{2}\tau\,\cdot^2}\|_{L^1(\R)}\leq C\tau^{-\frac{3}{2}}$ and Young's inequality
$L^1*L^2\to L^2$ gives the $L^2\to L^2$ estimate of $K_{\tau,0}$ in \eqref{K-high-1}. By Lemma \ref{Lemma-kernels-bounds} and \eqref{psi-modulus}, the kernels of $\nabla K_{\tau,0}$ and $K_{\tau,j}$, and of $\nabla K_{\tau,j}$ (without the term $\Pi_{\mathrm{high}}f$ of $\partial_1K_{\tau,1}\Pi_{\mathrm{high}}f$), are bounded by $C\big(e_\lambda+(\lambda+\tau)\kappa_\lambda\big)$ and $C\big((\lambda+\tau)e_\lambda+(\lambda+\tau)^2\kappa_\lambda\big)$ respectively, so that \eqref{kappa-e-L1} and \eqref{reduction-L2} give the other $L^2\to L^2$ estimates of \eqref{K-high-1} and \eqref{K-high-3}. 
\item \eqref{N-L1} together with the Hardy--Littlewood--Sobolev inequality
\cite[Theorem 0.3.2]{sogge2017fourier} gives the $L^{\frac{2d}{d+2}}\to L^{\frac{2d}{d-2}}$ estimate of \eqref{K-high-2}. 
\item \eqref{N-L2} together with Young's inequality $L^{\frac d{d-1}}*L^2\subset L^{\frac{2d}{d-2}}$ gives the $L^{2}\to L^{\frac{2d}{d-2}}$ estimate of
\eqref{K-high-1} and, with $L^{\frac d{d-1}}*L^{\frac{2d}{d+2}}\subset L^2$, the $L^{\frac{2d}{d+2}}\to L^2$ estimate of $K_{\tau,0}$ in \eqref{K-high-2}. 
\item  \eqref{N-lambda-L2}, whose singular part $|x_1-y_1|^{-1+\frac1d}$ is handled by the Hardy--Littlewood--Sobolev inequality and whose regular part
$\check k_{\tau,0}$ by Young's inequality, gives the $L^{\frac{2d}{d+2}}\to L^2$ estimate of $\nabla K_{\tau,0}$ in \eqref{K-high-2} and the first term of \eqref{K-high-3}.
\end{itemize}
\end{proof}

We finally estimate the boundary operators $G_\tau$ \eqref{Def-G-tau-Fourier} and $H_\tau$
\eqref{Def-H-tau-Fourier}.
\begin{proposition}\label{Prop-GH-high}
There exists $C>0$, depending only on $d,\alpha_1,\alpha_2,\beta$, such that for all $\tau\geq1$ and all
$|\ell|\leq\ell_*$, one has, for all $g\in H^{\frac{1}{2}}(\R^{d-1})$,
\begin{equation}\label{Est-G-tau-high}
\begin{gathered}
\tau^{\frac{3}{4}}\big\|G_\tau\Pi_{\mathrm{high}}g\big\|_{L^2(\mathcal S)}
+\tau^{-\frac{1}{4}}\big\|\nabla G_\tau\Pi_{\mathrm{high}}g\big\|_{L^2(\mathcal S)} 
\leq C\left(\|g\|_{H^{\frac{1}{2}}(\R^{d-1})}+ \tau^{\frac{1}{2}} \|g\|_{L^2(\R^{d-1})}\right),\\
 \big\|G_\tau\Pi_{\mathrm{high}}g\big\|_{L^{\frac{2d}{d-2}}(\mathcal S)}
\leq C\|g\|_{H^{\frac{1}{2}}(\R^{d-1})},
\end{gathered}
\end{equation}
and, for all $\theta\in\{0,-\frac{1}{2}\}$ and all $h\in H^{\theta}(\R^{d-1})$,
\begin{equation}\label{Est-H-tau-high}
\tau^{\frac{5}{4}}\big\|H_\tau\Pi_{\mathrm{high}}h\big\|_{L^2(\mathcal S)}
+\tau^{\frac{1}{4}}\big\|\nabla H_\tau\Pi_{\mathrm{high}}h\big\|_{L^2(\mathcal S)}
+\tau^{\frac{1}{2}}\big\|H_\tau\Pi_{\mathrm{high}}h\big\|_{L^{\frac{2d}{d-2}}(\mathcal S)}
\leq C\,\tau^{-\theta}\,\|h\|_{H^{\theta}(\R^{d-1})} .
\end{equation}
\end{proposition}
\begin{proof}
As in the previous proposition, both operators are Fourier multipliers in $x'$, of the class
\eqref{Def-Fourier-mulptilier-K} for each fixed $x_1\in[0,L]$, associated with the functions
$(\lambda,\omega')\mapsto g_\tau(x_1,\lambda\omega')$ and
$(\lambda,\omega')\mapsto k_{\tau,0}(x_1,0,\lambda\omega')$ respectively, and, as there, $\partial_1$ is computed for $\lambda\neq-\tau\varphi'(x_1)$. We first recall the pointwise
bounds \eqref{pw-boundary},
\begin{equation}\label{pw-boundary-recall}
\big\|k_{\tau,0}(x_1,0,\lambda,\cdot)\big\|_{L^\infty(\Sigma)}\leq C\,\Theta_\tau\,\mathcal E_\tau,
\qquad
\big\|g_\tau(x_1,\lambda,\cdot)\big\|_{L^\infty(\Sigma)}
\leq C\big(1+(\lambda+\tau)\Theta_\tau\big)\,\mathcal E_\tau,
\end{equation}
where $\Theta_\tau=\Theta_\tau(x_1,0,\lambda)$ and $\mathcal E_\tau=\mathcal E_\tau(x_1,0,\lambda)$ are
given by \eqref{Def-Theta-E}. Note that, by \eqref{Def-Theta-E}, $(\lambda+\tau)\Theta_\tau\leq C$ on
$\{\lambda\geq\frac{1}{2}(-\tau\varphi'(x_1))\}$ and, on the complementary set,
$(\lambda+\tau)\Theta_\tau\leq C|\lambda+\tau\varphi'(x_1)|\,x_1$, since there $\Theta_\tau\leq x_1$ and
$|\lambda+\tau\varphi'(x_1)|\geq\frac{1}{2}(-\tau\varphi'(x_1))\geq\max\{\lambda,\frac{\alpha_1}2\tau\}$. Hence
$(1+u)e^{-\frac u4}\leq C$, with $u=|\lambda+\tau\varphi'(x_1)|\,x_1$, gives in both cases
\begin{equation}\label{g-tau-bounded}
\big\|g_\tau(x_1,\lambda,\cdot)\big\|_{L^\infty(\Sigma)}\leq C\,
e^{-\frac{1}{4}|\lambda+\tau\varphi'(x_1)|\,x_1-\frac{\beta}{2}\tau x_1^2} .
\end{equation}

\medskip
\textit{Proof of the estimates \eqref{Est-H-tau-high} on $H_\tau$.}
We begin with the Hilbertian estimates. From \eqref{pw-boundary-recall}, \eqref{Est-K-k-L2} and
\eqref{N-sup} at $y_1=0$, we obtain
$\|H_\tau\Pi_{\mathrm{high}}h\|_{L^2(\mathcal S)}\leq C\tau^{-\frac{5}{4}}\|h\|_{L^2(\R^{d-1})}$, which
gives the first term of \eqref{Est-H-tau-high} for $\theta=0$, and for $\theta=-\frac{1}{2}$ by
$\|\Pi_{\mathrm{high}}h\|_{L^2}\leq C\tau^{\frac{1}{2}}\|h\|_{H^{-\frac{1}{2}}}$ on the region
$\lambda\lesssim\tau$, the region $\lambda\gtrsim\tau$ being treated by the sharper bound
$\Theta_\tau\leq\frac{1}{2\lambda}$. The estimate on the gradient follows in the same way for $\theta=0$, the tangential
derivative costing a factor $|\xi'|\leq\frac{3}{2}\lambda$ and the normal one being handled by
\eqref{pw-d1k0}. For $\theta=-\frac{1}{2}$, since \eqref{pw-d1k0} carries no factor $\frac1\lambda$, we use instead Plancherel's identity in $x'$,
\begin{equation*}
\big\|\nabla H_\tau\Pi_{\mathrm{high}}h\big\|^2_{L^2(\mathcal S)}
=\int_{\mathcal R_{\mathrm{high}}}\big\|\big(\partial_1k_{\tau,0},\,i\xi'k_{\tau,0}\big)(\cdot,0,\xi')\big\|^2_{L^2(0,L)}\,\big|\widehat h(\xi')\big|^2\,d\xi',
\end{equation*}
together with $|\partial_1k_{\tau,0}(x_1,0,\xi')|+|\xi'||k_{\tau,0}(x_1,0,\xi')|\leq C\mathcal E_\tau(x_1,0,\lambda)$, which follows from \eqref{pw-d1k0}, \eqref{pw-boundary-recall}, $|\xi'|\leq\frac{3}{2}\lambda$ and $\lambda\Theta_\tau\leq\frac{1}{2}$, and with the bound, valid for $c\in\{\frac14,\frac{1}{2}\}$ and all $\lambda\geq2\tau|\ell|$,
\begin{equation}\label{E-tau-L2}
\Big\|e^{-c|\lambda+\tau\varphi'(x_1)|\,x_1-\frac{\beta}{2}\tau x_1^2}\Big\|_{L^2_{x_1}(0,L)}\leq C\,\tau^{\frac14}(\lambda+\tau)^{-\frac{1}{2}},
\end{equation}
obtained by distinguishing, as in the proof of Proposition \ref{Prop-K-high}, the ranges $\lambda\leq\frac{\alpha_1}2\tau$, $\frac{\alpha_1}2\tau<\lambda<2\alpha_2\tau$ and $\lambda\geq2\alpha_2\tau$, which give respectively at most $C\tau^{-\frac{1}{2}}$, $C\tau^{-\frac14}$ and $C\lambda^{-\frac{1}{2}}$. Since $1+|\xi'|\leq C(\lambda+\tau)$ by \eqref{psi-modulus}, we get $\|\nabla H_\tau\Pi_{\mathrm{high}}h\|_{L^2(\mathcal S)}\leq C\tau^{\frac14}\|h\|_{H^{-\frac{1}{2}}(\R^{d-1})}$, which is the second term of \eqref{Est-H-tau-high} for $\theta=-\frac{1}{2}$.

We now prove the $L^2-L^{\frac{2d}{d-2}}$ estimate, which corresponds to the case $\theta=0$. It suffices
to apply \eqref{Est-K-k-mixed} to the multiplier $k_{\tau,0}(x_1,0,\cdot)$, for each fixed $x_1\in[0,L]$,
which gives
\begin{equation*}
\big\|H_\tau\Pi_{\mathrm{high}}h(x_1,\cdot)\big\|_{L^{\frac{2d}{d-2}}(\R^{d-1})}
\leq C\left(\int_{2\tau|\ell|}^{\infty}
\big\|k_{\tau,0}(x_1,0,\lambda,\cdot)\big\|^2_{L^\infty(\Sigma)}\,\lambda^{1-\frac{2}{d}}\,d\lambda
\right)^{\frac{1}{2}}\,\|h\|_{L^2(\R^{d-1})} .
\end{equation*}
Since, by \eqref{pw-boundary-recall} and as for \eqref{N-L2} at $y_1=0$,
\begin{equation*}
\left(\int_{2\tau|\ell|}^{\infty}\big\|k_{\tau,0}(x_1,0,\lambda,\cdot)\big\|^2_{L^\infty(\Sigma)}\,
\lambda^{1-\frac{2}{d}}\,d\lambda\right)^{\frac{1}{2}}
\leq C\,\tilde k_{\tau,0}(x_1),
\end{equation*}
where the kernel $\tilde k_{\tau,0}$ can be taken, as in
\cite[(6.14)]{Dehman-Ervedoza-Thabouti-2023}, of the form
\begin{equation*}
\tilde k_{\tau,0}(z_1)=
\mathbf 1_{|z_1|<1/\tau}\,\tau^{-\frac{1}{d}}
+\mathbf 1_{|z_1|>1/\tau}\,\tau^{-\frac{1}{d}}\left((\tau|z_1|)^{\frac{1}{d}}e^{-\tau|z_1|/C}
+\frac{e^{-\tau|z_1|^2/C}}{(\tau|z_1|)^{\frac{1}{2}}}\right).
\end{equation*}
 Easy computations then give $\|\tilde k_{\tau,0}\|_{L^{\frac{2d}{d-2}}(\R)}\leq C\tau^{-\frac{1}{2}}$, which concludes the
$L^2-L^{\frac{2d}{d-2}}$ bound for the operator $H_\tau$.

There remains the case $\theta=-\frac{1}{2}$, for which we use the $TT^\star$ argument, following
\cite[Proposition 2.5.14]{Thabouti2024}. Denote by $(I+|D_{x'}|)^{\frac{1}{2}}$ the Fourier multiplier by
$(1+|\xi'|)^{\frac{1}{2}}$ and set
$\widetilde H_\tau:=H_\tau\Pi_{\mathrm{high}}\circ(I+|D_{x'}|)^{\frac{1}{2}}$, so that
\begin{equation*}
\widehat{\widetilde H_\tau h_0}(x_1,\xi')=\mathbf 1_{\mathcal R_{\mathrm{high}}}(\xi')\,
(1+|\xi'|)^{\frac{1}{2}}\,k_{\tau,0}(x_1,0,\xi')\,\widehat{h_0}(\xi') .
\end{equation*}
The $\mathcal L\big(H^{-\frac{1}{2}}(\R^{d-1}),L^{\frac{2d}{d-2}}(\mathcal S)\big)$-norm of
$H_\tau\Pi_{\mathrm{high}}$ then coincides with the
$\mathcal L\big(L^2(\R^{d-1}),L^{\frac{2d}{d-2}}(\mathcal S)\big)$-norm of $\widetilde H_\tau$, hence, by
duality, with the square root of the
$\mathcal L\big(L^{\frac{2d}{d+2}}(\mathcal S),L^{\frac{2d}{d-2}}(\mathcal S)\big)$-norm of
$\widetilde H_\tau\widetilde H_\tau^\star$, and it therefore suffices to prove that the latter is bounded
independently of $\tau$.

Recalling from \eqref{Def-k-tau-0-x1-y1-xi'} that
$k_{\tau,0}(x_1,0,\cdot)=\mathbf 1_{\lambda<-\tau\varphi'(x_1)}k_{\tau,01}(x_1,0,\cdot)
+\mathbf 1_{\lambda\geq-\tau\varphi'(x_1)}k_{\tau,02}(x_1,0,\cdot)$, where $k_{\tau,01}$ and $k_{\tau,02}$
denote the two lines of \eqref{Def-k-tau-0-x1-y1-xi'}, a direct computation gives, for
$v\in L^{\frac{2d}{d+2}}(\mathcal S)$, writing $m(D')^\star$ for the Fourier multiplier by $\overline{m(\xi')}$,
\begin{multline}\label{TTstar-H}
\widetilde H_\tau\widetilde H_\tau^\star v(x_1,x')
=\int_0^L\int_{2\tau|\ell|}^{-\tau\max\{\varphi'(x_1),\varphi'(z_1)\}}
E_\lambda\Big[k_{\tau,01}(x_1,0,D')\,k_{\tau,01}(z_1,0,D')^\star\,(1+|D'|)\,
v(z_1,\cdot)\Big](x')\,d\lambda\,dz_1\\
+\int_0^L\int_{-\tau\min\{\varphi'(x_1),\varphi'(z_1)\}}^{+\infty}
E_\lambda\Big[k_{\tau,02}(x_1,0,D')\,k_{\tau,02}(z_1,0,D')^\star\,(1+|D'|)\,
v(z_1,\cdot)\Big](x')\,d\lambda\,dz_1\\
+\int_0^L\int_{-\tau\max\{\varphi'(x_1),\varphi'(z_1)\}}^{-\tau\min\{\varphi'(x_1),\varphi'(z_1)\}}
E_\lambda\Big[k_{\tau,0}(x_1,0,D')\,k_{\tau,0}(z_1,0,D')^\star\,(1+|D'|)\,
v(z_1,\cdot)\Big](x')\,d\lambda\,dz_1 ,
\end{multline}
the third term corresponding to the region of frequencies lying between the two turning levels, on which
the multiplier at $\min\{x_1,z_1\}$ is on the forward branch, given by $k_{\tau,01}$, and the one at
$\max\{x_1,z_1\}$ on the backward branch, given by $k_{\tau,02}$. Here we used that
$\varphi'$ is increasing, so that
$-\tau\varphi'(\max\{x_1,z_1\})=-\tau\max\{\varphi'(x_1),\varphi'(z_1)\}$. Note also that, on
$\mathcal R_{\mathrm{high}}$, one has $1+|\xi'|\leq C(1+\lambda)$ by \eqref{psi-modulus}. Applying
\eqref{Est-K-k-Lp-Lp'} of Proposition \ref{Lemma-E-lambda-multiplier} to each of the three terms of
\eqref{TTstar-H} and Minkowski's integral inequality in $z_1$, we obtain
\begin{multline}\label{TTstar-H-bound}
\big\|\widetilde H_\tau\widetilde H_\tau^\star v(x_1,\cdot)\big\|_{L^{\frac{2d}{d-2}}(\R^{d-1})}\\
\leq C\int_0^L\Bigg(\int_{0}^{-\tau\max\{\varphi'(x_1),\varphi'(z_1)\}}
\big\|k_{\tau,01}(x_1,0,\lambda,\cdot)\big\|_{L^\infty(\Sigma)}
\big\|k_{\tau,01}(z_1,0,\lambda,\cdot)\big\|_{L^\infty(\Sigma)}(1+\lambda)\lambda^{1-\frac{2}{d}}\,d\lambda\\
+\int_{-\tau\min\{\varphi'(x_1),\varphi'(z_1)\}}^{+\infty}
\big\|k_{\tau,02}(x_1,0,\lambda,\cdot)\big\|_{L^\infty(\Sigma)}
\big\|k_{\tau,02}(z_1,0,\lambda,\cdot)\big\|_{L^\infty(\Sigma)}(1+\lambda)\lambda^{1-\frac{2}{d}}\,d\lambda\\
+\int_{-\tau\max\{\varphi'(x_1),\varphi'(z_1)\}}^{-\tau\min\{\varphi'(x_1),\varphi'(z_1)\}}
\big\|k_{\tau,01}(\min\{x_1,z_1\},0,\lambda,\cdot)\big\|_{L^\infty(\Sigma)}
\big\|k_{\tau,02}(\max\{x_1,z_1\},0,\lambda,\cdot)\big\|_{L^\infty(\Sigma)}
(1+\lambda)\lambda^{1-\frac{2}{d}}\,d\lambda\Bigg)\\
\times\big\|v(z_1,\cdot)\big\|_{L^{\frac{2d}{d+2}}(\R^{d-1})}\,dz_1 .
\end{multline}
By \eqref{pw-boundary}, the norms $\|k_{\tau,0i}(x_1,0,\lambda,\cdot)\|_{L^\infty(\Sigma)}$,
$i\in\{1,2\}$, are bounded by $C\,\Theta_\tau(x_1,0,\lambda)\,\mathcal E_\tau(x_1,0,\lambda)$, which are of
the same form as the bounds on $k_{1,\tau}$ and $k_{2,\tau}$ used in
\cite[Lemme 2.5.19]{Thabouti2024}, where the turning level also depends on $x_1$. Consequently, the proof of \cite[Lemme 2.5.19]{Thabouti2024} yields that the bracket
in \eqref{TTstar-H-bound} is bounded by $C\,|x_1-z_1|^{-1+\frac{2}{d}}$, whence
\begin{equation*}
\big\|\widetilde H_\tau\widetilde H_\tau^\star v(x_1,\cdot)\big\|_{L^{\frac{2d}{d-2}}(\R^{d-1})}
\leq C\int_0^L|x_1-z_1|^{-1+\frac{2}{d}}\,\big\|v(z_1,\cdot)\big\|_{L^{\frac{2d}{d+2}}(\R^{d-1})}\,dz_1 .
\end{equation*}
The Hardy--Littlewood--Sobolev theorem in the $x_1$-variable then implies
\begin{equation*}
\big\|\widetilde H_\tau\widetilde H_\tau^\star v\big\|_{L^{\frac{2d}{d-2}}(\mathcal S)}
\leq C\,\|v\|_{L^{\frac{2d}{d+2}}(\mathcal S)} .
\end{equation*}
It follows that the $\mathcal L\big(L^2(\R^{d-1}),L^{\frac{2d}{d-2}}(\mathcal S)\big)$-norm of
$\widetilde H_\tau$ is bounded independently of $\tau$, as announced, and thus this is also the case for
the $\mathcal L\big(H^{-\frac{1}{2}}(\R^{d-1}),L^{\frac{2d}{d-2}}(\mathcal S)\big)$-norm of
$H_\tau\Pi_{\mathrm{high}}$. This concludes the proof of \eqref{Est-H-tau-high}.

\medskip
\textit{Proof of the estimates \eqref{Est-G-tau-high} on $G_\tau$.}
The Hilbertian estimates follow from \eqref{g-tau-bounded} and \eqref{Est-K-k-L2}, which give
\begin{equation*}
\big\|G_\tau\Pi_{\mathrm{high}}g(x_1,\cdot)\big\|_{L^2(\R^{d-1})}
\leq C\,e^{-\frac{\beta}{2}\tau x_1^2}\,\|g\|_{L^2(\R^{d-1})},
\quad\text{whence}\quad
\big\|G_\tau\Pi_{\mathrm{high}}g\big\|_{L^2(\mathcal S)}\leq C\tau^{-\frac{1}{4}}\|g\|_{L^2(\R^{d-1})},
\end{equation*}
and hence the first term of \eqref{Est-G-tau-high}. For the gradient, differentiating \eqref{Def-g-tau} and using
\eqref{psi-modulus}, \eqref{exponent-forward} and \eqref{pw-d1k0} at $y_1=0$, together with \eqref{g-tau-bounded}, we get
$|\partial_1g_\tau(x_1,\xi')|+|\xi'||g_\tau(x_1,\xi')|\leq C(\lambda+\tau)\,e^{-\frac{1}{4}|\lambda+\tau\varphi'(x_1)|\,x_1-\frac{\beta}{2}\tau x_1^2}$,
so that Plancherel's identity in $x'$ and \eqref{E-tau-L2} give, with $\lambda=\operatorname{Re}\psi(\xi')$,
\begin{equation*}
\big\|\nabla G_\tau\Pi_{\mathrm{high}}g\big\|_{L^2(\mathcal S)}
\leq C\tau^{\frac14}\big\|(\lambda+\tau)^{\frac{1}{2}}\,\widehat g\big\|_{L^2(\R^{d-1})}
\leq C\tau^{\frac14}\Big(\|g\|_{H^{\frac{1}{2}}(\R^{d-1})}+\tau^{\frac{1}{2}}\|g\|_{L^2(\R^{d-1})}\Big),
\end{equation*}
which is the second term of \eqref{Est-G-tau-high}.

We now prove the $H^{\frac{1}{2}}-L^{\frac{2d}{d-2}}$ estimate, again by the $TT^\star$ argument. Set
$\widetilde G_\tau:=G_\tau\Pi_{\mathrm{high}}\circ(I+|D_{x'}|)^{-\frac{1}{2}}$, so that
\begin{equation*}
\widehat{\widetilde G_\tau g_0}(x_1,\xi')=g_{\tau,0}(x_1,\xi')\,\widehat{g_0}(\xi'),
\qquad
g_{\tau,0}(x_1,\xi'):=\mathbf 1_{\mathcal R_{\mathrm{high}}}(\xi')\,
\frac{g_\tau(x_1,\xi')}{(1+|\xi'|)^{\frac{1}{2}}} ,
\end{equation*}
and the $\mathcal L\big(H^{\frac{1}{2}},L^{\frac{2d}{d-2}}(\mathcal S)\big)$-norm of
$G_\tau\Pi_{\mathrm{high}}$ coincides with the
$\mathcal L\big(L^2,L^{\frac{2d}{d-2}}(\mathcal S)\big)$-norm of $\widetilde G_\tau$, hence with the square
root of the $\mathcal L\big(L^{\frac{2d}{d+2}}(\mathcal S),L^{\frac{2d}{d-2}}(\mathcal S)\big)$-norm of
$\widetilde G_\tau\widetilde G_\tau^\star$.

Note that, although $g_\tau$ carries the branch indicators at the level $-\tau\varphi'(x_1)$, so that the
product of the multipliers at two distinct heights splits into three regions of $\lambda$ as in
\eqref{TTstar-H}, the bound \eqref{g-tau-bounded} is, in contrast with $k_{\tau,0}$, the same on the two
branches; this splitting therefore need not be carried through the estimates. Computing the adjoint with
respect to the $L^2(\mathcal S)$ scalar product and applying \eqref{Est-K-k-Lp-Lp'} together with
Minkowski's integral inequality in $z_1$, we thus obtain
\begin{align*}
&\big\|\widetilde G_\tau\widetilde G_\tau^\star v(x_1,\cdot)\big\|_{L^{\frac{2d}{d-2}}(\R^{d-1})}\\
&\leq C\int_0^L\int_{2\tau|\ell|}^{\infty}
\big\|g_{\tau,0}(x_1,\lambda,\cdot)\big\|_{L^\infty(\Sigma)}
\big\|g_{\tau,0}(z_1,\lambda,\cdot)\big\|_{L^\infty(\Sigma)}\,\lambda^{1-\frac{2}{d}}\,d\lambda\,
\big\|v(z_1,\cdot)\big\|_{L^{\frac{2d}{d+2}}(\R^{d-1})}\,dz_1
\\
&\leq C\int_0^L\int_0^{\infty}\frac{\lambda^{1-\frac{2}{d}}}{1+\lambda}\,
e^{-\frac{1}{4}|\lambda+\tau\varphi'(x_1)|\,x_1-\frac{\beta}{2}\tau x_1^2}
e^{-\frac{1}{4}|\lambda+\tau\varphi'(z_1)|\,z_1-\frac{\beta}{2}\tau z_1^2}\,d\lambda\,
\big\|v(z_1,\cdot)\big\|_{L^{\frac{2d}{d+2}}(\R^{d-1})}\,dz_1 ,
\end{align*}
where we used $|\xi'|\geq\lambda$, by \eqref{psi-bounds}. Setting $S:=x_1+z_1$ and splitting the
$\lambda$-integral into $\big]0,\frac{\alpha_1}{2}\tau\big]$,
$\big]\frac{\alpha_1}{2}\tau,2\alpha_2\tau\big]$ and $\big]2\alpha_2\tau,+\infty\big[$, on which
$|\lambda+\tau\varphi'(x_1)|$ is respectively bounded below by $\frac{\alpha_1}{2}\tau$, possibly zero, and
bounded below by $\frac{\lambda}{2}$, we obtain
\begin{multline*}
\int_0^{\infty}\frac{\lambda^{1-\frac{2}{d}}}{1+\lambda}\,
e^{-\frac{1}{4}|\lambda+\tau\varphi'(x_1)|\,x_1-\frac{\beta}{2}\tau x_1^2}
e^{-\frac{1}{4}|\lambda+\tau\varphi'(z_1)|\,z_1-\frac{\beta}{2}\tau z_1^2}\,d\lambda
\\\leq C\Big(\tau^{1-\frac{2}{d}}e^{-c\tau S}+\tau^{-\frac{2}{d}}\min\{\tau,S^{-1}\}+S^{-1+\frac{2}{d}}\Big)
\leq\frac{C}{S^{1-\frac{2}{d}}}
\leq\frac{C}{|x_1-z_1|^{1-\frac{2}{d}}},
\end{multline*}
where the three terms come from the three ranges of $\lambda$: in the intermediate one, we discarded the exponential at $\min\{x_1,z_1\}$ and integrated the one at $\max\{x_1,z_1\}\geq\frac S2$,
which gives a factor $\min\{C\tau,C/S\}$, and in the last one we used $\int_0^\infty\lambda^{-\frac{2}{d}}e^{-\frac{\lambda S}8}d\lambda=CS^{-1+\frac{2}{d}}$; the last inequalities follow by distinguishing the
cases $\tau S\leq1$ and $\tau S\geq1$. The Hardy--Littlewood--Sobolev inequality in the $x_1$-variable then
concludes the proof of \eqref{Est-G-tau-high}.
\end{proof}
%
%

\subsubsection{Proof of Proposition \ref{Prop-Carleman-Strip-2}}
\label{Subsubsec-Proof-Prop-Strip}

Within the framework of Proposition \ref{Prop-Carleman-Strip-2}, we rely on the parametrix representation
given in Proposition \ref{Prop-parametrix-Strip-1} for the solution $w$ of \eqref{Elliptic-Strip-w},
namely,  
\begin{equation*}
w=K_{\tau,0}\big(f_2+f_{2_*'}\big)+K_{\tau,1}(F_1)+\sum_{j=2}^dK_{\tau,j}(F_j)
+G_\tau(g)+H_\tau\big((\partial_1w-F_1)(0,\cdot)\big).
\end{equation*}
Decomposing each operator as $K_{\tau,j}=K_{\tau,j}\Pi_{\mathrm{low}}+K_{\tau,j}\Pi_{\mathrm{high}}$, and similarly for $G_\tau$ and
$H_\tau$, and combining the low-frequency estimates of Proposition \ref{Prop-low-frequencies} with the
high-frequency ones of Propositions \ref{Prop-K-high} and \ref{Prop-GH-high}, we obtain, for $\theta \in\{0,-\frac{1}{2}\}$, on the one hand,
\begin{multline}
\label{Carleman-Strip-w-1-proof}
\tau^{\frac{3}{2}}\|w\|_{L^2(\mathcal S)}
+\tau^{\frac{1}{2}}\|\nabla w\|_{L^2(\mathcal S)}
\leq
C\left(
\|f_2\|_{L^2(\mathcal S)}
+\tau^{\frac{3}{4}-\frac{1}{2d}}\|f_{2_*'}\|_{L^{\frac{2d}{d+2}}(\mathcal S)}
+\tau\|F\|_{L^2(\mathcal S)}
\right.
\\
\left.
+\tau^{\frac{3}{4}}\|g\|_{H^{\frac{1}{2}}(\R^{d-1})}
+\tau^{\frac{5}{4}}\|g\|_{L^{2}(\R^{d-1})}
+\tau^{\frac{1}{4}-\theta}\|(\partial_1w-F_1)(0,\cdot)\|_{H^{\theta}( \R^{d-1})}
\right),
\end{multline}
and, on the other hand,
\begin{multline}
\label{Carleman-Strip-w-2-proof}
\tau^{\frac{3}{4}+\frac{1}{2d}}\|w\|_{L^{\frac{2d}{d-2}}(\mathcal S)}
\leq
C\left(
\|f_2\|_{L^2(\mathcal S)}
+\tau^{\frac{3}{4}+\frac{1}{2d}}\|f_{2_*'}\|_{L^{\frac{2d}{d+2}}(\mathcal S)}
+\tau\|F\|_{L^2(\mathcal S)}
\right.
\\
\left.
+\tau^{\frac{3}{4}+\frac{1}{2d}}\|g\|_{H^{\frac{1}{2}}(\R^{d-1})}
+\tau^{\frac{1}{4}+\frac{1}{2d} - \theta}\|(\partial_1w-F_1)(0,\cdot)\|_{H^{\theta}( \R^{d-1})}
\right).
\end{multline}
 This ends the the proof of Proposition \ref{Prop-Carleman-Strip-2}.

 \subsection{Proof of Theorem \ref{Thm-Local-Carl-x-0}}\label{Sec-Proof-Carl-x-0}
In this section, we are going to show how Theorem \ref{Thm-Local-Carl-x-0} can be deduced from
Proposition \ref{Prop-Carleman-Strip-2}. We shall use the following multiplication estimate (proven in Appendix \ref{Proof-Lemma-est-product-of-function-in-H-1/2-and--1/2}), in which the
essential point is that the constant does not depend on the scale $\rho$.
\begin{lemma}\label{Lemma-est-product-of-function-in-H-1/2-and--1/2}
Let $n\geq2$, $\rho>0$, $x_0'\in\R^{n}$ and $M\geq1$, and let $\widetilde m\in W^{1,\infty}(\R^{n})$
satisfy
\begin{equation*}
\operatorname{supp}\widetilde m\subset B(x_0',2\rho),
\qquad
\|\widetilde m\|_{L^\infty(\R^{n})}\leq M,
\qquad
\|\nabla\widetilde m\|_{L^\infty(\R^{n})}\leq\frac{M}{\rho} .
\end{equation*}
Then there exists $C>0$, depending only on $n$, such that, for all $h$,
\begin{equation}\label{Est-product-scale}
\|\widetilde m\,h\|_{H^{\frac{1}{2}}(\R^{n})}\leq C\,M\,\|h\|_{H^{\frac{1}{2}}(\R^{n})},
\qquad
\|\widetilde m\,h\|_{H^{-\frac{1}{2}}(\R^{n})}\leq C\,M\,\|h\|_{H^{-\frac{1}{2}}(\R^{n})} .
\end{equation}
\end{lemma}

\begin{proof}[Proof of Theorem \ref{Thm-Local-Carl-x-0}]
The proof is achieved by splitting the weight function at $x_0$, allowing us to reduce
\eqref{Laplacian-conjugated-w} to the model problem treated by the Carleman estimate of
Proposition \ref{Prop-Carleman-Strip-2}. Fix $x_0=(a,x_0')\in\overline{B^+_{r_0}}$ and $\theta\in\{0,-\frac{1}{2}\}$.

\medskip\noindent\textit{Step 1: linearisation of the weight function and reduction to the model problem.}
Expanding the square,
\begin{equation}\label{splitting-face}
\Phi(x)=\varphi(x_1)-\frac\varepsilon2|x'|^2
=\underbrace{\varphi(x_1)+\ell\cdot(x'-x_0')-\frac\varepsilon2|x_0'|^2}_{=:\ \phi_\ell(x)}
\;+\;q(x'),
\end{equation}
where $\ell=-\varepsilon x_0'$ and $q(x')=-\tfrac\varepsilon2|x'-x_0'|^2$. Note that $\phi_\ell$ differs
from the weight $\varphi(x_1)+\ell\cdot x'$ of Proposition \ref{Prop-Carleman-Strip-2} by an additive constant. So we have 
\begin{equation}\label{Conj-const}
e^{\tau\phi_\ell}\,\Delta\big(e^{-\tau\phi_\ell}\,\cdot\,\big)
=e^{\tau(\varphi+\ell\cdot x')}\,\Delta\big(e^{-\tau(\varphi+\ell\cdot x')}\,\cdot\,\big).
\end{equation}
On $B_{x_0}((\varepsilon\tau)^{-1/2})$, where $|x'-x_0'|\leq (\varepsilon\tau)^{-1/2}$, we have
\begin{equation}\label{Bounds-ell-and-q}
|\ell|\leq \varepsilon\,r_0\leq \ell_*,\qquad
\tau\,\|q\|_{L^\infty}\leq \frac{1}{2},\qquad
\tau\,\|\nabla'q\|_{L^\infty}\leq (\varepsilon\tau)^{1/2}.
\end{equation}
 Consequently, for $w$ solving
\eqref{Laplacian-conjugated-w}, we set
\begin{equation}\label{Def-check-w}
\check w(x)=e^{-\tau q(x')}\,w(x),
\end{equation}
so that, on $B_{x_0}((\varepsilon\tau)^{-1/2})$, $|\check w|$ and $|w|$ are pointwise comparable. Indeed, writing $v=e^{-\tau\Phi}w$, so that $\check w=e^{\tau\phi_\ell}v$ by
\eqref{splitting-face}, and commuting $e^{-\tau q}$ with the divergence (only the tangential part
contributes, as $\partial_1q=0$), $\check w$ satisfies
\begin{equation}\label{equation-check-w}
e^{\tau\phi_\ell}\Delta\big(e^{-\tau\phi_\ell}\check w\big)
=\widetilde f_2+\widetilde f_{2_*'}+\operatorname{div}\widetilde F
\qquad\text{in }\Omega\cap B_{x_0}\big((\varepsilon\tau)^{-1/2}\big),
\end{equation}
with source terms
\begin{equation}\label{Def-sources-terms}
\widetilde F=e^{-\tau q}\,F,\qquad
\widetilde f_{2_*'}=e^{-\tau q}\,f_{2_*'},\qquad
\widetilde f_2=e^{-\tau q}\,f_2+\tau\,\nabla'q\cdot e^{-\tau q}F' .
\end{equation}
By \eqref{Conj-const}, and since $\nabla\phi_\ell=(\varphi'(x_1),\ell)$, $\Delta\phi_\ell=\varphi''(x_1)$,
the operator in \eqref{equation-check-w} expands as
\begin{equation}\label{expanded-operator}
e^{\tau\phi_\ell}\Delta\big(e^{-\tau\phi_\ell}\,\cdot\,\big)
=\Delta-2\tau\varphi'(x_1)\,\partial_1-2\tau\,\ell\cdot\nabla'
+\tau^2\varphi'(x_1)^2+\tau^2|\ell|^2-\tau\varphi''(x_1),
\end{equation}
which is the operator \eqref{Conjugated-Operator-strip} of Proposition \ref{Prop-Carleman-Strip-2}. 

Recall $r_0\leq \tfrac{1}{2}$ and set $L=2r_0\leq 1$, so that the strip \eqref{Def-Strip} is
$\mathcal S=(0,L)\times\R^{d-1}$. Since $w\in\overline{C_c^\infty}^{\,\Omega}\big(B_{x_0}((\varepsilon\tau)^{-1/2})\cap B_0(r_0)\big)$, $w$ vanishes near $\partial\Omega\setminus\Gamma$ and near $\{x_1=L\}$, so that we can extend $\check w$ and all the source terms \eqref{Def-sources-terms} by $0$ outside
$\Omega\cap B_{x_0}((\varepsilon\tau)^{-1/2})$, and $\check w$ then satisfies
\begin{equation}\label{Elliptic-check-w-strip}
\left\{
\begin{array}{ll}
\Delta\check w-2\tau\varphi'(x_1)\,\partial_1\check w-2\tau\,\ell\cdot\nabla'\check w
+\big(\tau^2\varphi'(x_1)^2+\tau^2|\ell|^2-\tau\varphi''(x_1)\big)\check w
=\widetilde f_2+\widetilde f_{2_*'}+\operatorname{div}\widetilde F & \text{in }\mathcal S,\\[1mm]
\check w(0,x')=\widetilde g(x'), & x'\in\R^{d-1},\\[1mm]
\check w(L,x')=0, & x'\in\R^{d-1},
\end{array}
\right.
\end{equation}
where, since $\partial_1q=0$, the boundary terms are
\begin{equation}\label{boundary-terms}
\widetilde g(x')=e^{-\tau q(x')}\,w_{|\Gamma}(x'),\qquad
\big(\partial_1\check w-\widetilde F_1\big)(0,x')
=-\,e^{-\tau q(x')}\,\big((\nabla w-F)\cdot\nu\big)_{|\Gamma}(x'),
\end{equation}
using $\nu=-e_1$, both extended by $0$ outside $\Gamma_{x_0}=\Gamma\cap B_{x_0}((\varepsilon\tau)^{-1/2})$.
If $\Gamma_{x_0}=\emptyset$, both boundary data vanish.

\medskip\noindent\textit{Step 2: applying the Carleman estimate of Proposition \ref{Prop-Carleman-Strip-2}.}
The hypotheses of Proposition \ref{Prop-Carleman-Strip-2} hold on $\mathcal S$: $0<L\le1$; the conditions
\eqref{Weight-function-varphi} on $[0,L]$; and $|\ell|\leq \ell_*$ by \eqref{Bounds-ell-and-q}. Applying
\eqref{Carleman-Strip-w-1} and \eqref{Carleman-Strip-w-2} to $\check w$, we obtain, for all
$\tau\geq\tau_0$,
\begin{multline}\label{Carleman-Strip-check-w-L2}
\tau^{\frac{3}{2}}\|\check w\|_{L^2(\mathcal S)}+\tau^{\frac{1}{2}}\|\nabla\check w\|_{L^2(\mathcal S)}
+\tau^{\frac{3}{4}-\frac{1}{2d}}\|\check w\|_{L^{\frac{2d}{d-2}}(\mathcal S)}
\leq C\Big(\|\widetilde f_2\|_{L^2(\mathcal S)}
+\tau^{\frac{3}{4}-\frac{1}{2d}}\|\widetilde f_{2_*'}\|_{L^{\frac{2d}{d+2}}(\mathcal S)}
+\tau\|\widetilde F\|_{L^2(\mathcal S)}\\
+\tau^{\frac{3}{4}}\|\widetilde g\|_{H^{\frac{1}{2}}(\R^{d-1})}
+\tau^{\frac{5}{4}}\|\widetilde g\|_{L^{2}(\R^{d-1})}
+\tau^{\frac{1}{4}-\theta}\|\partial_1\check w-\widetilde F_1\|_{H^{\theta}(\{0\}\times\R^{d-1})}\Big),
\end{multline}
and
\begin{multline}\label{Carleman-Strip-check-w-L2_*}
\tau^{\frac{3}{4}+\frac{1}{2d}}\|\check w\|_{L^{\frac{2d}{d-2}}(\mathcal S)}
+\tau^{\frac{3}{2}}\|\check w\|_{L^2(\mathcal S)}+\tau^{\frac{1}{2}}\|\nabla\check w\|_{L^2(\mathcal S)}
\leq C\Big(\|\widetilde f_2\|_{L^2(\mathcal S)}
+\tau^{\frac{3}{4}+\frac{1}{2d}}\|\widetilde f_{2_*'}\|_{L^{\frac{2d}{d+2}}(\mathcal S)}
+\tau\|\widetilde F\|_{L^2(\mathcal S)}\\
+\tau^{\frac{3}{4}+\frac{1}{2d}}\|\widetilde g\|_{H^{\frac{1}{2}}(\R^{d-1})}
+\tau^{\frac{5}{4}}\|\widetilde g\|_{L^{2}(\R^{d-1})}
+\tau^{\frac{1}{4}+\frac{1}{2d}-\theta}\|\partial_1\check w-\widetilde F_1\|_{H^{\theta}(\{0\}\times\R^{d-1})}\Big),
\end{multline}
where the $L^{\frac{2d}{d-2}}$ term in the left-hand side of \eqref{Carleman-Strip-check-w-L2} comes from
\eqref{Carleman-Strip-w-2} multiplied by $\tau^{-\frac1d}$, and the $L^2$ terms in the left-hand side of
\eqref{Carleman-Strip-check-w-L2_*} from \eqref{Carleman-Strip-w-1}.

\medskip\noindent\textit{Step 3: pull-back.} We then simply remark that, from the expression of $\widetilde f_2, \widetilde f_{2_*'},$ and $ \widetilde F$ and \eqref{Bounds-ell-and-q},
 \begin{equation*}
     \|\widetilde f_2\|_{L^2(\mathcal S)}
+\tau^{\frac{3}{4}\mp\frac{1}{2d}}\|\widetilde f_{2_*'}\|_{L^{\frac{2d}{d+2}}(\mathcal S)}
+\tau\|\widetilde F\|_{L^2(\mathcal S)}
\leq C\Big(\|f_2\|_{L^2(\Omega)}
+\tau^{\frac{3}{4}\mp\frac{1}{2d}}\|f_{2_*'}\|_{L^{\frac{2d}{d+2}}(\Omega)}
+\tau\|F\|_{L^2(\Omega)}\Big).
 \end{equation*}
On the other hand, for the boundary terms, we observe that $w_{|\Gamma}$ and
$\big((\nabla w-F)\cdot\nu\big)_{|\Gamma}$ are supported in $\Gamma_{x_0}$, hence in the ball
$B(x_0',(\varepsilon \tau)^{-1/2})$ of $\R^{d-1}$. Let therefore $\chi\in C^\infty_c\big(B(x_0',2 (\varepsilon \tau)^{-1/2})\big)$ be such
that $\chi=1$ on $B(x_0',(\varepsilon \tau)^{-1/2})$, $0\leq\chi\leq1$ and $\|\nabla_{x'}\chi\|_{L^\infty}\leq C(\varepsilon \tau)^{1/2}$,
and set $\widetilde m:=\chi\,e^{-\tau q}$, so that, by \eqref{boundary-terms},
\begin{equation*}
\widetilde g=\widetilde m\,w_{|\Gamma},
\qquad
\big(\partial_1\check w-\widetilde F_1\big)(0,\cdot)
=-\,\widetilde m\,\big((\nabla w-F)\cdot\nu\big)_{|\Gamma} .
\end{equation*}
By \eqref{Bounds-ell-and-q}, the function $\widetilde m$ satisfies the hypotheses of
Lemma \ref{Lemma-est-product-of-function-in-H-1/2-and--1/2} in $\R^{d-1}$, with $\rho=(\varepsilon \tau)^{-1/2}$ and $M=C$, indeed
$\operatorname{supp}\widetilde m\subset B(x_0',2(\varepsilon \tau)^{-1/2})$, $\|\widetilde m\|_{L^\infty}\leq e^{2}$ and
\begin{equation*}
\|\nabla_{x'}\widetilde m\|_{L^\infty}
\leq\|\nabla_{x'}\chi\|_{L^\infty}\|e^{-\tau q}\|_{L^\infty}
+\tau\|\nabla_{x'}q\|_{L^\infty}\|e^{-\tau q}\|_{L^\infty}\leq C (\varepsilon \tau)^{1/2} .
\end{equation*}
Lemma \ref{Lemma-est-product-of-function-in-H-1/2-and--1/2}, together with $\|\widetilde m\|_{L^\infty}\leq e^2$ for the $L^2$ norms, therefore yields
\begin{align}
&\|\widetilde g\|_{H^{\frac{1}{2}}(\R^{d-1})}
\leq C\,\|w_{|\Gamma}\|_{H^{\frac{1}{2}}(\Gamma_{x_0})},
\quad
\|\widetilde g\|_{L^{2}(\R^{d-1})}\leq C\,\|w_{|\Gamma}\|_{L^{2}(\Gamma_{x_0})},
\label{bdata-bound-1}\\
&\|\partial_1\check w-\widetilde F_1\|_{H^{\theta}(\{0\}\times\R^{d-1})}
\leq C\,\|(\nabla w-F)\cdot\nu_{|\Gamma}\|_{H^{\theta}(\Gamma_{x_0})},
\qquad \theta\in\Big\{0,-\frac{1}{2}\Big\},
\label{bdata-bound-2}
\end{align}
with constants independent of $\tau$ and $\varepsilon$. Consequently, for all $\tau\geq\tau_0$, we get
\begin{multline}\label{Carleman-Strip-check-w-L2-1}
\tau^{\frac{3}{2}}\|\check w\|_{L^2(\mathcal S)}+\tau^{\frac{1}{2}}\|\nabla\check w\|_{L^2(\mathcal S)}
+\tau^{\frac{3}{4}-\frac{1}{2d}}\|\check w\|_{L^{\frac{2d}{d-2}}(\mathcal S)}
\leq C\left(\|f_2\|_{L^2(\Omega)}
+\tau^{\frac{3}{4}-\frac{1}{2d}}\|f_{2_*'}\|_{L^{\frac{2d}{d+2}}(\Omega)}
+\tau\|F\|_{L^2(\Omega)}\right.\\
\left.
+\tau^{\frac{3}{4}}\|w\|_{H^{\frac{1}{2}}(\Gamma_{x_0})}
+\tau^{\frac{5}{4}}\|w\|_{L^{2}(\Gamma_{x_0})}
+\tau^{\frac{1}{4}-\theta}\|(\nabla w-F)\cdot\nu_{|\Gamma}\|_{H^{\theta}(\Gamma_{x_0})}\right),
\end{multline}
and
\begin{multline}\label{Carleman-Strip-check-w-L2_*-1}
\tau^{\frac{3}{4}+\frac{1}{2d}}\|\check w\|_{L^{\frac{2d}{d-2}}(\mathcal S)}
+\tau^{\frac{3}{2}}\|\check w\|_{L^2(\mathcal S)}+\tau^{\frac{1}{2}}\|\nabla\check w\|_{L^2(\mathcal S)}
\leq C\left(\|f_2\|_{L^2(\Omega)}
+\tau^{\frac{3}{4}+\frac{1}{2d}}\|f_{2_*'}\|_{L^{\frac{2d}{d+2}}(\Omega)}
+\tau\|F\|_{L^2(\Omega)}\right.\\
\left.
+\tau^{\frac{3}{4}+\frac{1}{2d}}\|w\|_{H^{\frac{1}{2}}(\Gamma_{x_0})}
+\tau^{\frac{5}{4}}\|w\|_{L^{2}(\Gamma_{x_0})}
+\tau^{\frac{1}{4}+\frac{1}{2d}-\theta}\|(\nabla w-F)\cdot\nu_{|\Gamma}\|_{H^{\theta}(\Gamma_{x_0})}\right).
\end{multline}

It remains to return to $w=e^{\tau q}\check w$. Since $|\check w|$ and $|w|$ are pointwise comparable, one
has $\|\check w\|_{L^p(\mathcal S)}\geq C^{-1}\|w\|_{L^p(\Omega)}$ for $p\in\{2,\frac{2d}{d-2}\}$.
Moreover $\nabla w=e^{\tau q}\big(\nabla\check w+\tau\nabla q\,\check w\big)$, so that, by
\eqref{Bounds-ell-and-q}, we get
\begin{equation*}
\tau^{\frac{1}{2}}\|\nabla w\|_{L^2(\Omega)}
\leq C\Big(\tau^{\frac{1}{2}}\|\nabla\check w\|_{L^2(\mathcal S)}
+\tau^{\frac{1}{2}}(\varepsilon\tau)^{\frac{1}{2}}\|\check w\|_{L^2(\mathcal S)}\Big)
\leq C\Big(\tau^{\frac{1}{2}}\|\nabla\check w\|_{L^2(\mathcal S)}
+\tau^{\frac{3}{2}}\|\check w\|_{L^2(\mathcal S)}\Big).
\end{equation*}
Hence the left-hand sides of
\eqref{Carleman-Strip-check-w-L2-1}--\eqref{Carleman-Strip-check-w-L2_*-1} dominate, up to constants, those
of \eqref{Carleman-Local-w-x0-L2} and \eqref{Carleman-Local-w-x0-L2_*}. This ends the proof
of Theorem \ref{Thm-Local-Carl-x-0}.
\end{proof}
\section{Variable Coefficient Carleman Estimate : Proof of Theorem \ref{Thm-Carleman-General}} \label{Sec-Carleman-variable}
This section is devoted to the proof of Theorem \ref{Thm-Carleman-General}. The proof relies on two main ingredients:
\begin{itemize} 
    \item Step 1 : extending the results proved in the previous section for the flat case to the variable-coefficient operator $\partial_i\!\big(a^{ij}(\varepsilon x)\partial_j\big)$.
    \item Step 2 : a gluing procedure based on a partition of unity to combine the resulting local Carleman estimates.
\end{itemize} 
\subsection{Step 1: Local Carleman estimate with variable coefficients} \label{Sec-perturbation-argument}
For $\tau\geq1$ and $\varepsilon\in(0,1]$, define the conjugated operator
\begin{equation*}
\mathcal L_{\tau,\varepsilon}^A w
:=e^{\tau\Phi}\,\partial_i\!\Big(a^{ij}(\varepsilon x)\,\partial_j\big(e^{-\tau\Phi}w\big)\Big),
\end{equation*}
where $\Phi$ is the weight \eqref{Weight-function-Phi}. A direct computation gives
\begin{equation}\label{Expression-LtauA}
\mathcal L_{\tau,\varepsilon}^A w
=\partial_i\big(a^{ij}(\varepsilon x)\partial_j w\big)
-\tau\,\partial_i\!\big(a^{ij}(\varepsilon x)(\partial_j\Phi)\big)\,w
-2\tau\,a^{ij}(\varepsilon x)(\partial_i\Phi)\,\partial_j w
+\tau^2 a^{ij}(\varepsilon x)(\partial_i\Phi)(\partial_j\Phi)\,w .
\end{equation}
\begin{proposition}\label{Prop-Local-Carl-x-0-variable}
Let $d\geq3$, $r_0\in(0,\frac{1}{2}]$, and let $\Omega=B^+_{2r_0}$, $\Gamma$ and $\nu$ be as in
Theorem \ref{Thm-Carleman-General}. Let $A$ satisfy \eqref{Ellipticity-A}--\eqref{Lipschitz-condition-A},
let $\varphi$ satisfy \eqref{Weight-function-varphi} and let $\Phi$ be given by \eqref{Weight-function-Phi}. There exist $C>0$, $\varepsilon_1\in(0,1]$ and $\tau_0\geq1$, depending only on
$d,\Lambda_0,\Lambda_1,\alpha_1,\alpha_2,\beta,$ and $r_0$, such that, for all
$0<\varepsilon\leq\varepsilon_1$, all $\tau\geq\tau_0$, all
$x_0\in\overline{B^+_{r_0}}$, all $w\in\overline{C_c^\infty}^{\,\Omega}\big(B_{x_0}((\varepsilon\tau)^{-1/2})\cap B_0(r_0)\big)$, and every decomposition
\begin{equation}\label{equation-w-LtauA}
\mathcal L_{\tau,\varepsilon}^A w=f_2+f_{2_*'}+\div F
\qquad\text{in }\Omega\cap B_{x_0}((\varepsilon\tau)^{-1/2}),
\end{equation}
with $f_2\in L^2(\Omega)$, $f_{2_*'}\in L^{\frac{2d}{d+2}}(\Omega)$, $F\in L^2(\Omega;\C^d)$ vanishing outside $B_{x_0}((\varepsilon\tau)^{-1/2})\cap B_0(\frac{3}{2}r_0)$, the following
estimates hold
\begin{multline}\label{Carleman-Local-A-L2}
\tau^{\frac{3}{2}}\|w\|_{L^2(\Omega)}+\tau^{\frac{1}{2}}\|\nabla w\|_{L^2(\Omega)}
+\tau^{\frac{3}{4}-\frac{1}{2d}}\|w\|_{L^{\frac{2d}{d-2}}(\Omega)}
\leq C\Big(\|f_2\|_{L^2(\Omega)}
+\tau^{\frac{3}{4}-\frac{1}{2d}}\|f_{2_*'}\|_{L^{\frac{2d}{d+2}}(\Omega)}
+\tau\|F\|_{L^2(\Omega)}\\
+\tau^{\frac{3}{4}}  \|w\|_{H^{\frac{1}{2}}(\Gamma_{x_0})}
+\tau^{\frac{5}{4}}\|w\|_{L^{2}(\Gamma_{x_0})}
+\tau^{\frac{3}{4} - \theta}\big\|\big(A(\varepsilon\cdot)\nabla w-F\big)\cdot\nu\big\|_{H^{\theta}(\Gamma_{x_0})} \Big),
\end{multline}
and
\begin{multline}\label{Carleman-Local-A-L2star}
\tau^{\frac{3}{4}+\frac{1}{2d}}\|w\|_{L^{\frac{2d}{d-2}}(\Omega)}
+\tau^{\frac{3}{2}}\|w\|_{L^2(\Omega)}+\tau^{\frac{1}{2}}\|\nabla w\|_{L^2(\Omega)}
\leq C\Big(\|f_2\|_{L^2(\Omega)}
+\tau^{\frac{3}{4}+\frac{1}{2d}}\|f_{2_*'}\|_{L^{\frac{2d}{d+2}}(\Omega)}
+\tau\|F\|_{L^2(\Omega)}\\
+\tau^{\frac{5}{4}}\|w\|_{L^{2}(\Gamma_{x_0})}
+\tau^{\frac{3}{4}+\frac{1}{2d}} \left(\|w\|_{H^{\frac{1}{2}}(\Gamma_{x_0})}
+\tau^{-\frac{1}{2}-\theta}\big\|\big(A(\varepsilon\cdot)\nabla w-F\big)\cdot\nu\big\|_{H^{\theta}(\Gamma_{x_0})}\right)\Big),
\end{multline}
where  $\Gamma_{x_0}=\Gamma\cap B_{x_0}(\varepsilon\tau)^{-1/2})$ and $\theta \in \{0,-1/2\}$. 
\end{proposition}
Let us explain the choice of the scale $(\varepsilon\tau)^{-1/2}$. Freezing the coefficients on a ball of radius $r$ produces, by \eqref{Expression-LtauA}, errors of size $\Lambda_1 r\,(\tau^2\|w\|_{L^2}+\tau\|\nabla w\|_{L^2})$, which are absorbed by the left-hand side of \eqref{Carleman-Local-A-L2} only if $r\ll\tau^{-1/2}$. On the other hand, the cut-off functions at scale $r$ used in the gluing argument of Subsection \ref{Gluing-argument} produce errors of size $r^{-1}(\|\nabla w\|_{L^2}+\tau\|w\|_{L^2})$, which are absorbed only if $r\gg\tau^{-1/2}$. Both requirements are critical at $r\sim\tau^{-1/2}$. For this reason, we replace $A$ by $A(\varepsilon\cdot)$, whose Lipschitz constant is $\Lambda_1\varepsilon$, and we localize at the scale $(\varepsilon\tau)^{-1/2}$. The freezing errors are then of size $\Lambda_1\varepsilon^{\frac{1}{2}}(\tau^{\frac{3}{2}}\|w\|_{L^2}+\tau^{\frac{1}{2}}\|\nabla w\|_{L^2})$, the cut-off errors of size $\varepsilon^{\frac{1}{2}}(\tau^{\frac{1}{2}}\|\nabla w\|_{L^2}+\tau^{\frac{3}{2}}\|w\|_{L^2})$, and both are absorbed for $\varepsilon$ small.

\begin{proof} \noindent\textit{Step 1: normalization.}
We first want to normalise the coefficients at $x_0$, and, since Theorem \ref{Thm-Local-Carl-x-0} is stated
on a half-ball with boundary terms on $\{x_1=0\}$, we must use a change of variables preserving the
hyperplanes $\{x_1=\mathrm{const}\}$. To do that, we use a change of coordinates somewhat similar to the one
introduced in \cite{koch2001carleman}. Multiplying $A$ by $\big(a^{11}(\varepsilon x_0)\big)^{-1}$, we may assume that $a^{11}(\varepsilon x_0)=1$. Moreover, a translation in the tangential directions replaces $-\frac\varepsilon2|x'|^2$ by $-\frac\varepsilon2|x'+x_0'|^2$, that is, adds to $\Phi$ a linear tangential term $\ell_0\cdot x'$, with $\ell_0=-\varepsilon x_0'$, and a constant, which are handled as in \eqref{splitting-face}; we can therefore assume that $x_0'=0$, the weight being $\Phi(x)=\varphi(x_1)+\ell_0\cdot x'-\frac\varepsilon2|x'|^2$ with $|\ell_0|\leq\varepsilon r_0$. For the normal component of $x_0=(a,x_0')$, either $\Gamma_{x_0}=\emptyset$ and the boundary terms are zero, or
$\Gamma_{x_0}\neq\emptyset$, in which case $0\leq a<(\varepsilon\tau)^{-1/2}$ and, up to replacing $(\varepsilon\tau)^{-1/2}$ by
$2(\varepsilon\tau)^{-1/2}$, we can also assume $a=0$. We may thus assume that $x_0=0$ (if $\Gamma_{x_0}=\emptyset$, the argument below applies with $y_1$ replaced by $y_1-a$ in the tangential components of $S$ and in $\varphi_2$), and we represent $A$ in the form
\begin{equation*}
A(0)=\begin{pmatrix}1&B^{T}\\ B&A'\end{pmatrix},
\qquad B^{j}=a^{1j}(0),\qquad A'=\big(a^{ij}(0)\big)_{i,j\geq2},
\end{equation*}
and we look for $S$ such that
\begin{equation}\label{Normalisation-S}
S^{-1}A(0)\,(S^{T})^{-1}=\operatorname{Id},
\qquad\text{that is}\qquad
A(0)=S\,S^{T}.
\end{equation}
Since $A(0)$ is uniformly positive definite, it follows that the $(d-1)\times(d-1)$ matrix
$A'-B\,B^{T}$ is uniformly positive definite. Now set
\begin{equation}\label{Def-S}
S=\begin{pmatrix}1&0\\ B&\big(A'-B\,B^{T}\big)^{\frac{1}{2}}\end{pmatrix}.
\end{equation}
By \eqref{Ellipticity-A}, the norms of $S$ and $S^{-1}$, therefore
the ratio of the radii of the balls exchanged by $S$, are controlled by $\Lambda_0$ and $\Lambda_1$,
uniformly in $\varepsilon$ and in $x_0$. From now on we work in the $y$-variables, we set $\check w(y)=w(S(y))$ and
$\check A(y)=S^{-1}A\big(S(y)\big)(S^{ T})^{-1}$, so that the coefficients become $\check A(\varepsilon\cdot)$, $\check A(0)=\operatorname{Id}$ and
$\check A$ is Lipschitz with constant $C\Lambda_1$, and to simplify notation, we continue writing $w$, $A$  for $\check w$, $\check A$. The equation \eqref{equation-w-LtauA} keeps the same form, with $F$ replaced by $S^{-1}F\circ S$, and, since the first rows of $S$ and $S^{-1}$ are $(1,0)$, $S$ preserves $\Gamma$, $\nu$ and the boundary terms $w_{|\Gamma}$ and $\big(A(\varepsilon\cdot)\nabla w-F\big)\cdot\nu$, all the norms changing only by constants depending on $\Lambda_0$ and $\Lambda_1$.

In the new coordinates, $x=S(y)=\big(y_1,\ Ry'+y_1B\big)$ with
$R=\big(A'-B\,B^{T}\big)^{\frac{1}{2}}$, we have
\begin{equation}\label{Transported-weight}
\Phi(x)=\underbrace{\varphi_0(y_1)+\big(R^T\ell_0\big)\cdot y'-\frac{\varepsilon}{2}\big|Ry'\big|^2}_{=:\ \Phi_0(y)}+\varphi_2(y),
\qquad\text{with}\qquad
\varphi_2(y)=-\varepsilon\,y_1\,\big(B\cdot Ry'\big)-\frac{\varepsilon}{2}|B|^2y_1^2 ,
\end{equation}
and $\varphi_0(y_1)=\varphi(y_1)+(\ell_0\cdot B)\,y_1$. The weight $\Phi_0$ is not exactly of the form \eqref{Weight-function-Phi}, but Theorem \ref{Thm-Local-Carl-x-0} applies to it. Indeed, since $|\ell_0\cdot B|\leq C\varepsilon r_0$, $\varphi_0$ satisfies \eqref{Weight-function-varphi} with $\frac{\alpha_1}2$, $2\alpha_2$ and $\beta$ for $\varepsilon$ small, and the proof of Theorem \ref{Thm-Local-Carl-x-0} only uses the splitting \eqref{splitting-face} of the tangential part of the weight at the centre of the ball, here $(R^T\ell_0)\cdot y'-\frac\varepsilon2|Ry'|^2$, with $|R^T\ell_0|\leq C\varepsilon r_0$, and the bounds \eqref{Bounds-ell-and-q}, which hold for $-\frac\varepsilon2|Ry'|^2$ with constants depending only on $|R|$. The function $\varphi_2$
satisfies
\begin{equation}\label{Bounds-phi-2}
\tau\|\varphi_2\|_{L^\infty(B_0((\varepsilon\tau)^{-1/2}))}\leq C,
\qquad
\tau\|\nabla\varphi_2\|_{L^\infty(B_0((\varepsilon\tau)^{-1/2}))}\leq C\,(\varepsilon\tau)^{\frac{1}{2}},
\qquad
|\nabla^2\varphi_2|\leq C\,\varepsilon .
\end{equation}
Thus $\varphi_2$ satisfies the same bounds as the quadratic remainder $q$ of \eqref{Bounds-ell-and-q} in
the proof of Theorem \ref{Thm-Local-Carl-x-0}. The only difference here is that $\varphi_2$ depends also on the normal variable, and $\partial_1\varphi_2(0,y')=-\varepsilon\,B\cdot Ry'$ does not vanish, so that a conjugation by $e^{-\tau\varphi_2}$ as in \eqref{Def-check-w} would modify the conormal derivative on $\Gamma_0$. We therefore treat $\varphi_2$ directly as a perturbation of the weight and, for the sake of completeness, we write the argument in full.

\medskip\noindent\textit{Step 2: perturbation of the flat conjugated operator.} Since $\Phi=\Phi_0+\varphi_2$, we write
\begin{equation}\label{Decomposition-L}
\mathcal L^A_{\tau,\varepsilon}w
=\mathcal L_\tau w+\mathcal R^{A}_{\tau,\varepsilon}w+\mathcal R^{\Phi}_{\tau,\varepsilon}w ,
\end{equation}
where $\mathcal L_\tau=e^{\tau\Phi_0}\Delta\big(e^{-\tau\Phi_0}\,\cdot\,\big)$ is the flat conjugated
operator of \eqref{Laplacian-conjugated-w} associated with $\Phi_0$, where $\mathcal R^{A}_{\tau,\varepsilon}=\mathcal L^A_{\tau,\varepsilon}-e^{\tau\Phi}\Delta\big(e^{-\tau\Phi}\,\cdot\,\big)$, that is, by \eqref{Expression-LtauA},
\begin{align*}
\mathcal R^{A}_{\tau,\varepsilon}w
&=\partial_i\!\big((a^{ij}(\varepsilon\cdot)-\delta^{ij})\partial_jw\big)
-\tau\,\partial_i\!\big((a^{ij}(\varepsilon\cdot)-\delta^{ij})(\partial_j\Phi)\big)\,w\\
&\qquad
-2\tau\,(a^{ij}(\varepsilon\cdot)-\delta^{ij})(\partial_i\Phi)\,\partial_jw
+\tau^2\,(a^{ij}(\varepsilon\cdot)-\delta^{ij})(\partial_i\Phi)(\partial_j\Phi)\,w ,
\end{align*}
and where $\mathcal R^{\Phi}_{\tau,\varepsilon}=e^{\tau\Phi}\Delta\big(e^{-\tau\Phi}\,\cdot\,\big)-\mathcal L_\tau$, that is,
\begin{equation}\label{Def-R-Phi}
\mathcal R^{\Phi}_{\tau,\varepsilon}w
=-2\tau\,\nabla\varphi_2\cdot\nabla w-\tau\,(\Delta\varphi_2)\,w
+\tau^2|\nabla\varphi_2|^2\,w+2\tau^2\,\nabla\Phi_0\cdot\nabla\varphi_2\,w.
\end{equation}
Since $A$ is Lipschitz, we have
\begin{equation}\label{Smallness-local-a}
\big\|A(\varepsilon\cdot)-\operatorname{Id}\big\|_{L^\infty(B_0((\varepsilon \tau)^{-1/2}))}
\leq C \varepsilon^{1/2}\tau^{-1/2},
\qquad
\big\|\nabla\big(A(\varepsilon\cdot)\big)\big\|_{L^\infty(B_0((\varepsilon \tau)^{-1/2}))}\leq C\Lambda_1\,\varepsilon .
\end{equation}

\medskip\noindent\textit{Step 3: interior perturbation bound.}
We put the first term of $\mathcal R^{A}_{\tau,\varepsilon}$ in divergence form, $\div(F_\flat)$ with
$F_{\flat,i}=\big(a^{ij}(\varepsilon\cdot)-\delta^{ij}\big)\partial_jw$, while the remaining terms of
$\mathcal R^{A}_{\tau,\varepsilon}$ and all the terms of $\mathcal R^{\Phi}_{\tau,\varepsilon}$ are
estimated directly in the $L^2$-norm. We therefore set
\begin{multline}\label{Def-f-2-flat}
f_{2,\flat}
=\tau^2\big(a^{ij}(\varepsilon\cdot)-\delta^{ij}\big)(\partial_i\Phi)(\partial_j\Phi)\,w
-2\tau\big(a^{ij}(\varepsilon\cdot)-\delta^{ij}\big)(\partial_i\Phi)\,\partial_jw \\
-\tau\,\partial_i \big((a^{ij}(\varepsilon\cdot)-\delta^{ij})(\partial_j\Phi)\big)\,w
+\mathcal R^{\Phi}_{\tau,\varepsilon}w .
\end{multline}
By \eqref{Smallness-local-a}, \eqref{Bounds-phi-2} and $|\nabla\Phi_0|+|\nabla^2\Phi_0|\leq C$, we have
\begin{equation}\label{estmation-F-flat-f-2-flat}
\tau\big\|F_\flat\big\|_{L^2(\Omega)}+\big\|f_{2,\flat}\big\|_{L^2(\Omega)}
\leq C\,\varepsilon^{\frac{1}{2}}\Big(\tau^{\frac{3}{2}}\|w\|_{L^2(\Omega)}
+\tau^{\frac{1}{2}}\|\nabla w\|_{L^2(\Omega)}\Big).
\end{equation}

\medskip\noindent\textit{Step 4: application of the flat estimate and boundary terms.}
By \eqref{Decomposition-L}, the function $w$ satisfies
\begin{equation*}
\mathcal L_\tau w=\big(f_2- f_{2,\flat}\big)+f_{2_*'}+\div\big(F-F_\flat\big)
\qquad\text{in }\Omega\cap B_{0}((\varepsilon\tau)^{-1/2}).
\end{equation*}
We therefore apply Theorem \ref{Thm-Local-Carl-x-0}, with the weight $\Phi_0$, to $w$ with
$\big(f_2- f_{2,\flat},\,f_{2_*'},\,F-F_\flat\big)$. Since $F_\flat=\big(A(\varepsilon\cdot)-\operatorname{Id}\big)\nabla w$, its boundary terms are $w_{|\Gamma_0}$ and
\begin{equation*}
\big(\nabla w-F+F_\flat\big)\cdot\nu
=\big(A(\varepsilon\cdot)\nabla w-F\big)\cdot\nu
\quad\text{in }H^{-\frac{1}{2}}(\Gamma_{0}).
\end{equation*}
Hence, using \eqref{estmation-F-flat-f-2-flat}, we obtain, for all
$\tau\geq\tau_0$,
\begin{multline*}
\tau^{\frac{3}{2}}\|w\|_{L^2(\Omega)}
+\tau^{\frac{1}{2}}\|\nabla w\|_{L^2(\Omega)}
+\tau^{\frac{3}{4}-\frac{1}{2d}}\|w\|_{L^{\frac{2d}{d-2}}(\Omega)}
\leq C\Big(\|f_2\|_{L^2(\Omega)}
+\tau^{\frac{3}{4}-\frac{1}{2d}}\|f_{2_*'}\|_{L^{\frac{2d}{d+2}}(\Omega)}
+\tau\|F\|_{L^2(\Omega)}\\
+\tau^{\frac{3}{4}}\|w\|_{H^{\frac{1}{2}}(\Gamma_{0})}
+\tau^{\frac{5}{4}}\|w\|_{L^{2}(\Gamma_{0})}
+\tau^{\frac{1}{4}-\theta}\big\|\big(A(\varepsilon\cdot)\nabla w-F\big)\cdot\nu\big\|_{H^{\theta}(\Gamma_{0})}\Big)
+C\,\varepsilon^{\frac{1}{2}}\Big(\tau^{\frac{3}{2}}\|w\|_{L^2(\Omega)}
+\tau^{\frac{1}{2}}\|\nabla w\|_{L^2(\Omega)}\Big),
\end{multline*}
choosing $\varepsilon\leq\varepsilon_1$, where $\varepsilon_1\leq\min\big\{(2C)^{-2},\varepsilon_0\big\}$ is also small enough for Step~1, and where $\varepsilon_0$ and
$\tau_0$ are the constants of Theorem \ref{Thm-Local-Carl-x-0}, the last two terms are absorbed into the
left-hand side, which yields \eqref{Carleman-Local-A-L2} in the $y$-variables for
all $\tau\geq\tau_0$ and all $\varepsilon\leq\varepsilon_1$. Undoing the change of variables
$S$, whose Jacobian and inverse are bounded by constants depending only on $\Lambda_0$ and $\Lambda_1$, we
obtain \eqref{Carleman-Local-A-L2} in the original variables. The estimate
\eqref{Carleman-Local-A-L2star} is obtained in the same way from \eqref{Carleman-Local-w-x0-L2_*}, which
completes the proof of Proposition \ref{Prop-Local-Carl-x-0-variable}.
\end{proof}

\subsection{Step 2: Gluing by partition of unity}\label{Gluing-argument}
We now derive \eqref{Carleman-General-1} from the local estimate \eqref{Carleman-Local-A-L2}; the
derivation of \eqref{Carleman-General-2} from \eqref{Carleman-Local-A-L2star} is identical, replacing the
exponent $\frac{3}{4}-\frac{1}{2d}$ by $\frac{3}{4}+\frac{1}{2d}$ throughout.

Set $\rho=\min\big\{(\varepsilon\tau)^{-1/2},\frac{r_0}{2}\big\}$, and let $(\eta_k)_{0\leq k\leq k_1}$ be a smooth partition of unity subordinated to a covering of
$\overline{\Omega\cap B_0(r_0)}$ by balls $B_k=B_{x_k}(\rho)$, with $x_k\in\overline{B^+_{r_0}}$,
\begin{equation}\label{Partition-unity}
\operatorname{supp} \eta_k\subset B_k,\qquad 0\leq\eta_k\leq1,\qquad\sum_k\eta_k=1\ \text{ on }
\overline{\Omega\cap B_0(r_0)},\qquad|\partial^\alpha\eta_k|\leq c_\alpha\rho^{-|\alpha|},
\end{equation}
the covering being chosen so that every point of $\R^d$ belongs to at most $N$ balls $B_k$, where $N$ depends only on $d$. In particular, $\sum_k\eta_k^2\leq N$ on $\R^d$ and $\sum_k\eta_k^2\leq\sum_k\eta_k=1$ on $\overline{\Omega\cap B_0(r_0)}$. Note also that $B_k\subset B_{x_k}((\varepsilon\tau)^{-1/2})\cap B_0(\frac{3}{2}r_0)$ and that, for $\tau\geq16r_0^{-2}$,
\begin{equation}\label{Bound-rho}
\rho^{-1}\leq(\varepsilon\tau)^{\frac{1}{2}}+2r_0^{-1}\leq2\tau^{\frac{1}{2}},
\qquad
\rho^{-2}\leq2\varepsilon\tau+8r_0^{-2}\leq\tau^2 .
\end{equation}

We introduce
\begin{equation*}
w=e^{\tau\Phi}u,\qquad
\widetilde f_2:=e^{\tau\Phi}\big(f_2-\tau\nabla\Phi\cdot F\big),\qquad
\widetilde f_{2_*'}:=e^{\tau\Phi}f_{2_*'},\qquad
\widetilde F:=e^{\tau\Phi}F,
\end{equation*}
so that, by \eqref{Elliptic-Eq-general}, $\mathcal L_{\tau,\varepsilon}^Aw=\widetilde f_2+\widetilde f_{2_*'}+\div\widetilde F$ in $\Omega$, and we set $w_k=\eta_kw$. Since $u\in\overline{C_c^\infty}^{\,\Omega}(B_0(r_0))$, we have $w_k\in\overline{C_c^\infty}^{\,\Omega}\big(B_{x_k}((\varepsilon\tau)^{-1/2})\cap B_0(r_0)\big)$. A direct computation gives
\begin{equation*}
\mathcal L_{\tau,\varepsilon}^A(\eta_kw)
=\eta_k\,\mathcal L_{\tau,\varepsilon}^Aw
+\partial_i\big(a^{ij}(\varepsilon x)(\partial_j\eta_k)\,w\big)
+(\partial_i\eta_k)\,a^{ij}(\varepsilon x)\,\partial_jw
-2\tau\,a^{ij}(\varepsilon x)(\partial_i\Phi)(\partial_j\eta_k)\,w ,
\end{equation*}
so that $w_k$ solves
\begin{equation*}
\mathcal L_{\tau,\varepsilon}^A(\eta_kw)=f_{2,k}+f_{2_*',k}+\div(F_k)
\qquad\text{in }\Omega,
\end{equation*}
where
\begin{align*}
&f_{2,k}=\eta_k\widetilde f_2-\nabla\eta_k\cdot\widetilde F
+(\partial_i\eta_k)\,a^{ij}(\varepsilon x)\,\partial_jw
-2\tau\,a^{ij}(\varepsilon x)(\partial_i\Phi)(\partial_j\eta_k)\,w,\\
&f_{2_*',k}=\eta_k\widetilde f_{2_*'},
\qquad
F_{i,k}=\eta_k\widetilde F_i+a^{ij}(\varepsilon x)(\partial_j\eta_k)\,w ,
\end{align*}
all vanishing outside $B_k\subset B_{x_k}((\varepsilon\tau)^{-1/2})\cap B_0(\frac{3}{2}r_0)$.
Applying Proposition \ref{Prop-Local-Carl-x-0-variable} to $w_k$, with $x_0=x_k$, and squaring, we obtain, for every $k$,
\begin{multline}\label{Carleman-Local-A-L2-wk}
\tau^{3}\|w_k\|_{L^2(\Omega)}^2+\tau^{\frac{3}{2}-\frac{1}{d}}\|w_k\|^2_{L^{\frac{2d}{d-2}}(\Omega)}
+\tau\|\nabla w_k\|_{L^2(\Omega)}^2 \leq C\Big(\|\eta_k\widetilde f_2\|_{L^2(\Omega)}^2
+\tau^{\frac{3}{2}-\frac{1}{d}}\|\eta_k\widetilde f_{2_*'}\|^2_{L^{\frac{2d}{d+2}}(\Omega)} \\
+\tau^2\|\eta_k\widetilde F\|_{L^2(\Omega)}^2
+\tau^{\frac{3}{2}}\|\eta_kw\|^2_{H^{\frac{1}{2}}(\Gamma)}
+\tau^{\frac{5}{2}}\|\eta_kw\|^2_{L^{2}(\Gamma)}
+\tau^{\frac{1}{2}-2 \theta}\big\|\eta_k\big(A(\varepsilon\cdot)\nabla w-\widetilde F\big)\cdot\nu\big\|^2
_{H^{\theta}(\Gamma)}\Big)\\
+C\Big(\|\nabla\eta_k\cdot\widetilde F\|_{L^2(\Omega)}^2
+\|(\partial_i\eta_k)a^{ij}(\varepsilon\cdot)(\partial_jw)\|_{L^2(\Omega)}^2\\
+\tau^2\|a^{ij}(\varepsilon\cdot)(\partial_i\Phi)(\partial_j\eta_k)w\|_{L^2(\Omega)}^2
+\tau^2\|a^{ij}(\varepsilon\cdot)(\partial_j\eta_k)w\|_{L^2(\Omega)}^2\Big),
\end{multline}
where we used that, by the expression of $F_k$,
$\big(A(\varepsilon\cdot)\nabla w_k-F_k\big)\cdot\nu=\eta_k\big(A(\varepsilon\cdot)\nabla w-\widetilde F\big)\cdot\nu$
on $\Gamma$.

Before summing over $k$, we estimate separately the various contributions.

Using the bounds in \eqref{Partition-unity}, we have
\begin{multline}\label{estimate-interior-terms}
\|\nabla\eta_k\cdot\widetilde F\|_{L^2(\Omega)}^2
+\|(\partial_i\eta_k)a^{ij}(\varepsilon\cdot)(\partial_jw)\|_{L^2(\Omega)}^2
+\tau^2\|a^{ij}(\varepsilon\cdot)(\partial_i\Phi)(\partial_j\eta_k)w\|_{L^2(\Omega)}^2
+\tau^2\|a^{ij}(\varepsilon\cdot)(\partial_j\eta_k)w\|_{L^2(\Omega)}^2\\
\leq C\rho^{-2}\left(\|\widetilde F\|_{L^2(B_k\cap\Omega)}^2
+\|\nabla w\|_{L^2(B_k\cap\Omega)}^2
+\tau^2\|w\|_{L^2(B_k\cap\Omega)}^2\right).
\end{multline}

For the boundary terms, we start with the Dirichlet term. Using suitable interpolation estimates, as in \cite[Lemma 7.2]{Dehman-Ervedoza-Thabouti-2023}, one proves
\begin{equation}\label{Interpolation-partition}
\sum_k\|\eta_kw\|_{H^{\frac{1}{2}}(\Gamma)}^2
\leq C\Big(\|w\|_{H^{\frac{1}{2}}(\Gamma)}^2
+\rho^{-1}\|w\|_{L^2(\Gamma)}^2\Big),
\end{equation}
and, since $\sum_k\eta_k^2\leq1$ on the support of $w$,
\begin{equation}\label{Sum-boundary-L2}
\sum_k\|\eta_kw\|_{L^{2}(\Gamma)}^2
\leq\|w\|_{L^{2}(\Gamma)}^2.
\end{equation}
By \eqref{Bound-rho}, $\tau^{\frac{3}{2}}\rho^{-1}\leq2\tau^{2}\leq2\tau^{\frac52}$, so that the $L^2$-terms of \eqref{Interpolation-partition} and \eqref{Sum-boundary-L2} are bounded by $C\tau^{\frac52}\|w\|^2_{L^2(\Gamma)}$, which we keep in the right-hand side. 

On the other hand, for the Neumann term, we set $g=\big(A(\varepsilon\cdot)\nabla w-\widetilde F\big)\cdot\nu$. For $\theta=0$, since $\sum_k\eta_k^2\leq N$,
\begin{equation*}
\sum_k\|\eta_k\,g\|_{L^2(\Gamma)}^2\leq N\|g\|_{L^2(\Gamma)}^2 .
\end{equation*}
For $\theta=-\frac{1}{2}$, we prove
\begin{equation}\label{Gluing-Neumann}
\sum_k\|\eta_k\,g\|^2_{H^{-\frac{1}{2}}(\Gamma)}\leq C\,\|g\|^2_{H^{-\frac{1}{2}}(\Gamma)} .
\end{equation}
Arguing by duality, it suffices to prove that, for all $(v_k)_{0\leq k\leq k_1}$ in $H^{\frac{1}{2}}(\R^{d-1})$,
\begin{equation}\label{Sum-H12}
\Big\|\sum_k \eta_k v_k\Big\|^2_{H^{\frac{1}{2}}(\R^{d-1})}\leq C\sum_k\|v_k\|^2_{H^{\frac{1}{2}}(\R^{d-1})}.
\end{equation}
Indeed, choosing $v_k$ with $\|v_k\|_{H^{\frac{1}{2}}(\R^{d-1})}=\|\eta_kg\|_{H^{-\frac{1}{2}}(\Gamma)}$ and $\langle \eta_kg,v_k\rangle=\|\eta_kg\|^2_{H^{-\frac{1}{2}}(\Gamma)}$, \eqref{Sum-H12} gives
\begin{equation*}
\sum_k\|\eta_kg\|^2_{H^{-\frac{1}{2}}(\Gamma)}=\Big\langle g,\sum_k\eta_kv_k\Big\rangle\leq C^{\frac{1}{2}}\|g\|_{H^{-\frac{1}{2}}(\Gamma)}\Big(\sum_k\|\eta_kg\|^2_{H^{-\frac{1}{2}}(\Gamma)}\Big)^{\frac{1}{2}}.
\end{equation*}
To prove \eqref{Sum-H12}, we set $V_k=\eta_kv_k$ and use the Gagliardo seminorm on $\R^{d-1}$. For $|x'-y'|<\rho$, at most $N'$ indices $k$, with $N'$ depending only on $d$, satisfy $V_k(x')\neq0$ or $V_k(y')\neq0$, so that $\big|\sum_k(V_k(x')-V_k(y'))\big|^2\leq N'\sum_k|V_k(x')-V_k(y')|^2$, while, for $|x'-y'|\geq\rho$, we use $\big|\sum_kV_k(x')\big|^2\leq N\sum_k|V_k(x')|^2$ and $\int_{|x'-y'|\geq\rho}|x'-y'|^{-d}\,dy'=C\rho^{-1}$. Hence
\begin{equation*}
\Big\|\sum_k V_k\Big\|^2_{H^{\frac{1}{2}}(\R^{d-1})}\leq C\sum_k\Big(\|V_k\|^2_{H^{\frac{1}{2}}(\R^{d-1})}+\rho^{-1}\|V_k\|^2_{L^2(\R^{d-1})}\Big)\leq C\sum_k\|V_k\|^2_{H^{\frac{1}{2}}(\R^{d-1})},
\end{equation*}
where the last inequality follows, since $V_k$ is supported in a ball of radius $\rho$, from Hölder's inequality and the Sobolev embedding $H^{\frac{1}{2}}(\R^{d-1})\subset L^{\frac{2(d-1)}{d-2}}(\R^{d-1})$. Finally, Lemma \ref{Lemma-est-product-of-function-in-H-1/2-and--1/2} gives $\|V_k\|_{H^{\frac{1}{2}}(\R^{d-1})}\leq C\|v_k\|_{H^{\frac{1}{2}}(\R^{d-1})}$, which proves \eqref{Sum-H12}.
\medskip

Finally, we sum \eqref{Carleman-Local-A-L2-wk} over $k$. Since $w=\sum_kw_k$, the bounded overlap property gives $\|w\|^2_{L^2(\Omega)}\leq N\sum_k\|w_k\|^2_{L^2(\Omega)}$ and $\|\nabla w\|^2_{L^2(\Omega)}\leq N\sum_k\|\nabla w_k\|^2_{L^2(\Omega)}$. Using moreover $\sum_k\eta_k^2\leq N$, \eqref{Bound-rho}, \eqref{estimate-interior-terms}, \eqref{Interpolation-partition}, \eqref{Sum-boundary-L2} and \eqref{Gluing-Neumann}, we obtain, for all $\tau\geq\max\{\tau_0,16r_0^{-2}\}$,
\begin{multline*}
\tau^{3}\|w\|_{L^2(\Omega)}^2+\tau\|\nabla w\|_{L^2(\Omega)}^2\leq C\Big(\|\widetilde f_2\|_{L^2(\Omega)}^2+ \tau^{\frac{3}{2}-\frac{1}{d}}\|\widetilde f_{2_*'}\|_{L^{\frac{2d}{d+2}}(\Omega)}^2+\tau^2\|\widetilde F\|_{L^2(\Omega)}^2
+\tau^{\frac{3}{2}}\|w\|_{H^{\frac{1}{2}}(\Gamma)}^2\\
+\tau^{\frac{5}{2}}\|w\|_{L^2(\Gamma)}^2
+\tau^{\frac{1}{2}-2 \theta}\big\|\big(A(\varepsilon \cdot) \nabla w-\widetilde F\big)\cdot\nu\big\|_{H^{\theta}(\Gamma)}^2\Big)
+C\big(\varepsilon+r_0^{-2}\tau^{-1}\big)\big(\tau^3\|w\|_{L^2(\Omega)}^2+\tau\|\nabla w\|_{L^2(\Omega)}^2\big).
\end{multline*}
Choosing $\varepsilon\leq\varepsilon_2:=\min\{(4C)^{-1},\varepsilon_1\}$ and $\tau\geq\tau_2:=\max\{\tau_0,(16+4C)r_0^{-2}\}$, where $\varepsilon_1$ and $\tau_0$ are the
constants of Proposition \ref{Prop-Local-Carl-x-0-variable}, the last term is absorbed into the
left-hand side. Since $e^{\tau\Phi}\nabla u=\nabla w-\tau w\nabla\Phi$, $\|\widetilde f_2\|_{L^2(\Omega)}\leq\|e^{\tau\Phi}f_2\|_{L^2(\Omega)}+C\tau\|e^{\tau\Phi}F\|_{L^2(\Omega)}$, $w_{|\Gamma}=(e^{\tau\Phi}u)_{|\Gamma}$ and $\big(A(\varepsilon\cdot)\nabla w-\widetilde F\big)\cdot\nu=\big(A(\varepsilon\cdot)\nabla(e^{\tau\Phi}u)-e^{\tau\Phi}F\big)\cdot\nu$, this yields \eqref{Carleman-General-1} for all $\tau\geq\tau_2$ and all
$\varepsilon\leq\varepsilon_2$. The estimate \eqref{Carleman-General-2} is obtained in the same way from
\eqref{Carleman-Local-A-L2star}, using moreover $\|w\|^2_{L^{\frac{2d}{d-2}}(\Omega)}\leq C\sum_k\|w_k\|^2_{L^{\frac{2d}{d-2}}(\Omega)}$, $\sum_k\|\eta_k\widetilde f_{2_*'}\|^2_{L^{\frac{2d}{d+2}}(\Omega)}\leq C\|\widetilde f_{2_*'}\|^2_{L^{\frac{2d}{d+2}}(\Omega)}$ (as $\frac{2d}{d+2}\leq2$) and $\tau^{\frac{3}{2}+\frac1d}\rho^{-1}\leq2\tau^{\frac52}$, which completes the proof of Theorem \ref{Thm-Carleman-General}.
\section{Proof of Theorem \ref{Thm2-Carleman-General}}\label{Sec-proof-Thm2} 
Under the reverse sign condition $0 < \alpha_1 \leq \varphi' \leq \alpha_2$, $\varphi''\geq \beta$, the only new ingredient is the following analogue of Proposition \eqref{Prop-Carleman-Strip-2}. 
\begin{proposition}\label{Prop-Carleman-Strip-reverse}
    Let $\mathcal{S}$ be as in \eqref{Def-Strip} with $0<L\leq 1,$ let $\varphi \in C^2([0,L])$ satisfy $\varphi(0), \, 0 < \alpha_1 \leq  \varphi' \leq \alpha_2 $ and $\varphi''\geq \beta$ on $[0,L]$, and let $\ell \in \R^{d-1}$ satisfy  \eqref{Def-ell-star}. There exist $C>0$ and $\tau_0 \geq 1$, depending only on $d, \alpha_1, \alpha_2, \beta$ (uniform in $L\leq 1$ and in $|\ell | \leq \ell_*$), such that for all $(f_2,f_{2_*'},F,g) $ as in \eqref{Reg-Source-Terms-Strip}, if $w$ solves 
\begin{equation}\label{Elliptic-Strip-w-reverse}
\left\{
\begin{array}{ll}
e^{\tau(\varphi+\ell\cdot x')}\,\Delta\big(e^{-\tau(\varphi+\ell\cdot x')}w\big)=f_2+f_{2_*'}+\div F & \text{in }\mathcal S,\\
w(0,x')=g(x'), & x'\in\R^{d-1},\\
w(L,x')=0,\qquad(\partial_1w-F_1)(L,x')=0, & x'\in\R^{d-1},
\end{array}
\right.
\end{equation}
then, for all $\tau \geq \tau_0$, 
\begin{equation}\label{Carleman-Strip-reverse-1}
\tau^{\frac{3}{2}}\|w\|_{L^2(\mathcal S)}+\tau^{\frac{1}{2}}\|\nabla w\|_{L^2(\mathcal S)}
\leq C\Big(\|f_2\|_{L^2(\mathcal S)}+\tau^{\frac{3}{4}-\frac1{2d}}\|f_{2_*'}\|_{L^{\frac{2d}{d+2}}(\mathcal S)}+\tau\|F\|_{L^2(\mathcal S)}+\tau^{\frac{3}{4}}\|g\|_{H^{\frac{1}{2}}(\R^{d-1})}\Big),
\end{equation}
and
\begin{equation}\label{Carleman-Strip-reverse-2}
\tau^{\frac{3}{4}+\frac1{2d}}\|w\|_{L^{\frac{2d}{d-2}}(\mathcal S)}
\leq C\Big(\|f_2\|_{L^2(\mathcal S)}+\tau^{\frac{3}{4}+\frac1{2d}}\|f_{2_*'}\|_{L^{\frac{2d}{d+2}}(\mathcal S)}+\tau\|F\|_{L^2(\mathcal S)}+\tau^{\frac{3}{4}+\frac1{2d}}\|g\|_{H^{\frac{1}{2}}(\R^{d-1})}\Big).
\end{equation}
\end{proposition}
\begin{remark}\label{Remark-parametrix-for the second Carleman theorem}
    With the same notation of Subsection \ref{Sec-Construction-parametrix}, the construction of the parametrix for Theorem  \ref{Thm2-Carleman-General} is obtained as follows. We write the factorization \eqref{factorisation} in the reverse order,   
    \begin{equation*}
        (\partial_1-\tau\varphi')^2-\psi^2
=\big(\partial_1-\tau\varphi'-\psi\big)\big(\partial_1-\tau\varphi'+\psi\big),
    \end{equation*}
    and we set $\widehat z = \big(\partial_1-\tau\varphi'+\psi\big)\widehat w.$ Since $\tau\varphi'(x_1)+\operatorname{Re}\psi(\xi') \geq\tau\alpha_1>0$, the equation for $\widehat z-\widehat F_1$, associated with the first factor, is solved by integrating backward from $ \{ x_1=L\}$. This uses the boundary term $(\widehat z-\widehat F_1)(L,\xi')=(\partial_1\widehat w-\widehat F_1)(L,\xi')=0$. We then solve $ \big(\partial_1-\tau\varphi'+\psi\big)\widehat w
=\widehat z,$ choosing the direction of integration according to the sign of $\tau\varphi'(x_1)-\operatorname{Re}\psi(\xi')$, which is increasing in $x_1$ since $\varphi''\geq \beta >0$: we integrate forward from $x_1=0$, using the Dirichlet trace $\widehat w(0,\xi')$, where $\operatorname{Re}\psi(\xi')>\tau\varphi'(x_1)$, and backward from $x_1=L$, using $\widehat w(L,\xi')=0$, in the complementary region. The resulting formula expresses $w$ in terms of the source terms and its Dirichlet trace on $\{x_1=0\}$, without any Neumann trace; it is  \eqref{Parametrix-reverse} below.
\end{remark}
\begin{proof}[Proof of Proposition \ref{Prop-Carleman-Strip-reverse}]
    Let $\mathcal T(x_1,x')=(L-x_1,x')$, and set $\widetilde w=w\circ\mathcal T$, $\widetilde\varphi(x_1)=\varphi(L-x_1)-\varphi(L)$, $\widetilde f_2=f_2\circ\mathcal T$, $\widetilde f_{2_*'}=f_{2_*'}\circ\mathcal T$ and $\widetilde F=(-F_1,F_2,\dots,F_d)\circ\mathcal T$. Since $\Delta$ commutes with $\mathcal T$ and an additive constant in the weight does not change the conjugated operator, $\widetilde w$ solves
\begin{equation*}
\left\{
\begin{array}{ll}
e^{\tau(\widetilde\varphi+\ell\cdot x')}\,\Delta\big(e^{-\tau(\widetilde\varphi+\ell\cdot x')}\widetilde w\big)=\widetilde f_2+\widetilde f_{2_*'}+\div\widetilde F & \text{in }\mathcal S,\\
\widetilde w(0,x')=0,\qquad(\partial_1\widetilde w-\widetilde F_1)(0,x')=0, & x'\in\R^{d-1},\\
\widetilde w(L,x')=g(x'), & x'\in\R^{d-1},
\end{array}
\right.
\end{equation*}
and $\widetilde\varphi$ satisfies \eqref{Weight-function-varphi} on $[0,L]$. We may therefore follow the construction of Proposition \ref{Prop-parametrix-Strip-1} for $\widetilde w$. Since the boundary terms of $\widetilde w$ on $\{x_1=0\}$ vanish, the boundary terms of \eqref{Explicit-z} and \eqref{w-forward} vanish as well, while the backward integration \eqref{w-backward} now uses $\widehat{\widetilde w}(L,\xi')=\widehat g(\xi')$, which adds to its right-hand side the term
\begin{equation*}
\mathbf 1_{\operatorname{Re}\psi(\xi')\geq-\tau\widetilde\varphi'(x_1)}\,e^{\tau(\widetilde\varphi(x_1)-\widetilde\varphi(L))+(x_1-L)\psi(\xi')}\,\widehat g(\xi').
\end{equation*}
Replacing $x_1$ by $L-x_1$, we obtain
\begin{equation}\label{Parametrix-reverse}
w=\Big(K_{\tau,0}\big(\widetilde f_2+\widetilde f_{2_*'}\big)+\sum_{j=1}^dK_{\tau,j}(\widetilde F_j)\Big)\circ\mathcal T+G^\sharp_\tau g,
\end{equation}
where $K_{\tau,j}$ are the operators of Proposition \ref{Prop-parametrix-Strip-1} associated with the weight $\widetilde\varphi+\ell\cdot x'$, and $G^\sharp_\tau$ is given by
\begin{equation}\label{Def-G-sharp}
\widehat{G^\sharp_\tau g}(x_1,\xi')=g^\sharp_\tau(x_1,\xi')\,\widehat g(\xi'),\qquad
g^\sharp_\tau(x_1,\xi')=\mathbf 1_{\operatorname{Re}\psi(\xi')\geq\tau\varphi'(x_1)}\,e^{\tau\varphi(x_1)-x_1\psi(\xi')} .
\end{equation}
Since the norms are invariant under $\mathcal T$, the operators $K_{\tau,j}$ are estimated by Proposition \ref{Prop-low-frequencies} and Proposition \ref{Prop-K-high}, exactly as in the proof of Proposition \ref{Prop-Carleman-Strip-2}, and $G^\sharp_\tau$ by Lemma \ref{Lemma-G-sharp} below, multiplied by $\tau^{\frac{3}{4}}$ and $\tau^{\frac{3}{4}+\frac1{2d}}$ respectively. This gives \eqref{Carleman-Strip-reverse-1} and \eqref{Carleman-Strip-reverse-2}.
\end{proof}
\begin{lemma}\label{Lemma-G-sharp}
    Under the assumptions of Proposition \ref{Prop-Carleman-Strip-reverse}, there exists $C>0$, depending only on $d, \alpha_1,\alpha_2,\beta$, such that, for all $\tau\geq 1$ and all $g\in H^{\frac{1}{2}}(\R^{d-1})$, 
\begin{equation}\label{Est-G-sharp}
\tau^{\frac{3}{4}}\big\|G^\sharp_\tau g\big\|_{L^2(\mathcal S)}+\tau^{-\frac14}\big\|\nabla G^\sharp_\tau g\big\|_{L^2(\mathcal S)}+\big\|G^\sharp_\tau g\big\|_{L^{\frac{2d}{d-2}}(\mathcal S)}\leq C\|g\|_{H^{\frac{1}{2}}(\R^{d-1})}.
\end{equation}
\end{lemma}

\begin{proof}
Set $\lambda=\operatorname{Re}\psi(\xi')$. Since $\varphi'$ is increasing, $g^\sharp_\tau(x_1,\xi')\neq0$ implies $\lambda\geq\tau\varphi'(x_1)\geq\alpha_1\tau\geq4\tau|\ell|$, so that $G^\sharp_\tau=G^\sharp_\tau\Pi_{\mathrm{high}}$ and, since $|\xi'|^2=\lambda^2-(\operatorname{Im}\psi)^2+\tau^2|\ell|^2$ with $|\operatorname{Im}\psi|\leq\tau|\ell|$, $\lambda\leq|\xi'|\leq\frac{3}{2}\lambda$. Moreover, since $\tau\varphi'(x_1)-\tau\varphi'(r)\geq\beta\tau(x_1-r)$ for $0\leq r\leq x_1$, we have, on the support of $g^\sharp_\tau$,
\begin{equation}\label{Bound-g-sharp}
\big|g^\sharp_\tau(x_1,\xi')\big|=e^{-\int_0^{x_1}(\lambda-\tau\varphi'(r))\,dr}\leq e^{-|\lambda-\tau\varphi'(x_1)|\,x_1-\frac\beta2\tau x_1^2},
\end{equation}
which is the bound on $g_\tau$ used in the proof of \eqref{Est-G-tau-high}, with $-\tau\varphi'(x_1)$ replaced by $\tau\varphi'(x_1)$. In particular, $\int_0^L|g^\sharp_\tau(x_1,\xi')|^2dx_1\leq C\tau^{-\frac{1}{2}}$, and $\int_0^L|g^\sharp_\tau(x_1,\xi')|^2dx_1\leq\lambda^{-1}$ if $\lambda\geq2\alpha_2\tau$, since then $\lambda-\tau\varphi'(x_1)\geq\frac\lambda2$. By Plancherel's theorem, and since $\alpha_1\tau\leq\lambda\leq|\xi'|$ on the support of $g^\sharp_\tau$,
\begin{equation*}
\tau^{\frac{3}{2}}\big\|G^\sharp_\tau g\big\|^2_{L^2(\mathcal S)}\leq C\tau\int_{\lambda\geq\alpha_1\tau}|\widehat g(\xi')|^2\,d\xi'\leq C\int_{\R^{d-1}}|\xi'|\,|\widehat g(\xi')|^2\,d\xi'.
\end{equation*}
For the gradient, on the support of $g^\sharp_\tau$ and away from the turning point $\lambda=\tau\varphi'(x_1)$, as for the other operators of the parametrix, $|\partial_1g^\sharp_\tau|=|\tau\varphi'(x_1)-\psi(\xi')|\,|g^\sharp_\tau|\leq C\lambda\,|g^\sharp_\tau|$ and $|\xi'|\leq\frac{3}{2}\lambda$, so that
\begin{equation*}
\big\|\nabla G^\sharp_\tau g\big\|^2_{L^2(\mathcal S)}\leq C\int_{\R^{d-1}}\lambda^2\Big(\tau^{-\frac{1}{2}}\mathbf 1_{\lambda<2\alpha_2\tau}+\lambda^{-1}\mathbf 1_{\lambda\geq2\alpha_2\tau}\Big)|\widehat g(\xi')|^2\,d\xi'\leq C\tau^{\frac{1}{2}}\int_{\R^{d-1}}|\xi'|\,|\widehat g(\xi')|^2\,d\xi'.
\end{equation*}
Finally, the $H^{\frac{1}{2}}-L^{\frac{2d}{d-2}}$ estimate follows from the $TT^\star$ argument of the proof of \eqref{Est-G-tau-high}, which only uses the bound on the multiplier, here replaced by \eqref{Bound-g-sharp}: the $\lambda$-integral is split in the same way, the region $\lambda<\alpha_1\tau$ not contributing, and $|\lambda-\tau\varphi'(x_1)|$ being possibly zero on $[\alpha_1\tau,2\alpha_2\tau]$ and bounded below by $\frac\lambda2$ on $(2\alpha_2\tau,+\infty)$.
\end{proof}
\begin{proof}[Proof of Theorem \ref{Thm2-Carleman-General}.] The local estimates, analogous to Theorem \ref{Thm-Local-Carl-x-0}, their variable-coefficient version, analogous to Proposition \ref{Prop-Local-Carl-x-0-variable}, and the gluing are obtained exactly as in Subsection \ref{Sec-Proof-Carl-x-0} and Section \ref{Sec-Carleman-variable}, using Proposition \ref{Prop-Carleman-Strip-reverse} instead of Proposition \ref{Prop-Carleman-Strip-2}, with the following changes.
\begin{itemize}
\item There are no Neumann terms and no $L^2$ boundary terms in the local estimates.
\item Proposition \ref{Prop-Carleman-Strip-reverse} is applied to functions which, together with the source terms, vanish near $\{x_1=2r_0\}$, so that the last line of \eqref{Elliptic-Strip-w-reverse} holds with $L=2r_0$.
\item In Step 1 of the proof of Proposition \ref{Prop-Local-Carl-x-0-variable}, one has $\frac{\alpha_1}2\leq\varphi_0'=\varphi'+\ell_0\cdot B\leq2\alpha_2$ for $\varepsilon$ small enough.
\item In the gluing, the last term of \eqref{Interpolation-partition}, multiplied by $\tau^{\frac{3}{2}}$, is absorbed into the left-hand side thanks to the trace inequality $\|w\|^2_{L^2(\Gamma)}\leq2\|w\|_{L^2(\Omega)}\|\nabla w\|_{L^2(\Omega)}$ and \eqref{Bound-rho}, which give
\begin{equation*}
\tau^{\frac{3}{2}}\rho^{-1}\|w\|^2_{L^2(\Gamma)}\leq2\big(\varepsilon^{\frac{1}{2}}+2r_0^{-1}\tau^{-\frac{1}{2}}\big)\,\tau^{\frac{3}{2}}\|w\|_{L^2(\Omega)}\,\tau^{\frac{1}{2}}\|\nabla w\|_{L^2(\Omega)}.
\end{equation*}
This gives \eqref{Carleman-General-1-Thm2}. For \eqref{Carleman-General-2-Thm2}, the corresponding term $\tau^{\frac{3}{2}+\frac1d}\rho^{-1}\|w\|^2_{L^2(\Gamma)}$ cannot be absorbed in this way, since the trace inequality loses a factor $\tau^{\frac1d}$; it is bounded by $C\tau^{\frac52}\|w\|^2_{L^2(\Gamma)}$, which is why \eqref{Carleman-General-2-Thm2} contains the term $\tau^{\frac54}\|e^{\tau\Phi}u\|_{L^2(\Gamma)}$.
\end{itemize}
\end{proof}

\appendix

\section{Proof of Lemma \ref{Lemma-est-product-of-function-in-H-1/2-and--1/2}} \label{Proof-Lemma-est-product-of-function-in-H-1/2-and--1/2}
\begin{proof}[Proof of Lemma \ref{Lemma-est-product-of-function-in-H-1/2-and--1/2}]
We first recall that the norm $H^{1/2}(\R^n)$ defined by $\|f\|_{H^{1/2}(\R^n)}= \|(1+|\xi|^2)^{\frac{1}{4}}  \widehat f  \|_{L^2(\R^n)}$  is equivalent to 
$$ \|f\|_{L^2(\R^n)}
+
[f]_{\dot H^{\frac{1}{2}}(\R^n)},\quad \text{where } \quad
[f]^2_{\dot H^{\frac{1}{2}}(\R^n)}=
\iint_{\R^n\times\R^n}
\frac{|f(z)-f(z')|^2}{|z-z'|^{n+1}}
\,dz\,dz',$$
We shall therefore estimate separately the $L^2$ norm and the
$\dot H^{1/2}$ seminorm.

First, for all $R>0$ and all $f\in\dot H^{\frac{1}{2}}(\R^n)$, Hölder's inequality and by the Sobolev embedding $\dot H^{\frac{1}{2}}(\R^{n})\hookrightarrow L^{\frac{2n}{n-1}}(\R^{n})$,
\begin{equation}\label{local-Sobolev}
\|f\|_{L^2(B(x_0',R))}\leq|B(x_0',R)|^{\frac{1}{2n}}\,\|f\|_{L^{\frac{2n}{n-1}}(\R^{n})}
\leq C\,R^{\frac{1}{2}}\,\|f\|_{\dot H^{\frac{1}{2}}(\R^{n})}.
\end{equation}
We next estimate the $\dot H^{1/2}$ seminorm of $\widetilde m f$. To do that, we first estimate  
\begin{equation*}
   [\widetilde m f]_{\dot H^{\frac{1}{2}}(\R^{n})}^2
\leq
2\|\widetilde m\|_{L^\infty}^2
[f]_{\dot H^{\frac{1}{2}(\R^{n})}}^2
+
2\int_{\R^n} \int_{\R^n}
\frac{
|\widetilde m(z)-\widetilde m(z')|^2
}{
|z-z'|^{n+1}
}
\,dz |f(z')|^2\,dz'. 
\end{equation*}
We then claim that for all $z'\in\R^n$
\begin{equation} \label{estimate-integral-m-1}
    \int_{\R^n}
\frac{
|\widetilde m(z)-\widetilde m(z')|^2
}{
|z-z'|^{n+1}
}
\,dz\leq C M^2/ \rho.
\end{equation}
To prove this, we use the fact that 
$|\widetilde m(z)-\widetilde m(z')| \leq \min\left\{
2M,\,
\frac{M}{\rho}|z-z'|
\right\},$
and splitting the integral 
\begin{equation*}
  \int_{\R^n}
\frac{
|\widetilde m(z)-\widetilde m(z')|^2
}{
|z-z'|^{n+1}
}
\,dz
\leq
C M^2
\int_{|z-z'|\leq\rho}
\frac{|z-z'|^2 }{|z-z'|^{n+1} \rho^2}
\,dz
+
C M^2
\int_{|z-z'|>\rho}
\frac{1}{|z-z'|^{n+1}}
\,dz. 
\end{equation*}
Then, using polar coordinates, both integrals are bounded by $CM^2/\rho$. 

On the other hand, if $|z'-x_0'|\geq4\rho$, then
$\widetilde m(z')=0$ and $|z-z'|\geq\frac{1}{2}|z'-x_0'|$ on $\operatorname{supp}\widetilde m$, so that
\begin{equation} \label{estimate-integral-m-2} 
  \int_{\R^n}
\frac{
|\widetilde m(z)-\widetilde m(z')|^2
}{
|z-z'|^{n+1}
}
\,dz \leq C M^2\rho^{n}|z'-x_0'|^{-n-1}. 
\end{equation}

Consequently, decomposing $\R^n$ into $B(x_0',4\rho)$ and the
dyadic annuli $B(x_0',2^{j+1}\rho)\setminus B(x_0',2^{j}\rho)$, $j\geq2$, and using
\eqref{local-Sobolev}, \eqref{estimate-integral-m-1} and \eqref{estimate-integral-m-2} 
\begin{align*}
\int_{\R^{n}}  \int_{\R^n}
\frac{
|\widetilde m(z)-\widetilde m(z')|^2
}{
|z-z'|^{n+1}
}
\,dz \,|f(z')|^2\,dz'
&\leq\frac{CM^2}{\rho}\|f\|^2_{L^2(B(x_0',4\rho))}
+CM^2\rho^{n}\sum_{j\geq2}\big(2^{j}\rho\big)^{-n-1}\|f\|^2_{L^2(B(x_0',2^{j+1}\rho))}
\\& \leq CM^2\Big(1+\sum_{j\geq2}2^{-jn}\Big)\|f\|^2_{\dot H^{\frac{1}{2}}(\R^n)} .
\end{align*}
Together with $\|\widetilde m\,h\|_{L^2}\leq M\|h\|_{L^2}$, this gives the first estimate of
\eqref{Est-product-scale}. The second follows by duality, which ends the poof the lemma. 
\end{proof}

\noindent{\bf Acknowledgments.} The authors warmly thank Sylvain Ervedoza for fruitful discussions on a preliminary version of this paper.

\medskip
\noindent{\bf Use of AI tools.} The results, proofs and writing of this paper are the authors' own. The authors used Claude (Anthropic) and ChatGPT (OpenAI) for detecting typos, editing the language, improving the exposition, helping draw Figures \ref{Geometry-local-at-the-boundary}, \ref{Geometry near the boundary-interior observation}, \ref{Geometry far from the boundary-interior observation}, and \ref{Three-ball-step}, and in identifying relevant references. All the content produced with these tools, including the references, was checked by the authors, who take the full responsibility for the content of this paper.  

\bibliographystyle{plain}

\end{document}